\documentclass[11pt]{amsart}
\usepackage[T1]{fontenc}
\usepackage[utf8]{inputenc}
\usepackage{amssymb}
\usepackage{comment}
\usepackage{thmtools,xcolor}
\usepackage{thm-restate}
\usepackage{graphicx}
\usepackage[all]{xy}
\usepackage{tikz}
\usetikzlibrary{arrows,calc,positioning,decorations.pathreplacing}
\usepackage{mathtools}
\usepackage{geometry} 
\usepackage[mathscr]{euscript}
\usepackage{accents}
\usepackage{combelow}

\numberwithin{equation}{section}
\numberwithin{figure}{section}

\declaretheorem[name=Theorem,numberwithin=section]{thm}
\declaretheorem[name=Lemma,numberlike=thm]{lem}
\declaretheorem[name=Proposition,numberlike=thm]{prop}
\declaretheorem[name=Corollary,numberlike=thm]{coro}

\theoremstyle{definition}
\declaretheorem[name=Definition,numbered=yes,numberlike=thm]{defn}
\declaretheorem[name=Remark,numbered=yes,numberlike=thm]{rmk}
\declaretheorem[name=Notation,numbered=yes,numberlike=thm]{notation}

\usepackage[colorlinks,citecolor=blue,linkcolor=blue,urlcolor=blue,filecolor=blue,breaklinks]{hyperref}

\hypersetup{
    pdftitle={Quantum braid Floer homology},
    pdfauthor={Cristina Anghel, András Juhász},
    pdfkeywords={knot theory, braid, categorification, Floer homology, topological model, quantum invariant, Alexander polynomial, Jones polynomial}
}

\newcommand{\cB}{\mathcal{B}}
\newcommand{\cC}{\mathcal{C}}
\newcommand{\cD}{\mathcal{D}}
\newcommand{\cF}{\mathcal{F}}
\newcommand{\cH}{\mathcal{H}}
\newcommand{\cI}{\mathcal{I}}
\newcommand{\cJ}{\mathscr{J}}
\newcommand{\cK}{\mathscr{K}}
\newcommand{\cM}{\mathcal{M}}

\newcommand{\cs}{\mathscr{S}}
\newcommand{\ct}{\mathscr{T}}

\newcommand{\bD}{\mathbb{D}}
\newcommand{\bF}{\mathbb{F}}
\newcommand{\bI}{\mathbb{I}}
\newcommand{\bN}{\mathbb{N}}
\newcommand{\bQ}{\mathbb{Q}}
\newcommand{\bR}{\mathbb{R}}
\newcommand{\bT}{\mathbb{T}}
\newcommand{\bZ}{\mathbb{Z}}
\newcommand{\ds}{\boldsymbol{d}}
\newcommand{\mK}{\overline{K}}
\newcommand{\ph}{\widehat{\partial}}
\newcommand{\sI}{\mathscr{I}}
\newcommand{\fb}{\mathfrak{b}}
\newcommand{\fc}{\mathfrak{c}}
\newcommand{\fq}{\mathfrak{q}}

\DeclareMathOperator{\Id}{Id}
\DeclareMathOperator{\rank}{rank}
\DeclareMathOperator{\gr}{gr}
\DeclareMathOperator{\Sym}{Sym}

\DeclareMathOperator{\Conf}{Conf}
\DeclareMathOperator{\qHF}{qHF}

\DeclareMathOperator{\qHFB}{qHFB}
\DeclareMathOperator{\HFB}{HFB}
\DeclareMathOperator{\qCFB}{qCFB}
\DeclareMathOperator{\CFK}{CFK}
\newcommand{\cMh}{\widehat{\cM}}
\newcommand{\CFKh}{\widehat{\CFK}}
\DeclareMathOperator{\CF}{CF}
\newcommand{\CFh}{\widehat{\CF}}
\DeclareMathOperator{\HF}{HF}
\newcommand{\HFh}{\widehat{\HF}}
\DeclareMathOperator{\HFK}{HFK}
\newcommand{\HFKh}{\widehat{\HFK}}
\DeclareMathOperator{\sgn}{\text{sgn}}
\DeclareMathOperator{\Pd}{Pd}
\DeclareMathOperator{\qPd}{qPd}

\newcommand{\as}{\boldsymbol{\alpha}}
\newcommand{\bs}{\boldsymbol{\beta}}
\newcommand{\one}{\boldsymbol{1}}
\newcommand{\ba}{\boldsymbol{a}}
\newcommand{\bb}{\boldsymbol{b}}
\newcommand{\hs}{\boldsymbol{h}}
\newcommand{\js}{\boldsymbol{j}}
\newcommand{\ks}{\boldsymbol{k}}
\newcommand{\qs}{\boldsymbol{q}}
\newcommand{\us}{\boldsymbol{u}}
\newcommand{\ws}{\boldsymbol{w}}
\newcommand{\xs}{\boldsymbol{x}}
\newcommand{\ys}{\boldsymbol{y}}
\newcommand{\zs}{\boldsymbol{z}}
\definecolor{dgreen}{RGB}{0,150,0}

\begin{document}

\title{Quantum braid Floer homology}

\author{Cristina Anghel}
\address{Laboratoire de math\'ematiques Blaise Pascal, Universit\'e Clermont-Auvergne, Campus des C\'ezeaux 3, place Vasarely, 63178 Aubi\`ere, France; Institute of Mathematics “Simion Stoilow” of the Romanian Academy, 21 Calea Grivi\cb{t}ei Street, 010702 Bucharest, Romania.}
\email{cristina.anghel@uca.fr; cranghel@imar.ro}

\author{Andr\'{a}s Juh\'{a}sz}
\address{Mathematical Institute, University of Oxford, Andrew Wiles Building,
	Radcliffe Observatory Quarter, Woodstock Road, Oxford, OX2 6GG, UK}
\email{juhasza@maths.ox.ac.uk}

\date{}

\begin{abstract}
    Given a braid whose closure is a knot, we associate to it a doubly pointed Heegaard diagram with additional marked points that we call a quantum Heegaard diagram. This is obtained by surgery from the unified intersection model for the Alexander and Jones polynomials due to the first author. By translating the Alexander and quantum gradings from the disc model to the Heegaard diagram, we obtain a bifiltration on the hat knot Floer complex. The degree-$(0,0)$ component of the differential with respect to this bifiltration defines a triply graded invariant of the braid that we call quantum braid Floer homology. It is generally not invariant under Markov conjugation or stabilisation. It contains a distinguished non-zero element, and admits a spectral sequence to knot Floer homology. Furthermore, its graded Euler characteristic is a two-variable polynomial that specialises to the Alexander and Jones polynomials.
\end{abstract}

\maketitle

\section{Introduction}\label{introduction}
The Alexander and Jones polynomials are important invariants in knot theory that have different natures. The Alexander polynomial is genuinely geometric, arising from knot complements, whereas the geometry of the Jones polynomial is a fundamental open problem in quantum topology. Categorifications of knot polynomials have become powerful tools in knot theory. Khovanov homology was the first categorification of the Jones polynomial, defined via combinatorial tools. Seidel and Smith~\cite{seidel-smith} and Manolescu~\cite{manolescu} have constructed symplectic versions of Khovanov homology via Hilbert schemes. Ozsv{\'a}th and Szab{\'o}~\cite{OS2}, and independently Rasmussen~\cite{Rass}, have defined knot Floer homology, a categorification of the Alexander polynomial using Heegaard diagrams. Rasmussen~\cite{Rass} conjectured that there is a spectral sequence from reduced Khovanov homology of the mirror of a knot to its $\delta$-graded knot Floer homology, which was later constructed by Dowlin~\cite{Dow} with coefficients in $\bQ$. There are still important open questions about geometric categorifications of the Jones polynomial. Aganagi\'c~\cite{Ag, Ag2} proposed such a categorification using immersed curves predicted by mirror symmetry, and this was developed by LePage and Shende~\cite{LePageShende}. In the configuration space model of the first author~\cite{Anghel-cjp, Anghel2023TAMS}, the Jones polynomial is recovered from the intersection of two explicit embedded Lagrangians graded by a local system, each obtained as the image of a product of pairwise disjoint arcs or circles in a punctured disc. It remains open to promote this intersection model to a Lagrangian Floer homology that is invariant under the two Markov moves and whose graded Euler characteristic recovers the Jones polynomial. It would also be interesting to construct a geometric spectral sequence from such a theory to knot Floer homology.

\subsection{Main results} Let $n \ge 1$ and $\beta \in B_{n+1}$ be a braid whose closure is a knot $K$. Using $\beta$, we construct a doubly pointed Heegaard diagram $\cH_\beta = (\Sigma, \as, \bs, w, z)$ with additional marked points $\qs$ and $\hs$, and a distinguished intersection point $\xs^0$ that we call a \emph{quantum Heegaard diagram}. We construct this as a stabilisation of the disc model associated with $\beta$ by the first author~\cite{Anghel2023TAMS, Anghel2024AIF}. We show that $\cH_\beta$ is a Heegaard diagram for $K$. Furthermore, we identify the Alexander grading on the disc model with the Alexander grading $A^{\HF}$ on the knot Floer complex $(\CFKh(\cH_\beta), \ph)$. We define a \emph{quantum Alexander grading} $A^{\qHF}$ on $(\CFKh(\cH_\beta), \ph)$ using $\qs$, $\hs$, and $\xs^0$. Using $A^{\HF}$ and $A^{\qHF}$, we associate a two-variable polynomial $\Omega^q_\beta(x,d)$ to the intersection of the Heegaard tori that recovers both the Alexander and Jones polynomials as specialisations of the variables. This provides an intersection model for the Jones polynomial using quantum Heegaard diagrams.

We also define a bifiltration on $(\CFKh(\cH_\beta), \ph)$ and let $\partial^q$ be the degree-$(0,0)$ part of $\ph$. The homology with respect to $\partial^q$ is a triply graded group $\qHFB(\beta)$ called \emph{quantum braid Floer homology}, whose isomorphism class depends only on $\beta$ as an element of the braid group. It admits a spectral sequence to $\HFKh(K)$. Furthermore, its bigraded Euler characteristic is $\Omega^q_\beta(x,d)$. The element $\xs^0$ is a cycle whose homology class in $\qHFB(\beta)$ is a non-zero braid invariant $\Theta(\beta)$. The polynomial $\Omega^q_\beta(x,d)$ is not invariant under the two Markov moves (see Section~\ref{sec:examples}); hence, neither is $\qHFB(\beta)$.

\subsection{Topological models in the disc} Bigelow~\cite{Bigelow2002} gave an intersection model for the Jones polynomial using immersed Lagrangians. The first author~\cite{Anghel2023TAMS, Anghel2024AIF} provided a topological model that recovers both the Alexander and Jones polynomials as specialisations of a graded intersection between two embedded Lagrangians in a configuration space; see Section~\ref{S:unifmodel}. The geometric supports of these Lagrangians are arcs and simple closed curves in the punctured disc. The curves are obtained by untwisting the figure eights in Bigelow's model. This is a geometric framework that allows us to see both the Alexander and Jones polynomials within a single configuration-space framework. Our construction identifies the Alexander specialisation of the disc model with the graded Euler characteristic of knot Floer homology.
 
\subsection{The Alexander polynomial as a graded intersection in the disc model}\label{ss1} 
\begin{figure}
\centering
\def\svgwidth{0.9\textwidth}
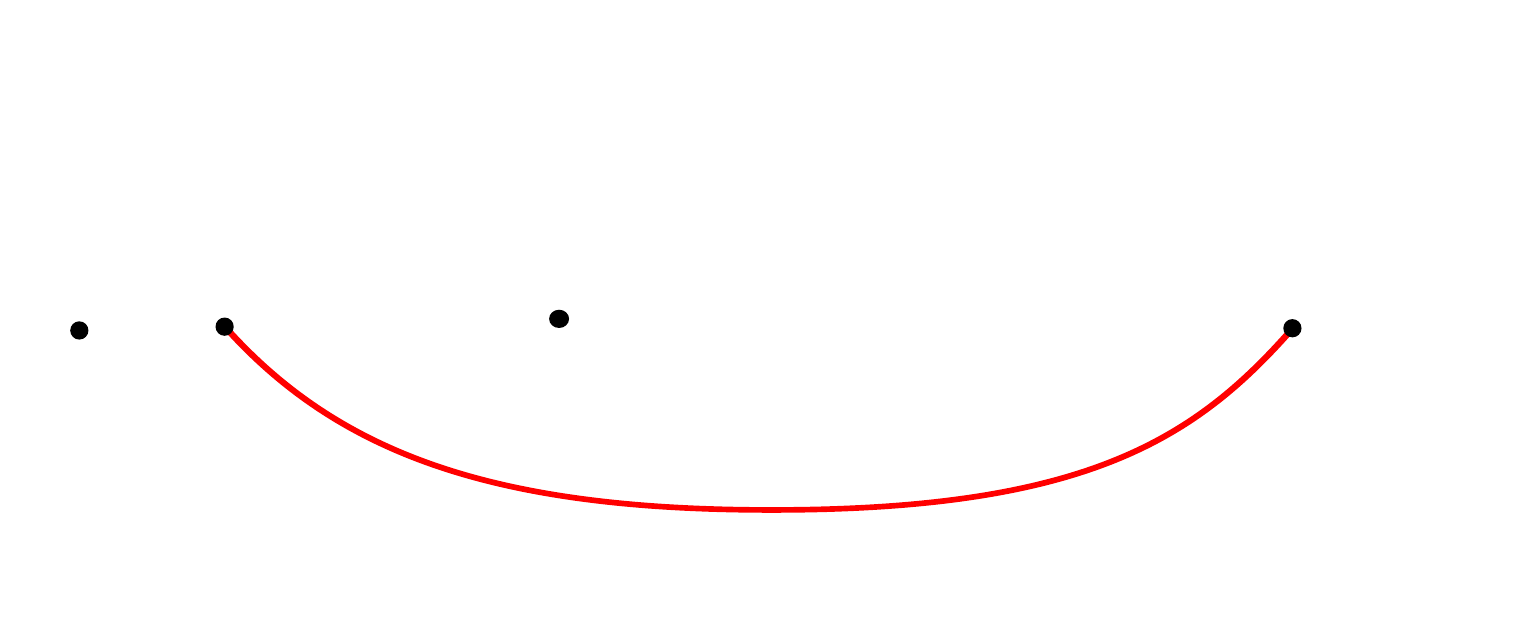
\caption{The arcs $a_1,\dots,a_n$ and the curves $b_1,\dots,b_n$ in the punctured disc $\cD_n$.}
\label{fig:disc-model}
\end{figure}

We review a construction of the first author~\cite{Anghel2024AIF} that exhibits the Alexander polynomial of a knot as a graded intersection in a configuration space. Let $K$ be a knot in $S^3$, obtained as the braid closure of the $(n+1)$-strand braid $\beta$. For $i \in \{0,\dots,n\}$, consider the base points 
\[
z_i := \left(-\frac{n-i+1}{n+2}, 0 \right) \text{ and } w_i := \left(\frac{n-i+1}{n+2}, 0 \right) 
\]
in $D^2$, and write
\[
\zs := \{z_0,\dots,z_n\} \text{ and } \ws := \{w_0,\dots,w_n\}.
\]
We endow the punctured disc 
\[
\cD_n := D^2 \setminus (\zs \cup \ws)
\]
with the restriction of the standard symplectic structure on $\bR^2$. Let $S^1(p, r)$ be the circle in $\bR^2$ with centre $p$ and radius $r$. For $i \in \{0,\dots,n\}$, let 
\begin{equation}\label{eqn:a_i-c_i}
a_i := S^1(\underline{0}, |z_i|) \cap \{y<0\} \text{ and } c_i := S^1(\underline{0}, |z_i|) \cap \{y > 0\}
\end{equation}
be semicircular arcs in $\cD_n$ such that $a_i$ is oriented from $z_i$ to $w_i$, as in Figure~\ref{fig:disc-model}. Furthermore, let $N(c_i)$ be a small regular neighbourhood of $c_i$ such that $N(c_i) \cap (\zs \cup \ws) = \{z_i, w_i\}$. Then
\[
b_i := \partial N(c_i)
\]
is a simple closed curve around punctures $z_i$ and $w_i$, oriented anticlockwise. We write
\[
\ba := \{a_1,\dots,a_n\} \text{ and } \bb := \{b_1,\dots,b_n\}.
\]
Note that no two curves in $\bb$ are nested. Let $C_n$ be the configuration space of $n$ distinct unordered points in $\cD_n$. The product symplectic form on $\cD_n^n$ descends to a symplectic form on $C_n$. We use a local system $\varphi \colon \pi_1(C_n) \to \bZ[x^{\pm 1}]^\times$, which counts monodromies around punctures in $\zs$ anticlockwise, punctures in $\ws$ clockwise, and linking with the fat diagonal with sign $-1$ when $n \ge 2$. See the lower-right panel of Figure~\ref{fig:specialisations} for a schematic picture of $\varphi$. Consider the Lagrangian submanifolds $\cs = a_1 \times \dots \times a_n$ and $\ct := b_1 \times \dots \times b_n$ of $C_n$ (note that the quotient map $\cD_n^n \to C_n$ is injective on these products, and we endow $\cs$ and $\ct$ with the product orientation).

Let $\bI$ be the trivial $(n+1)$-braid. The braid group $B_{2n+2}$ is isomorphic to the mapping class group of $\cD_n$ relative to $\partial \cD_n$. Using this isomorphism, let
\[
\ba' := \{a_1',\dots,a_n'\}
\]
be the action of the braid $\beta \cup \mathbb I \in B_{2n+2}$ on $\ba$ (the inclusion $B_{n+1} \times B_{n+1} \to B_{2n+2}$ is obtained by placing two $(n+1)$-braids side-by-side). We choose a representative of the action of $\beta \cup \bI$ that is supported in the half-plane $\{x<0\}$ and such that $\ba'$ is transverse to $\bb$. The \emph{disc model associated with $\beta$} is
\[
\bD_\beta := (\cD_n, \ba', \bb).
\]
We define the Lagrangian 
\[
(\beta \cup \mathbb I)\cs := a_1' \times \dots \times a_n' \subseteq C_n.
\]
Then we obtain the transverse intersection
\begin{equation}
\cI_{\bD_\beta} := (\beta \cup \mathbb I) \cs \cap \ct.
\end{equation}
Our graded intersection $\langle  (\beta \cup \mathbb I) \cs, \ct \rangle$ is parametrised by $\cI_{\bD_\beta}$ and graded by the local system $\varphi \colon \pi_1(C_n) \to \bZ[x^{\pm 1}]^\times$ introduced in Definition~\ref{def:varphi}, as follows. To each point $\xs \in \cI_{\bD_\beta}$, we associate a loop $l_{\xs}$ in $C_n$ in Definition~\ref{def:loop}, which is graded by the local system $\varphi$:
\[
\xs \in \cI_{\bD_\beta} \rightsquigarrow \text{loop } l_{\xs} \rightsquigarrow \text{grading }  \varphi(l_{\xs}).
\]
When the closure of $\beta$ is a knot, we have $w(\beta) \equiv n \pmod 2$.

\begin{defn}\label{def:grdisc} 
Let $\beta$ be a braid whose closure is a knot. For $\xs \in (\beta \cup \mathbb I) \cs \cap \ct$, we define the \emph{punctured disc Alexander grading}
\begin{equation}\label{ALphi}
A^{\Pd}(\xs):=\deg_x \varphi(l_{\xs}) + \frac{w(\beta)+n}{2} \in \bZ.
\end{equation}
\end{defn}

The first author~\cite{Anghel2024AIF} obtained an intersection formula for the Alexander polynomial. We will show in Section~\ref{ss:Alex} that it is equivalent to the following:

\begin{thm}\label{thm:Alexander-disc}
Let $K$ be a knot in $S^3$ and $\beta \in B_{n+1}$ a braid such that $K = \hat{\beta}$.
Then the symmetrised Alexander polynomial
\begin{equation}\label{eqn:Alexander}
\Delta_K(x)= (-1)^{\frac{n(n+1)}{2}} \sum_{\xs\in \cI_{\bD_\beta}} \varepsilon_{\xs}\cdot x^{A^{\Pd}(\xs)} \in \bZ[x^{\pm 1}].
\end{equation}
Here, $\varepsilon_{\xs}$ is the intersection sign between $(\beta\cup \mathbb I)\cs$ and $\ct$ in the configuration space.
\end{thm}

We will prove that this topological model can be lifted to a categorification of the Alexander polynomial as follows.

\subsection{A Heegaard diagram associated with a braid}
\begin{figure}
\centering
\def\svgwidth{\textwidth}
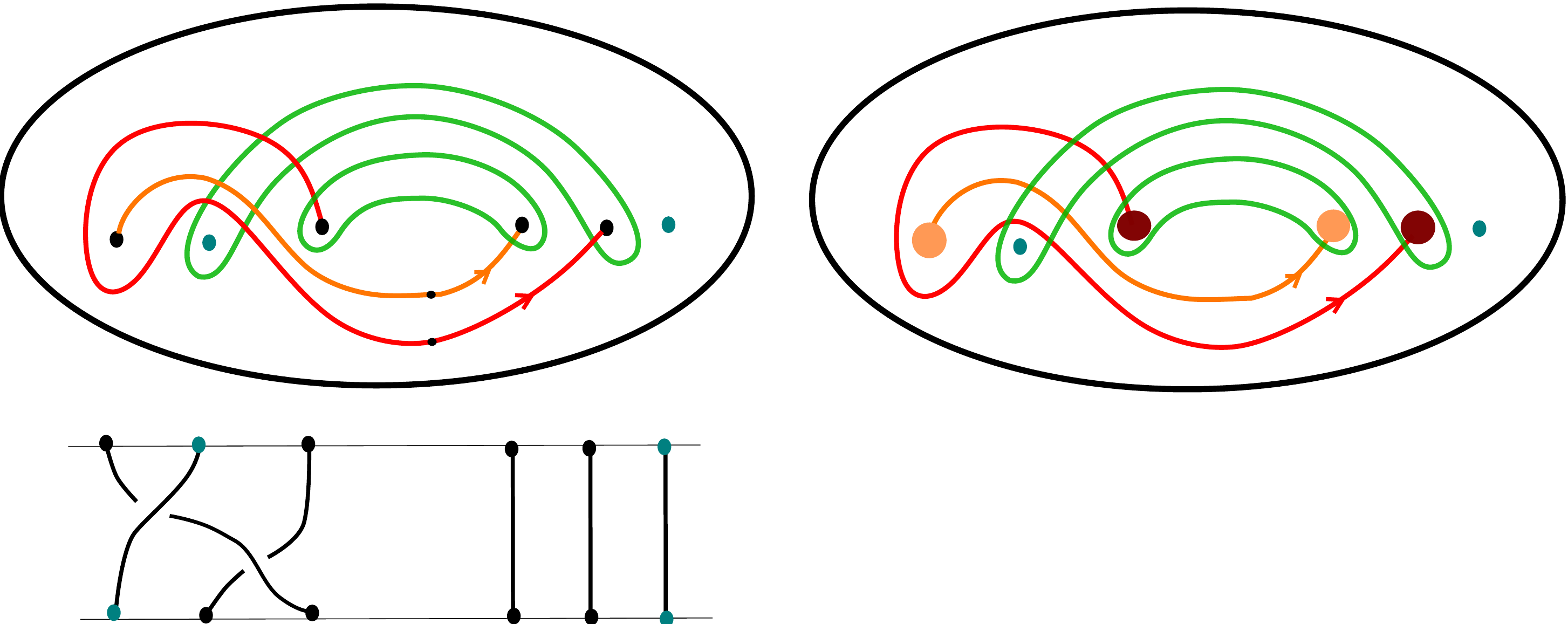
\caption{Top left: the disc model. Top right: the associated Heegaard diagram, where we perform surgery along the pairs of discs connected by arcs. Bottom left: the braid $\beta \cup \bI$.} 
\label{FigHS}
\end{figure}

We now modify the topological model in order to arrive at the setup for knot Floer homology. We start from the disc model $\bD_\beta$ and construct a Heegaard diagram $\cH_\beta$, as follows. Let $\pi$ be the permutation of the indices of the base points $\{z_0,\dots,z_n\}$ induced by the braid action of $\beta \cup \bI$. For example, the left-hand side of Figure~\ref{FigHS} shows the action of the braid $\sigma_2^{-1}\sigma_1 \cup \bI$ on $\cD_2$. We construct a Heegaard diagram from $\bD_\beta$ as follows.

\begin{defn}\label{def:HD}
Let $\Sigma$ be the surface obtained from $D^2$ by gluing a disc to $\partial D^2$ and performing surgery along the 0-spheres $\{w_i,z_{\pi(i)}\}$ for $i \in \{1,\dots,n\}$; see the right-hand side of Figure~\ref{FigHS}. We fix two base points $w := w_0$ and $z := z_{\pi{(0)}}$. Let $\as$ be the set of curves given by closing up the arcs $\ba'$ along the added handles. We orient $\as$ coherently with $\ba'$. Let $\bs$ be the image of the curves $\bb$ on $\Sigma$. We obtain the doubly pointed Heegaard diagram
\[
\cH_{\beta} := (\Sigma, \as, \bs, w, z).
\]
\end{defn}

For a fixed $i \in \{1,\dots,n\}$, up to isotopy, any two choices of $\alpha_i$ are related by a Dehn twist along the belt circle of the tube $T_i$ attached along $\partial N(\{w_i,z_{\pi(i)}\})$, which intersects $\as$ in a single point that lies on $\alpha_i$ and is disjoint from $\bs$. Hence, $\cH_\beta$ is unique up to diffeomorphism.

\begin{restatable}{thm}{mainthmHD}
\label{THD} 
Let $\beta$ be an $(n+1)$-braid with closure a knot $K$. Then the doubly pointed Heegaard diagram $\mathscr H_{\beta}$ in Definition~\ref{def:HD} represents the knot $K$.
\end{restatable}

We will prove Theorem~\ref{THD} in Section~\ref{Ss:surf}.

\subsection{Knot Floer homology categorifies the graded intersection} 
Let $\CFKh(\mathscr H_{\beta})$ be the knot Floer complex associated with the Heegaard diagram $\mathscr H_{\beta}$. Consider the tori $\mathbb T_{\alpha} := \alpha_1 \times \dots \times \alpha_n$ and $\mathbb T_{\beta} := \beta_1 \times \dots \times \beta_n$ in $\Sym^n(\Sigma)$. There is a canonical bijection between the generators $\cI_{\cH_\beta} := \bT_{\as} \cap \bT_{\bs}$ of this chain complex and the set of intersection points $\cI_{\bD_\beta}$ in the disc model. We let 
\[
\cI := \cI_{\bD_\beta} = \cI_{\cH_\beta}. 
\]
By Theorem~\ref{THD}, the homology of this complex is the knot Floer homology group $\HFKh(K)$. This has a bigrading provided by the Alexander grading $A^{\HF}$ and the Maslov grading $M$. In Section~\ref{S:identif}, we will show that the local system grading agrees with the Alexander grading:

\begin{restatable}{thm}{thmidgr}
\label{thm:Alexander-agree} 
Let $\beta \in B_{n+1}$ be a braid with closure a knot $K = \hat{\beta}$, and consider the corresponding disc model $\bD_\beta$ and doubly pointed Heegaard diagram $\cH_{\beta}$. The grading of an intersection point $\xs \in \cI$ in $\bD_\beta$ via the local system $\varphi$ and the Heegaard Floer Alexander grading agree:
\begin{equation}
A^{\Pd}(\xs) = A^{\HF}(\xs).
\end{equation}
\end{restatable}

Knot Floer homology is bigraded. Given $\xs \in \cI_{\cH_\beta}$, we write $M(\xs)$ for the \emph{Maslov grading} of $\xs$; this is the homological grading in the knot Floer complex.

\begin{prop}\label{prop:signs} 
Let $\beta$ be a braid whose closure is a knot. For $\xs \in \cI$, let $\varepsilon_{\xs}$ be the intersection sign of $(\beta \cup \bI)\cs$ and $\ct$ at $\xs$. Then we have
\[
(-1)^{\frac{n(n+1)}{2}} \varepsilon_{\xs} = (-1)^{M(\xs)}.
\]
\end{prop}

\begin{rmk}\label{rmk:signs}
Theorem~\ref{thm:Alexander-agree} and Proposition~\ref{prop:signs} together can be viewed as saying that knot Floer homology of $\cH_\beta$ categorifies, up to an overall sign, the intersection polynomial on the right-hand side of equation~\eqref{eqn:Alexander} defined using the disc model $\bD_\beta$.
\end{rmk}

\subsection{Unifying the Alexander and the Jones polynomials via quantum Heegaard diagrams} 
The topological model constructed by the first author~\cite{Anghel2024AIF} provides a unified model that captures both the Alexander and Jones polynomials through different specialisations of the variables.

\begin{defn}\label{def:canonical-intersection-point}
For $i \in \{1,\dots,n\}$, let $x_i^0$ be the point of $a_i' \cap b_i$ located on the right-hand side of $D^2$. The \emph{canonical intersection point} is
\[
\xs^0:=(x_1^0,\dots,x_n^0) \in \cI.
\]
For $i \in \{1,\dots,n\}$, let $q_i$ be a base point lying in the quadrant just below the canonical intersection point $x_i^0$, outside the disc $O_i$ bounded by the curve $b_i$ in $D^2$ and below the arc $a_i$, as shown in Figure~\ref{SurfJ}. Furthermore, let $q_0$ be a base point in the same region as $w_0$ and close to it. We write
\[
\qs := \{q_0,\dots,q_n\}.
\]
\end{defn}

\begin{defn}\label{def:quantum-disc-model}
    Let $\beta$ be an $(n+1)$-braid. Consider the punctured disc $\cD^q_n := \cD_n \setminus \qs$. The \emph{quantum disc model} associated with $\beta$ is
    \[
    \bD^q_\beta := (\cD^q_n, \ba', \bb).
    \]
    We write 
    \[
    C_n^q := \Conf_n(\cD_n^q) \subseteq C_n.
    \]
\end{defn}

\begin{figure}
\centering
\def\svgwidth{1.2\textwidth}
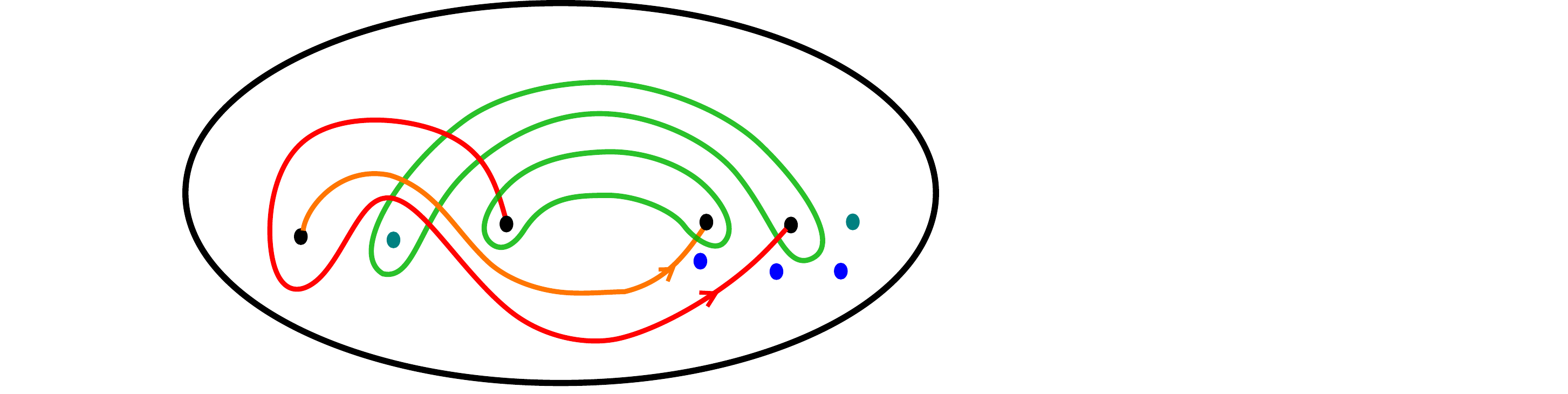
\caption{The unified disc model $\bD_\beta^q$ for the Alexander and Jones polynomials.}\label{fig:quantum-disc-model}
\end{figure}

\begin{figure}
\centering
\def\svgwidth{\textwidth}
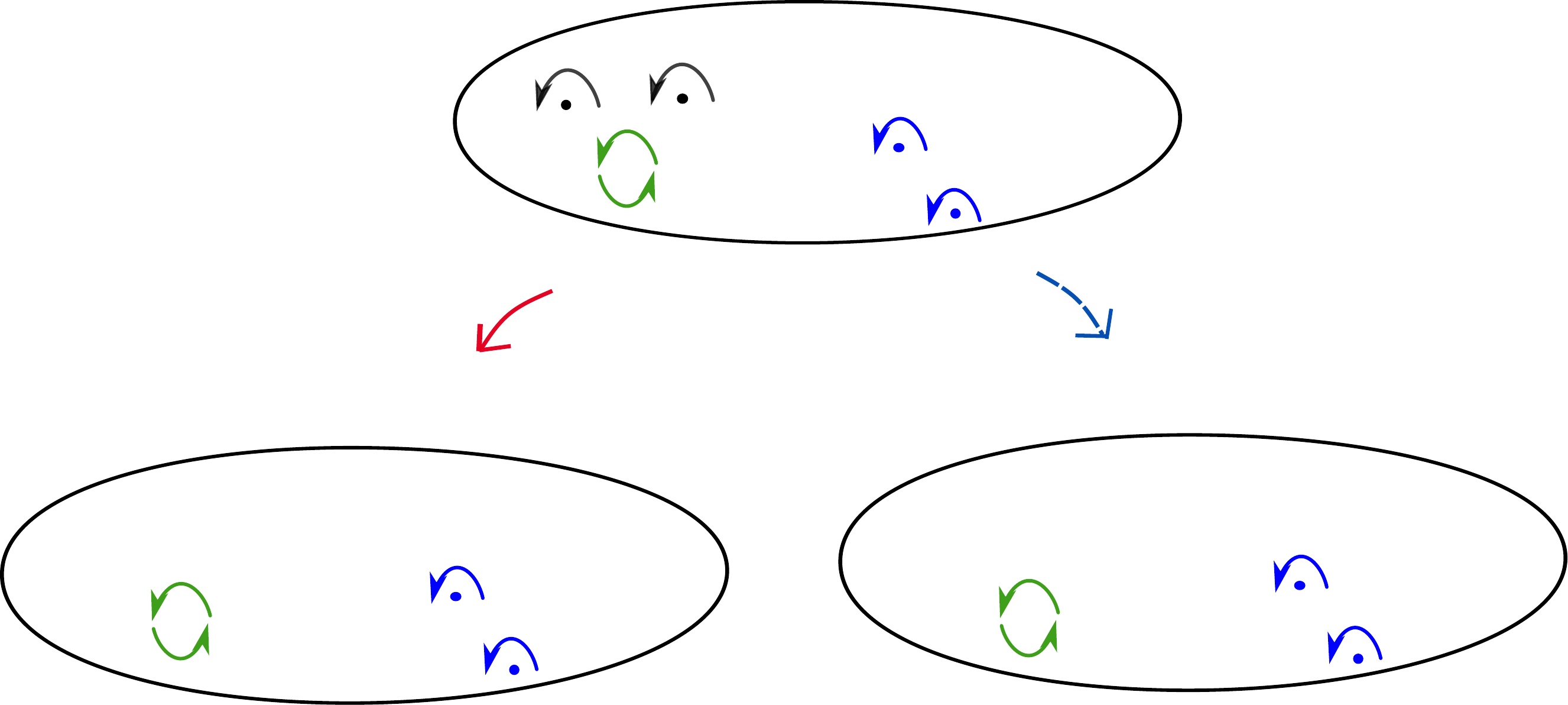
\caption{One can obtain from the quantum intersection polynomial $\Omega'_\beta(x,d)$ the Jones polynomial via the substitution $d=x^{-1}$ and the Alexander polynomial via the substitution $d=-1$. This figure shows the corresponding local systems: $\phi$ (top) and $\varphi$ (lower right).}
\label{fig:specialisations}
\end{figure}

For a schematic picture of the quantum disc model $\bD_\beta^q$, see Figure~\ref{fig:quantum-disc-model}. Let $\beta$ be an $(n+1)$-braid whose closure is a knot. The first author~\cite{Anghel2024AIF} defined an intersection polynomial
\[
\Omega'_\beta(x,d) \in \bZ\left[x^{\pm 1},d^{\pm 1}\right] 
\]
using the two-variable local system shown at the top of Figure~\ref{fig:specialisations} that recovers both the Alexander and the Jones polynomials; see Theorem~\ref{Topmod} and Section~\ref{S:unifmodel}. Based on this, we construct a unified model for the Alexander and Jones polynomials via a Heegaard diagram with additional data associated with $\beta$, as follows.

\begin{figure}
\centering
\def\svgwidth{1.5\textwidth}
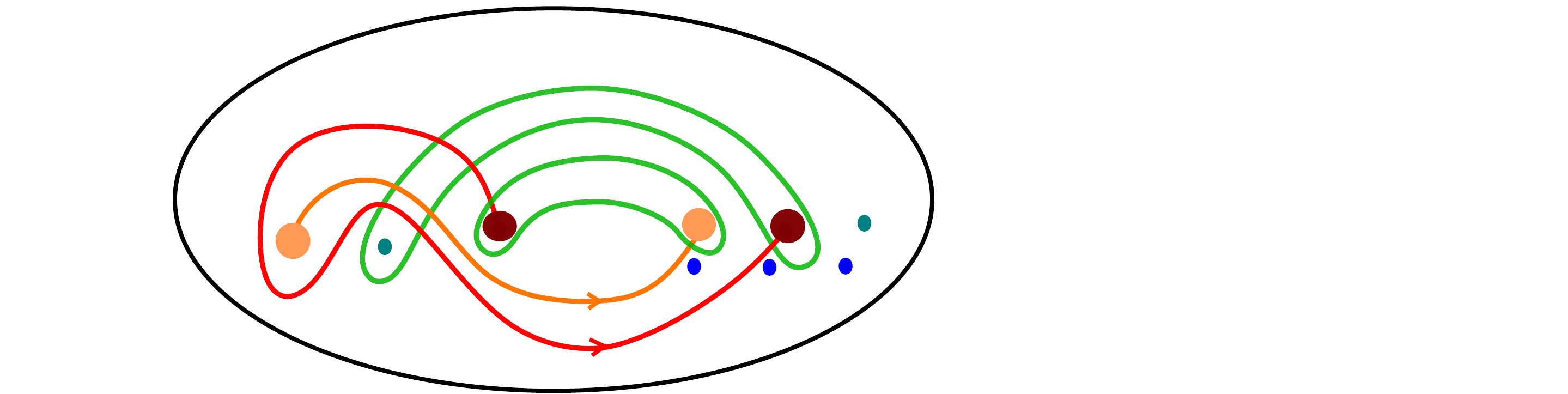
\caption{The quantum Heegaard diagram $\cH^q_\beta$ associated with a braid $\beta$.}\label{SurfJ}
\end{figure}

\begin{defn}\label{def:quantum-Heegaard-diagram}
    Let $\beta$ be an $(n+1)$-braid, and $\mathscr H_{\beta} = (\Sigma, \as, \bs, w, z)$ be the associated Heegaard diagram constructed in Definition~\ref{def:HD}. Pick a base point $h_i$ on the handle attached at $\{w_i,z_{\pi(i)}\}$ away from $\as$, and write $\qs = \{q_0,\dots,q_n\}$ and $\hs := \{h_1,\dots,h_n\}$. Then 
    \[
    \mathscr H^q_{\beta}:=(\Sigma,\as,\bs, w, z, \qs, \hs, \xs^0)
    \]
    is the \emph{quantum Heegaard diagram} associated with $\beta$.
\end{defn}

By construction, $\qs$, $\hs \in \Sigma \setminus (\as \cup \bs \cup \{w,z\})$. We will construct the \emph{quantum Alexander grading} $A^{\qHF}$ using the quantum Heegaard diagram $\cH^q_{\beta}$ in Definition~\ref{qAlex}.

\begin{defn}\label{def:quantum-intersection}
We define the \emph{quantum intersection polynomial}
\begin{equation}
\Omega^q_\beta(x,d):= (-1)^{\frac{n(n+1)}{2}}\sum_{\xs\in \cI} \varepsilon_{\xs} \, x^{A^{\HF}(\xs)} d^{A^{\qHF}(\xs)} \in \bZ[x^{\pm 1}, d^{\pm 1}].
\end{equation}
\end{defn}

In Definition~\ref{def:grdisc''}, we will introduce the quantum Alexander grading $A^{\qPd}$ using the punctured disc model. Analogously to Theorem~\ref{thm:Alexander-agree}, we will show that the quantum Alexander gradings defined via the punctured disc model and via quantum Heegaard diagrams agree.

\begin{restatable}{thm}{thmidqgr}
\label{thm:Jones-agree} 
Let $\beta \in B_{n+1}$ be a braid whose closure is a knot, and consider the corresponding quantum disc model $\bD_\beta^q$ and quantum Heegaard diagram $\mathscr H^q_{\beta}$. Then, for $\xs \in \cI$, the quantum Alexander grading in $\bD_\beta^q$ from Definition~\ref{def:grdisc''} and in $\cH^q_\beta$ from Definition~\ref{qAlex} agree:
\begin{equation}
A^{qPd}(\xs) = A^{\qHF}(\xs).
\end{equation}
\end{restatable}

For a knot $K$, we can unify the symmetrised Alexander polynomial $\Delta_K(x)$ and the Jones polynomial $V_{\mK}(x)$ of the mirror $\mK$ using quantum Heegaard diagrams:

\begin{restatable}{thm}{thmAlexanderJones}\label{qINT} 
Let $\beta \in B_{n+1}$ be a braid such that $K=\hat{\beta}$ is a knot. The quantum intersection polynomial $\Omega^q_\beta(x,d)$ defined using the quantum Heegaard diagram $\mathscr H^q_{\beta}$ satisfies
\[
\Omega^q_\beta(x,d) = \Omega'_\beta(x,-d).
\]
In particular, it unifies the Alexander and Jones polynomials as follows: 
\begin{equation}
\begin{aligned}
&\Delta_K(x) = \Omega^q_\beta(x,1),\\
&V_{\mK}(x) = \Omega^q_\beta(x,-x^{-1}),\\
\end{aligned}
\end{equation}
where $\mK$ is the mirror of $K$.
\end{restatable}

We will prove Theorem~\ref{qINT} in Section~\ref{sec:unification}. We summarise our models in Figure~\ref{fig:specialisations}.

\subsection{The Jones polynomial via quantum Heegaard diagrams} 
The intersection formula of Theorem~\ref{thm:Alexander-disc} for the Alexander polynomial arises from $\Omega^q_\beta(x,d)$ via the substitution $d=1$. By Remark~\ref{rmk:signs}, knot Floer homology categorifies this polynomial. On the other hand, the substitution $d=-x^{-1}$ leads to an intersection model for the Jones polynomial. Theorem~\ref{qINT} immediately implies the following:

\begin{restatable}{coro}{JonesviaqHD}
\label{TintJ}
Let $\beta$ be a braid whose closure is a knot $K$. The quantum Heegaard diagram $\mathscr H^q_{\beta}$ together with the classical and quantum Alexander gradings provides a Lagrangian intersection model for the Jones polynomial:
\begin{equation}
V_{\mK}(x)= (-1)^{\frac{n(n+1)}{2}} \sum_{\xs\in \cI} (-1)^{A^{\qHF}(\xs)} \varepsilon_{\xs} \, x^{A^{\HF}(\xs)-A^{\qHF}(\xs)}.
\end{equation}
\end{restatable}

\subsection{Quantum braid Floer homology}

Let $\beta$ be a braid whose closure is a knot $K$. The quantum Alexander grading $A^{\qHF}$ is not a filtration on the knot Floer complex 
\[
(\CFKh(\cH_\beta), \ph). 
\]
However, if we write $A^{\qHF}(\xs) = P(\xs) + W(\xs)$ for $\xs \in \cI$, where $W$ is introduced in Definition~\ref{def:W}, then $\cF(\xs) := (P(\xs), -W(\xs))$ is a bifiltration. We define the \emph{quantum braid Floer differential} $\partial^q$ to be the degree-$(0,0)$ part of $\ph$ with respect to this bifiltration. Its homology $\qHFB(\beta)$ is triply graded by $A^{\HF}$, $A^{\qHF}$, and the Maslov grading $M$. In Theorem~\ref{thm:braid-invariance}, we will show that $\qHFB(\beta)$ is an invariant of $\beta$ as a braid element. However, it is not invariant under conjugation or stabilisation. In particular, it is not an annular braid or a knot invariant. 

By Proposition~\ref{prop:categorification}, the graded Euler characteristic of $\qHFB(\beta)$ is $\Omega_\beta^q(x,d)$. In Theorem~\ref{thm:spectral-sequence}, we will show that $P - W$ is a filtration on the knot Floer complex, and the degree-zero part of $\ph$ is $\partial^q$; hence, there is a spectral sequence from $\qHFB(\beta)$ to $\HFKh(K)$. The differential $\partial^q$ preserves the set of canonical coordinates of a generator, which implies that $\xs^0$ determines a distinguished non-zero class $\Theta(\beta) \in \qHFB(\beta)$; see Proposition~\ref{prop:theta}.

\subsection{Related work} The term ``braid Floer homology'' is also used in the unrelated work of Ghrist, van den Berg, Vandervorst, and W\'ojcik~\cite{BFH} in Hamiltonian dynamics. Baldwin and Grigsby~\cite{Baldwin-Grigsby} studied braid conjugacy class invariants obtained from link Floer homology. 

Boninger~\cite{Boninger} associates to a braid $\rho \in B_n$ and multi-indices $\js, \ks \subseteq \{1, \ldots, n\}$ such that $|\js| = |\ks|$ a triply graded group $\widehat{\HFB}(\rho, \js, \ks)$ whose graded Euler characteristic is the determinant of the corresponding square submatrix of $\psi_n(\rho) - \lambda I_n$, where $\psi_n$ is the unreduced Burau representation. For codimension-one submatrices, the specialisation $\lambda = 1$ recovers the Alexander polynomial, and the corresponding groups admit a spectral sequence to the knot Floer homology of the mirror of the braid closure. Boninger also relates Heegaard Floer theory to the $U_q(\mathfrak{gl}(1|1))$ braid representation, whose associated link invariant is the Alexander polynomial. Although the Burau representation is closely related to the Hecke algebra framework underlying the Jones polynomial, no general Jones specialisation of Boninger's polynomials is established.

Our construction has a similar formal structure but a different geometric origin and decategorification. We define a single group $\qHFB(\beta)$ as the homology of the degree-$(0,0)$ part of the knot Floer differential with respect to the quantum bifiltration $\cF$. Its additional grading is the quantum Alexander grading $A^{\qHF}$, which involves intersection with the fat diagonal, multiplicities at the additional marked points $\qs$ and $\hs$, and the term $W(\xs)$ determined by the common coordinates of $\xs$ with the distinguished intersection point $\xs^0$. Its graded Euler characteristic is $\Omega_\beta^q(x,d)$, which specialises both to the Alexander polynomial and to the Jones polynomial of the mirror. No direct relationship between the two homology theories is presently known.

\subsection*{Acknowledgements} 
 The first author gratefully acknowledges the support of the ANR grant ANR-24-CPJ1-0026-01 at Universit\'e Clermont Auvergne - LMBP. She also acknowledges partial support from grants provided by the Ministry of Research, Innovation and Digitisation, CNCS - UEFISCDI, project numbers  PN-IV-P2-2.1-TE-2023-2040 and PN-IV-P1-PCE-2023-2001, within PNCDI IV. The second author used ChatGPT 5.6 Pro for proofreading the manuscript. 

\section{Unified topological model for the Alexander and Jones polynomials}\label{S:unifmodel}

The first author~\cite{Anghel2024AIF} defined a unified topological model that recovers the Alexander and Jones polynomials via intersections of embedded Lagrangians in the configuration space of the punctured disc. We outline the construction of this model in the case that will be relevant to us.

For $n \ge 1$, let $\beta \in B_{n+1}$ be a braid whose closure is a knot. Let $\bD^q_\beta$ be the associated quantum disc model from Definition~\ref{def:quantum-disc-model}. For $i \in \{1,\dots,n\}$, let
\[
d_i := \left(0,-\frac{n-i+1}{n+2}\right) \in a_i,
\]
and write $\ds := \{d_1,\dots,d_n\}$; see Figure~\ref{FigHS}.

\begin{figure}[h]
\begin{tikzpicture}[scale=0.95]
\foreach \x/\y in {-0.3/2,2/2,4/2,2/1,2.5/1,3/1.06} {\node at (\x,\y) [circle,fill,inner sep=1pt] {};}
\node at (2.2,2.6) [anchor=north east] {$z_i$};
\node at (3,2) [anchor=north east] {$\zeta_i$};
\node at (2.2,1) [anchor=north east] {\tiny $\mathrm d_1$};
\node at (2.8,1) [anchor=north east] {\tiny $\mathrm d_2$};
\node at (3.6,1) [anchor=north east] {\tiny $\mathrm d_n$};
\node at (2.63,2.3) [anchor=north east] {$\wedge$};
\draw (2,1.8) ellipse (0.4cm and 0.8cm);
\draw (1.9,1.8) ellipse (3.65cm and 1.65cm);
\foreach \x/\y in {7/2,9/2,11/2,8.9/1,9.5/1,10/1.07} {\node at (\x,\y) [circle,fill,inner sep=1pt] {};}
\node at (9.2,1) [anchor=north east] {\tiny $\mathrm d_1$};
\node at (9.8,1) [anchor=north east] {\tiny $\mathrm d_2$};
\node at (10.6,1) [anchor=north east] {\tiny $\mathrm d_n$};
\node at (9,1.5) [anchor=north east] {$\delta$};
\draw (9.4,1.8) ellipse (3.65cm and 1.65cm);
\draw (9.5,1)  arc[radius = 3mm, start angle= 0, end angle= 180];
\draw [->](9.5,1)  arc[radius = 3mm, start angle= 0, end angle= 90];
\draw (8.9,1) to[out=50,in=120] (9.5,1);
\draw [->](8.9,1) to[out=50,in=160] (9.25,1.12);
\end{tikzpicture}
\caption{The loop $\zeta_i$ (left) and $\delta$ (right).}\label{fig:loops}
\end{figure}

For $i \in \{0,\dots,n\}$, the loop $\zeta_i$ in $\cD^q_n$ is based at $d_1$ and goes around the $i$-th puncture once anticlockwise, as shown on the left of Figure~\ref{fig:loops}. We define the loop $\omega_i$ around $w_i$ and $\gamma_i$ around $q_i$ analogously. For $i \in \{0,\dots,n\}$, consider the loop
\[
\tilde{\zeta}_i(t) := \left(\zeta_i(t), d_2,\dots,d_n \right) 
\]
in $C_n^q$ for $t \in [0,1]$. We define the loops $\tilde{\omega}_i$ and $\tilde{\gamma}_i$ analogously. For $n=2$, the loop $\delta$ in the configuration space of two distinct unordered points swaps $d_1$ and $d_2$, as shown on the right of Figure~\ref{fig:loops}. For $n \ge 2$, consider the loop
\[
\tilde{\delta}(t) := \left(\delta(t), d_3,\dots,d_n \right) 
\]
in $C_n^q$ for $t \in [0,1]$. Let $\rho: \pi_1(C_n^q) \rightarrow H_1(C_n^q)$ be the abelianisation map. For $n \ge 2$, we have
\[
H_1(C_n^q) \cong \bZ^{3n+4},
\]
generated by the homology classes 
\[
[\tilde{\zeta}_0], \dots, [\tilde{\zeta}_n], [\tilde{\omega}_0], \dots, [\tilde{\omega}_n], [\tilde{\gamma}_0], \dots, [\tilde{\gamma}_n], [\tilde{\delta}].
\]
Furthermore, $H_1(C_1^q) \cong \bZ^6$, generated by
\[
[\tilde{\zeta}_0], [\tilde{\zeta}_1], [\tilde{\omega}_0], [\tilde{\omega}_1], [\tilde{\gamma}_0], [\tilde{\gamma}_1].
\]

\begin{defn}\label{localsystemg}  
For $n \ge 2$, we define the local system
\begin{equation}
\begin{aligned}
\phi \colon \pi_1\left(C_n^q \right) \stackrel{\rho}{\longrightarrow} & \ \bZ^{n+1} \oplus \ \bZ^{n+1} \oplus \bZ^{n+1} \oplus \ \bZ  \ \stackrel{f}{\longrightarrow}   \bZ[x^{\pm 1}, d^{\pm 1}]^\times, \\
& \langle [\tilde{\zeta}_i] \rangle \ \,\,\,\, \langle [\tilde{\omega}_i]\rangle \ \,\,\,\, \langle [\tilde{\gamma}_i]\rangle  \ \,\,\,\,  \langle [\tilde{\delta}]\rangle
\end{aligned}
\end{equation}
where $i \in \{0,\dots,n\}$, and $f$ is the group homomorphism
\begin{equation}\label{eqn:Localsystem}
f\left([\tilde{\zeta}_i]\right)= x \text{, } f\left([\tilde{\omega}_i]\right)= x^{-1} \text{, } f\left([\tilde{\gamma}_i]\right)= -d \text{, and } f\left([\tilde{\delta}]\right)= d;
\end{equation}
see the top of Figure~\ref{fig:specialisations}.

For $n = 1$, we let
\begin{equation}
\begin{aligned}
\phi \colon \pi_1\left(C_1^q \right) \stackrel{\rho}{\longrightarrow} & \ \bZ^2 \oplus \ \bZ^2 \oplus \bZ^2  \ \stackrel{f}{\longrightarrow}   \bZ[x^{\pm 1}, d^{\pm 1}]^\times, \\
& \langle [\tilde{\zeta}_i] \rangle \  \langle [\tilde{\omega}_i]\rangle \ \langle [\tilde{\gamma}_i]\rangle
\end{aligned}
\end{equation}
where $i \in \{0,1\}$, and $f$ is the group homomorphism
\begin{equation}
f\left([\tilde{\zeta}_i]\right)= x \text{, } f\left([\tilde{\omega}_i]\right)= x^{-1} \text{, and } f\left([\tilde{\gamma}_i]\right)= -d.
\end{equation}
\end{defn}

\subsection{Grading of the intersection points}
We consider paths $\eta_1,\dots,\eta_n$ in the punctured disc such that $\eta_i$ for $i \in \{1,\dots,n\}$ connects $d_i$ to the left-hand side of the curve $b_i$, as in Figure~\ref{fig:disc-model}. We also fix paths $\eta'_1,\dots,\eta'_n$ such that $\eta_i'$ connects $d_i$ to the right-hand side of the curve $b_i$.

The braid group $B_{3n+3}$ is the mapping class group of the punctured disc $\cD_n^q$, where the mapping classes fix the boundary pointwise. Consider a braid $\beta \in B_{n+1}$ and let $\mathbb I \in B_{2n+2}$ be the trivial $(2n+2)$-braid. We act on $\cD_n^q$ by $\beta \cup \mathbb I \in B_{3n+3}$ and write $a_i' := (\beta \cup \mathbb I) a_i$. Consider the submanifold
\[
(\beta \cup \mathbb I) \cs := a_1' \times \dots \times a_n' \subseteq C^q_n.
\]
For simplicity, we fix a representative of the action of $\beta \cup \bI$ that is supported in a disc containing the punctures $\{z_0,\dots,z_n\}$ that lies in the half-plane $\{x<0\}$. We choose a representative such that $(\beta \cup \mathbb I) \cs$ is a Lagrangian submanifold, as in \cite[Section~2.2]{Anghel2023TAMS}.

We now define the graded intersection $\langle (\beta\cup \bI) \cs,\ct \rangle$. It is parametrised by the set of intersection points 
\begin{equation}
\cI_{\bD_\beta}:=(\beta \cup \mathbb I) \cs\cap \ct.
\end{equation}
These are graded via the local system $\phi$ defined on the configuration space, as follows. To each intersection point $\xs \in \cI_{\bD_\beta}$, we first associate a loop $l_{\xs}$ in the configuration space, on which we evaluate the homomorphism $\phi$:
\[
\xs \in \cI_{\bD_\beta} \ \ \rightsquigarrow \text{ loop } l_{\xs} \ \ \rightsquigarrow \text{ grading }  \phi(l_{\xs}).
\]

\begin{figure}
\centering
\def\svgwidth{0.5\textwidth}
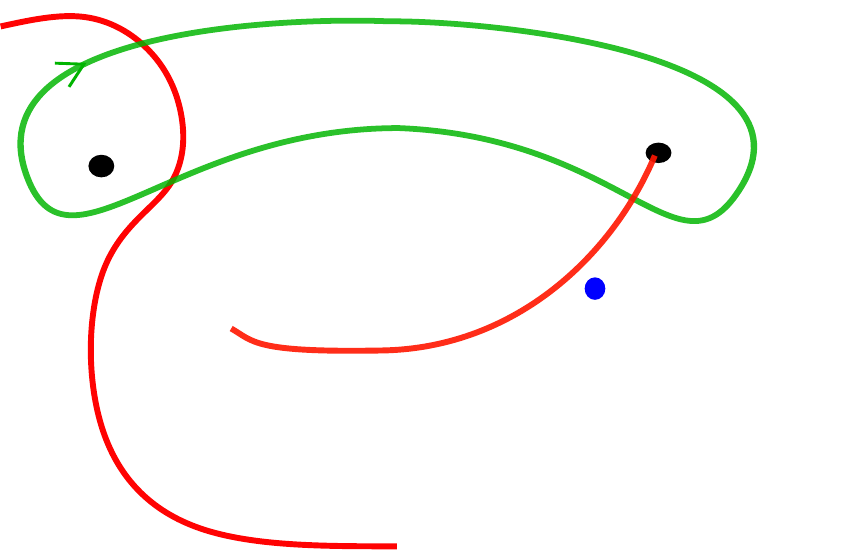
\caption{Construction of the loop $l_{\xs}$: follow $\eta_i$ or $\eta_i'$ from $d_i$ to $b_i$, then go along $b_i$ to $x_i$, and finally follow $\ba'$ back to $\ds$.}\label{fig:l_x}
\end{figure}

\begin{defn}\label{def:loop}
Let $\xs=(x_1,\dots,x_n) \in \cI_{\bD_\beta}$, where $x_i \in a_{\pi(i)}' \cap b_i$ for a permutation $\pi \in S_n$. For $i \in \{1,\dots,n\}$, we distinguish two cases: 
\begin{enumerate}
    \item If the point $x_i$ lies on the left-hand side of $b_i$, let $\tau_i$ be the path starting at $d_i$, continuing along $\eta_i$, and then following the left-hand side of $b_i$ to the point $x_i$; see Figure~\ref{fig:l_x}.
    \item If $x_i$ lies on the right-hand side of $b_i$, we define $\tau_i$ to be the path defined by following $\eta'_i$ and then continuing along the right-hand side of $b_i$ to $x_i$.
\end{enumerate}
We consider the path
\begin{equation}
\tau_{\xs} := \tau_1 \times \dots \times \tau_n
\end{equation}
in the configuration space $C^q_n$ from $\ds$ to $\xs$.

Let $\nu_i$ be the path that starts at $x_i$ and ends at $d_{\pi(i)}$, following the arc $a_{\pi(i)}'$. We consider the associated path
\begin{equation}
\nu_{\xs} := \nu_1 \times \dots \times \nu_n 
\end{equation}
in the configuration space $C^q_n$ from $\xs$ to $\ds$. The loop associated with the intersection point $\xs$ is given by the concatenation  
\begin{equation}
l_{\xs} := \tau_{\xs} \nu_{\xs}, 
\end{equation}
which is a loop based at $\ds$.
\end{defn}

\begin{rmk}
    In the definition of $l_{\xs}$, irrespective of how we connect $\eta_i(1)$ or $\eta_i'(1)$ with $x_i$ along $b_i$, the evaluation $\phi(l_{\xs})$ is always the same. Indeed, for $i \in \{1,\dots,n\}$, choose an anticlockwise parametrisation of $b_i$ based at either $\eta_i(1)$ or $\eta_i'(1)$, depending on whether $x_i$ is on the left or right of $b_i$. Any two choices of $l_{\xs}$ differ by the image of the product loop $b_1^{k_1} \times \dots \times b_n^{k_n}$ in $C_n^q$ for some $k_1, \dots, k_n \in \bZ$. In $H_1(C_n^q)$, we have
    \[
    [b_1^{k_1} \times \dots \times b_n^{k_n}] = \sum_{i=1}^n k_i\left([\tilde{\zeta}_i] + [\tilde{\omega}_i] \right),
    \]
    where the coefficient of $[\tilde{\delta}]$ is zero since the curves $b_1,\dots,b_n$ are not nested. As $f([\tilde{\zeta}_i] + [\tilde{\omega}_i]) = xx^{-1} = 1$, we obtain that $\phi(l_{\xs})$ is independent of the choice of connecting paths on $b_1,\dots,b_n$.
\end{rmk}

Now we present the definition of the graded intersection, which was introduced by the first author~\cite[Definition~3.8]{Anghel2024AIF}. 

\begin{defn}\label{def:graded-intersection} 
We define the \emph{graded intersection} between $(\beta\cup \mathbb I)\cs$ and $\ct$ as 
\begin{equation}\label{int0}
\langle (\beta\cup \bI) \cs,\ct \rangle:= \sum_{\xs = (x_1,\dots,x_n)\in \cI_{\bD_\beta}} \varepsilon_{x_1} \dots \varepsilon_{x_n} \cdot \phi(l_{\xs}) \in \bZ[x^{\pm 1}, d^{\pm 1}],
\end{equation}
where $\varepsilon_{x_k}$ is the intersection sign between $a_{\pi(k)}'$ and $b_k$ at $x_k$ for $k \in \{1,\dots,n\}$.
\end{defn}

\begin{defn}
    For $\xs = (x_1,\dots,x_n) \in \cI_{\bD_\beta}$, let $\pi_{\xs} \in S_n$ be the permutation such that $x_i \in a_{\pi_{\xs}(i)}' \cap b_i$ for $i \in \{1,\dots,n\}$.
\end{defn}

\begin{lem}\label{lem:signs}
    For $\xs = (x_1,\dots,x_n) \in \cI_{\bD_\beta}$, let $\varepsilon_{\xs}$ be the intersection sign between $(\beta \cup \bI) \cs$ and $\ct$ at $\xs$. Then
    \[
    \varepsilon_{\xs} = (-1)^{\frac{n(n-1)}{2}} \sgn(\pi_{\xs}) \varepsilon_{x_1} \dots \varepsilon_{x_n},
    \]
    where $\sgn(\pi_{\xs})$ is the sign of $\pi_{\xs}$. 
\end{lem}

\begin{proof}
    Choose the lift $(x_1,\dots,x_n)\in (\cD_n^q)^n$ of $\xs \in C_n^q$ ordered by the curves $b_1,\dots,b_n$ and write $\pi := \pi_{\xs}$. Let $e_i$ be a positive tangent vector to $a'_{\pi(i)}$ at $x_i$, and let $v_i$ be a positive tangent vector to $b_i$ at $x_i$. By definition of the local intersection signs, the interleaved basis $(e_1,v_1,\dots,e_n,v_n)$ contributes $\varepsilon_{x_1}\dots \varepsilon_{x_n}$ relative to the product orientation of $(\cD_n^q)^n$. The submanifold $(\beta\cup\bI)\cs$ is oriented as the image of $a_1' \times \cdots \times a_n'$ under the quotient map. In the chosen $\bb$-ordered lift, this product orientation is represented by the ordered basis
    $(e_{\pi^{-1}(1)},\dots,e_{\pi^{-1}(n)})$. Reordering this basis to $(e_1,\dots,e_n)$ contributes the sign $\sgn(\pi_{\xs})$. Therefore, the ordered basis defining the intersection sign $\varepsilon_{\xs}$ differs from $(e_1,\dots,e_n,v_1,\dots,v_n)$ by the factor $\sgn(\pi)$. Finally, to pass from $(e_1,\dots,e_n,v_1,\dots,v_n)$ to the interleaved basis
    $(e_1,v_1,\dots,e_n,v_n)$, one performs 
    \[
    (n-1)+\dots+1=\frac{n(n-1)}{2}
    \]
    transpositions. Hence
    \[
    \varepsilon_{\xs} = (-1)^{\frac{n(n-1)}{2}} \sgn(\pi) \varepsilon_{x_1}\cdots \varepsilon_{x_n},
    \]
    as claimed.
\end{proof}

\subsection{Topological unification of the Alexander and Jones polynomials}

The following is \cite[Definition~3.10]{Anghel2024AIF}.

\begin{defn}\label{defn} Let $\beta \in B_{n+1}$ be a braid whose closure is a knot. We define the \emph{intersection polynomial} for the braid $\beta \in B_{n+1}$ as 
\begin{equation}
\Omega'_\beta(x,d):= x^{\frac{w(\beta)+n}{2}} d^{w(\beta)} \langle  (\beta \cup \mathbb I) \cs, \ct \rangle \in \bZ[x^{\pm 1}, d^{\pm 1}], 
\end{equation}
where $w(\beta)$ is the writhe of $\beta$. 
\end{defn}

The first author~\cite{Anghel2024AIF} showed that the intersection polynomial unifies the Alexander and Jones polynomials.

\begin{restatable}{thm}{Topmod}\label{Topmod} Let $K$ be a knot and $\beta \in B_{n+1}$ a braid such that $K = \hat{\beta}$. The intersection polynomial $\Omega'_\beta(x,d)$ recovers both the Alexander and Jones polynomials through
\begin{equation*}
{\Delta_K(x)} = \Omega'_\beta(x,-1) \text{ and }
{V_{\mK}(x)} = \Omega'_\beta(x,x^{-1}),
\end{equation*}
where $\mK$ is the mirror of $K$.
\end{restatable}

\begin{rmk}\label{rmk:Jones-convention}
In this paper, we use the version of the Jones polynomial that satisfies the skein relation
\[
x^{-1} V_{L_+}(x) - xV_{L_-}(x) = (x^{1/2} - x^{-1/2}) V_{L_0}(x)
\]
for a link $L$, and is normalised such that $V_U(x)=1$. Under this convention, the right-handed trefoil $3_1$ satisfies $V_{3_1}(x) = -x^4+x^3+x$. In \cite{Anghel2023TAMS} and \cite{Anghel2024AIF}, the first author used the Jones polynomial $J_L(x) = V_{\overline{L}}(x)$, where $\overline{L}$ is the mirror of $L$.
\end{rmk}

We now introduce two gradings on $\cI_{\bD_\beta}$.

\begin{defn}\label{def:grdisc''} 
Let $\beta$ be a braid whose closure is a knot. Let $\phi$ be the local system introduced in Definition~\ref{localsystemg}. For $\xs \in \cI_{\bD_\beta}$, we define the \emph{Alexander grading}
\begin{equation}
A^{\Pd}(\xs) := \deg_x \phi(l_{\xs}) + \frac{w(\beta)+n}{2} \in \bZ,
\end{equation}
and the \emph{quantum Alexander grading}
\begin{equation}
A^{\qPd}(\xs) := \deg_d \phi(l_{\xs}) + w(\beta) \in \bZ.
\end{equation}
\end{defn}

\subsection{Intersection model for the Alexander polynomial}\label{ss:Alex}
In this section, we discuss the above model in the case of the Alexander polynomial. The local system $\phi$ counts monodromies around the punctures, as shown at the top of Figure~\ref{fig:specialisations}. The local system $\varphi := \phi|_{d=-1}$ shown at the bottom right of Figure~\ref{fig:specialisations} recovers the Alexander polynomial by Theorem~\ref{Topmod}, so the $q$-punctures no longer play a role. Hence, we will consider only the set of punctures $\zs \cup \ws$. Consider the punctured disc $\cD_n := D^2 \setminus (\zs \cup \ws)$ and let $C_n:=\Conf_n(\cD_n)$ be the unordered configuration space of $n$ distinct points in $\cD_n$; see Figure~\ref{fig:disc-model}.

\begin{defn}\label{def:varphi} 
For $n \ge 2$, we define the local system $\varphi: \pi_1\left(C_n \right) \ \rightarrow \ \bZ[x^{\pm 1}]^\times$ by
\begin{equation}
\varphi(\tilde{\zeta}_i)=x \text{, }
\varphi(\tilde{\omega}_i)=x^{-1} \text{, and }
\varphi(\tilde{\delta})=-1
\end{equation}  
for $i \in \{0,\dots,n\}$. For $n=1$, we omit the generator $\tilde{\delta}$.
\end{defn}

For $\xs \in \cI_{\bD_\beta}$, we can write 
\[
\varphi(l_{\xs}) = \sgn(\varphi(l_{\xs})) x^{A^{\Pd}(\xs) - \frac{w(\beta)+n}{2}},
\]
where $\sgn(\varphi(l_{\xs})) \in \{1,-1\}$ and we introduced $A^{\Pd}(\xs)$ in Definition~\ref{def:grdisc}.

\begin{lem}\label{lem:delta}
    For $\xs \in \cI_{\bD_\beta}$, we have
    \[
    \sgn(\varphi(l_{\xs})) = \sgn(\pi_{\xs}).
    \] 
\end{lem}

\begin{proof}
    Recall that $\rho \colon \pi_1(C_n) \to H_1(C_n)$ is the abelianisation map. Let $\ell$ be the coefficient of $[\delta]$ in $\rho([l_{\xs}])$. Then $\ell$ is the number of times $l_{\xs}$ links the diagonal in $C_n$. By the definition of $\varphi$, we have $\sgn(\varphi(l_{\xs})) = (-1)^\ell$. The loop $l_{\xs}$ is based at $\ds = (d_1,\dots,d_n)$. Put $\pi := \pi_{\xs}$. The monodromy around $l_{\xs}$ maps $\ds$ to $\pi(\ds) := (d_{\pi(1)},\dots,d_{\pi(n)})$. We claim $(-1)^\ell = \sgn(\pi)$, which implies the result. Indeed, given a point $\us = (u_1,\dots,u_n) \in C_n$, consider the discriminant 
    \[
    \Delta(\us) := \prod_{i<j}(u_i-u_j)^2.
    \]
    This is a well-defined function on $C_n$ whose zero locus is the diagonal, where two coordinates coincide. As $l_{\xs}$ avoids the diagonal, the linking number $\ell$ is the winding number of $\Delta \circ l_{\xs}$ around 0. Lift $l_{\xs}$ to the ordered configuration space from $\ds$. The lift $\tilde{l}$ ends at $\pi(\ds)$. Let 
    \[
    V(\us) := \prod_{i<j} (u_i-u_j)
    \]
    be the Vandermonde product. Then $\Delta = V^2$, and
    \[
    V(\pi(\ds)) = \sgn(\pi) V(\ds).
    \]
    So, along $\tilde{l}$, the Vandermonde $V$ changes by a factor of $\sgn(\pi)$. Squaring gives the discriminant winding: $\ell$ is even if $\pi$ is even, and $\ell$ is odd if $\pi$ is odd.
\end{proof}

We are now ready to prove Theorem~\ref{thm:Alexander-disc}.

\begin{proof}[Proof of Theorem~\ref{thm:Alexander-disc}]
    By Theorem~\ref{Topmod}, Definition~\ref{defn} (together with the fact that $w(\beta) \equiv n \pmod 2$ as $\hat{\beta}$ is a knot), Definition~\ref{def:graded-intersection}, Lemma~\ref{lem:signs}, and Lemma~\ref{lem:delta}, we have 
    \[
    \begin{split}
    {\Delta_K(x)} &=\Omega'_\beta|_{d=-1} \\
    &= x^{\frac{w(\beta)+n}{2}} \cdot (-1)^n \langle  (\beta \cup \mathbb I) \cs, \ct \rangle|_{d=-1} \\
    &= x^{\frac{w(\beta)+n}{2}} \cdot (-1)^n \sum_{\xs = (x_1,\dots,x_n)\in \cI_{\bD_\beta}} \varepsilon_{x_1} \dots \varepsilon_{x_n} \cdot \phi(l_{\xs})|_{d=-1} \\
    &= (-1)^{\frac{n(n-1)}{2}} (-1)^n x^{\frac{w(\beta)+n}{2}} \sum_{\xs \in \cI_{\bD_\beta}} \varepsilon_{\xs} \cdot \sgn(\pi_{\xs}) \varphi(l_{\xs}) \\
    &= (-1)^{\frac{n(n+1)}{2}} \sum_{\xs \in \cI_{\bD_\beta}} \varepsilon_{\xs} \cdot x^{A^{\Pd}(\xs)},
    \end{split}
    \]
    as claimed.
\end{proof}

\section{Heegaard diagrams}\label{Ss:surf}

Heegaard diagrams are a useful tool for encoding $3$-manifolds and knots in $3$-manifolds and form the main input for the construction of Heegaard Floer homology. In this section, we review the relevant definitions and prove Theorem~\ref{THD}.

\begin{defn} A \emph{Heegaard diagram} is a tuple $\mathscr H=(\Sigma,\as,\bs)$, where
\begin{itemize}
\item $\Sigma$ is a closed, connected, oriented surface of genus $g$,
\item $\as=\{\alpha_1, \dots, \alpha_g\}$ and $\bs=\{\beta_1, \dots, \beta_g\}$ are sets of pairwise disjoint simple closed curves on $\Sigma$ that are transverse to each other, and
\item both $\as$ and $\bs$ are linearly independent in $H_1(\Sigma)$ (equivalently, generate a rank $g$ subgroup).
\end{itemize}
\end{defn}

Given a Heegaard diagram, one can obtain a closed, connected, oriented $3$-manifold $Y$ by attaching 3-dimensional 2-handles to $\Sigma \times [0,1]$ along $\as \times \{0\}$ and $\bs \times \{1\}$, and attaching two $3$-balls to the resulting $S^2$ boundary components. This exhibits $Y$ as the union of two handlebodies, which we call the $\alpha$- and $\beta$-handlebodies.

\begin{defn}
Let $K \subseteq S^{3}$ be a knot. Then $\mathscr H=(\Sigma,\as,\bs, w, z)$ is a \emph{doubly pointed Heegaard diagram for $K$} if $\mathscr H=(\Sigma,\as,\bs)$ is a Heegaard diagram for $\mathbb S^3$. Furthermore, $w$ and $z \in \Sigma \setminus (\as \cup \bs)$ are base points, and the knot $K$ can be obtained by connecting $w$ and $z$ with an arc in $\Sigma \setminus \as$ and pushing it into the $\alpha$-handlebody, then joining $z$ to $w$ with an arc in $\Sigma \setminus \bs$ and pushing it into the $\beta$-handlebody.
\end{defn}

Suppose that $(\Sigma,\as,\bs, w, z)$ is a doubly pointed Heegaard diagram for the knot $K$. Let $\Sigma' := \Sigma \setminus N(\{w, z\})$. If we attach 3-dimensional 2-handles to $\Sigma' \times [0,1]$ along $\as \times \{0\}$ and $\bs \times \{1\}$, we obtain the exterior of $K$, and $\partial \Sigma' \times \{0\}$ is the union of two oppositely oriented meridians of $K$. We now prove Theorem~\ref{THD}, which we restate below.

\mainthmHD*

\begin{proof}
We start from the disc model $\bD_\beta = (\cD_n, \ba', \bb)$. We view $S^2$ as the equator of $S^3$ and $\cD_n \subseteq S^2$. Let $\Sigma' := \Sigma \setminus N(\{w, z\})$ and $S' := S^2 \setminus N(\ws \cup \zs)$.

\medskip\noindent {\it Lower hemisphere: the braid $\beta \cup \mathbb \bI$ and the cups.} For $i \in \{1,\dots,n\}$, we embed the tube $T_i$ attached to $S'$ along $\partial N(\partial a_i')$ in the lower hemisphere $S^3_-$, such that it goes vertically down over $\partial N(\partial a_i')$ and horizontally shadows $a_i'$ when projected onto $S^2$. 

We obtain $\alpha_i$ from $a_i'$ by going vertically down from $a_i'(0)$, following the bottom of the associated tube, and then going up to $a_i'(1)$. A $3$-dimensional $2$-handle attached to $\Sigma' \times \{0\}$ along $\alpha_i \times \{0\}$ is a wall with base $a_i' \times \{0\}$. Let $a_0' := (\beta \cup \bI)a_0$. Hence, the 3-manifold $\Sigma' \times [0,1]$ with 3-dimensional 2-handles attached along $\alpha_i \times \{0\}$ for $i \in \{1,\dots,n\}$ is the manifold $S^3_- \setminus N(\tilde{\ba}')$, where $\tilde{\ba}' = \{\tilde{a}_0',\dots,\tilde{a}_n'\}$ and $\tilde{a}_j'$ is $a_j'$ with its interior pushed into $S^3_-$ for $j \in \{0,\dots,n\}$. Here, the annulus $\partial N(\tilde{a}_j') \setminus S^2$ corresponds to $T_j \times \{1\}$ for $j \in \{1,\dots,n\}$ and $\partial N(\tilde{a}_0') \setminus S^2$ corresponds to $(\partial \Sigma') \times [0,1]$ together with $\Sigma' \times \{0\}$ compressed along $\as \times \{0\}$ (the lower boundary after attaching the 3-dimensional 2-handles).

Let $\tilde{\ba}$ be the arcs $\ba \cup \{a_0\}$ with their interiors pushed into $S^3_-$. This is the collection of cups in $S^3_-$. Since $\ba' \cup \{a_0'\} = (\beta \cup \bI) (\ba \cup \{a_0\})$, there is a diffeomorphism 
\[
d \colon (S^3_-, \tilde{\ba}) \to (S^3_-, \tilde{\ba}') 
\]
such that $d|_{S^2}$ is the action of $\beta \cup \bI$ on $(S^2, \zs \cup \ws)$. It follows that $S^3 \setminus N(\tilde{\ba}')$ is the complement of the braid $\beta \cup \bI$ together with the set of cups in $S^3_-$.

\medskip\noindent {\it Upper hemisphere: caps.} In equation~\eqref{eqn:a_i-c_i}, we defined the half-circle arcs $c_0, \dots, c_n$ connecting the punctures $\zs$ and $\ws$. If we attach $3$-dimensional $2$-handles to $S' \times [0,1]$ along $\bs \times \{1\}$ in $S^3_+$, then we obtain the 3-manifold $S^3_+ \setminus N(\tilde{\boldsymbol{c}})$, where $\tilde{\boldsymbol{c}} = \{\tilde{c_0}, \dots, \tilde{c}_n\}$ and $\tilde{c}_j$ is $c_j$ with its interior pushed into $S^3_+$ for $j \in \{0,\dots,n\}$. Note that $\tilde{\boldsymbol{c}}$ is a collection of caps. Here, the annulus $\partial N(\tilde{c_0}) \setminus S^2$ corresponds to $\partial N(\{w_0,z_0\}) \times [0,1]$ together with $S' \times \{1\}$ compressed along $\bs \times \{1\}$ (the upper boundary after attaching the 3-dimensional 2-handles).

\medskip\noindent {\it Identifying the knot complement.} Since $S' = \Sigma' \setminus (T_1 \cup \dots \cup T_n)$, the doubly pointed Heegaard diagram $\cH_{\beta}$ represents 
\[
\left(S^3_- \setminus N(\tilde{\ba}')\right) \cup \left(S^3_+ \setminus N(\tilde{\boldsymbol{c}})\right), 
\]
which is the complement of the knot obtained as $\tilde{\ba}' \cup \tilde{\boldsymbol{c}}$. We conclude that the knot represented by $\cH_\beta$ is the braid closure $\hat{\beta}$.
\end{proof}

\section{Domains and the Alexander grading}

In this section, we review domains in Heegaard diagrams and simultaneously extend the terminology to the disc model. We then review the Alexander grading $A^{\HF}$ in knot Floer homology. Heegaard Floer homology of closed 3-manifolds was defined by Ozsv\'ath and Szab\'o~\cite{OS, OS1}, extended to knots by Ozsv\'ath and Szab\'o~\cite{OS2} and independently by Rasmussen~\cite{Rass}, and to sutured manifolds by the second author~\cite{Juhasz1}.

\subsection{Whitney discs and their domains}

We first recall the notion of topological Whitney discs from Lagrangian intersection Floer homology.

\begin{defn} 
Let $L_0$ and $L_1$ be Lagrangians in a symplectic manifold $M$ that intersect transversely in a set $\cI$. For $x$ and $y \in \cI$, a \emph{topological Whitney disc} from $x$ to $y$ is a continuous map $u \colon D^2 \rightarrow M$
such that \begin{equation}
\begin{aligned}
&u(-i)=x, \ u(i)=y,\\
&u(\partial_+ D^2) \subseteq L_0, \ u(\partial_- D^2) \subseteq L_1,
\end{aligned}
\end{equation}
where $\partial_+ D^2 = S^1 \cap \{z \in \mathbb C : \text{Re}(z) \ge 0\}$ and $\partial_- D^2 = S^1 \cap \{z \in \mathbb C : \text{Re}(z) \le 0\}$. We denote the set of homology classes of Whitney discs from $x$ to $y$ by $\pi_2(x,y)$. 
\end{defn}

Let 
\[
\mathscr H = (\Sigma, \as = \{\alpha_1, \dots, \alpha_n\}, \bs = \{\beta_1,\dots,\beta_n\})
\]
be a Heegaard diagram, where $\Sigma$ is a closed, connected, and oriented surface of genus $n$. Then $\mathbb T_{\alpha} := \alpha_1 \times \dots \times \alpha_n$ and $\mathbb T_{\beta} := \beta_1 \times \dots \times \beta_n$ are Lagrangian tori in $\Sym^n(\Sigma)$ for a suitable choice of symplectic structure; see Perutz~\cite{Perutz}. If $\as$ and $\bs$ intersect transversely, then so do $\mathbb T_\alpha$ and $\mathbb T_\beta$. Let 
\[
\cI_{\cH} := \mathbb T_{\alpha} \cap \mathbb T_{\beta};
\]
this is a set of generators of the Heegaard Floer chain complex.

Furthermore, let $\beta$ be an $(n+1)$-braid, and consider the corresponding disc model 
\[
\bD_\beta = (\cD_n, \ba', \bb),
\]
where $\cD_n = D^2 \setminus (\zs \cup \ws)$. This gives rise to transverse Lagrangians $(\beta \cup \bI)\cs$ and $\ct$ in $C_n$.

\begin{defn} We call the closures of the connected components of $\Sigma \setminus (\as \cup \bs)$ and $D^2 \setminus (\ba' \cup \bb)$ \emph{regions}. Consider the free abelian groups
\[
\begin{split}
D(\cH) &:= \bZ\langle \pi_0(\Sigma \setminus (\as \cup \bs)) \rangle; \\
D(\bD_\beta) &:= \bZ\langle \pi_0(D^2 \setminus (\ba' \cup \bb)) \rangle.
\end{split}
\]
We call elements of $D(\cH)$ and $D(\bD_\beta)$ \emph{domains}. We say that a domain $D$ is \emph{positive} and write $D \ge 0$ if all of its coefficients are non-negative.
\end{defn}

Note that a domain can be viewed as a 2-chain in $\Sigma$ or $D^2$.

\begin{defn}
For $p \in \Sigma \setminus (\as \cup \bs)$, let
\[
V_p:=\{ (x_1,\dots,x_n) \in \Sym^n(\Sigma) \mid \exists \ 1\leq i \leq n : x_i=p
\}.
\]
This is a complex hypersurface in $\Sym^n(\Sigma)$. For $p \in D^2 \setminus (\ba' \cup \bb)$, we define $V_p \subseteq \Sym^n(D^2)$ analogously.
\end{defn}

For $\xs$, $\ys \in I_{\cD_\beta}$, we use $\pi_2(\xs,\ys)$ to denote homology classes of topological Whitney discs in $\Sym^n(D^2)$, with boundary on $(\beta \cup \bI)\cs$ and $\ct$. Although the boundary of such a disc lies in the configuration space $C_n$, its interior is allowed to meet the fat diagonal $\Delta$ and the divisors $V_p$ for punctures $p \in \zs \cup \ws$. These intersections compute the corresponding linking numbers of the boundary loop with $\Delta$ and $V_p$ in $C_n$. Following Ozsv\'ath and Szab\'o~\cite{OS2}, there is a correspondence between topological Whitney discs and domains.

\begin{defn}\label{def:multiplicity}
Let $\xs$, $\ys \in \cI_{\cH}$ or $\cI_{\bD_\beta}$ be intersection points and $u \in \pi_2(\xs, \ys)$ a Whitney disc. Then, for $p \in \Sigma \setminus (\as \cup \bs)$ or $D^2 \setminus (\ba' \cup \bb)$, the \emph{multiplicity of $u$ at $p$} is the algebraic intersection number
\[
n_p(u):= u \cdot V_p,
\]
which is defined by taking a smooth representative of $u$ transverse to $V_p$. Let $R_1,\dots,R_d$ be the regions and choose a point $p_i \in R_i$ for every $i \in \{1,\dots,d\}$. The \emph{domain of $u$} is 
\[
\mathcal D(u) := \sum_{i=1}^d n_{p_i}(u) \cdot R_i.
\]
\end{defn}

\begin{defn}\label{def:xy-domains}
    Given $\xs$, $\ys \in \cI_{\cH}$, we say that $D \in D(\cH)$ is a \emph{domain from $\xs$ to $\ys$} if $\partial D \cap \as$ is a 1-chain with boundary $\ys - \xs$ and $\partial D \cap \bs$ is a 1-chain with boundary $\xs - \ys$. We denote the set of domains from $\xs$ to $\ys$ by $D_{\cH}(\xs,\ys)$. If $\xs$, $\ys \in \cI_{\bD_\beta}$, we define when $D \in D(\bD_\beta)$ is a domain from $\xs$ to $\ys$ and the set of domains $D_{\beta}(\xs,\ys)$ from $\xs$ to $\ys$ analogously, but we also require that $n_{w_0}(D) = 0$ to ensure that $D$ has multiplicity zero along $\partial D^2$.
\end{defn}

\begin{rmk}
If $\xs$, $\ys \in \cI_{\cH}$ or $\cI_{\bD_\beta}$ and $u \in \pi_2(\xs,\ys)$, then $\cD(u)$ is a domain from $\xs$ to $\ys$. Conversely, every $D$ in $D_{\cH}(\xs,\ys)$ or $D_\beta(\xs,\ys)$ is the domain of some $u \in \pi_2(\xs,\ys)$.
\end{rmk}

\begin{defn}
    We say that the domain $D \in D(\cH)$ is \emph{periodic} if $\partial D$ is a linear combination of curves in $\as \cup \bs$. Similarly, $D \in D(\bD_\beta)$ is \emph{periodic} if $\partial D$ is a linear combination of curves in $\bb$.
\end{defn}

The following is a standard result from Heegaard Floer homology.

\begin{lem}
    Let $D$ and $D'$ be domains in $D_{\cH}(\xs,\ys)$ or $D_{\beta}(\xs,\ys)$. Then $D-D'$ is a periodic domain.
\end{lem}

\begin{lem}\label{lem:periodic}
    Let $O_i$ be the disc that the curve $b_i$ bounds in $D^2$ for $i \in \{1,\dots,n\}$. Then the periodic domains in $D(\bD_\beta)$ are exactly the $\bZ$-linear combinations of $O_1,\dots,O_n$.  
\end{lem}

\begin{proof}
    Suppose that $D$ is a periodic domain. Then $\partial D$ is a $\bZ$-linear combination of $b_1,\dots,b_n$. Note that $D^2 \setminus \bb$ consists of $O_1 \sqcup \dots \sqcup O_n$ and a component $E$ that contains $\partial D^2$. Hence, $D = \sum_{i=1}^n \lambda_i O_i + \lambda E$. But $\partial D^2$ appears with multiplicity zero in $\partial D$; hence $\lambda = 0$.
\end{proof}

Suppose that the braid $\beta$ induces the permutation $\pi$ of the indices of the punctures $z_0,\dots,z_n$. Then $a_i'(0) = z_{\pi(i)}$ and $a_i'(1) = w_i$.

\begin{defn}
    We say that the domain $D \in D(\bD_\beta)$ is \emph{matched} if $n_{w_i}(D) = n_{z_{\pi(i)}}(D)$ for $i \in \{1,\dots,n\}$.
\end{defn}

\begin{prop}\label{prop:matched}
    Let $\beta$ be a braid whose closure is a knot, and let $\xs$, $\ys \in \cI_{\bD_\beta}$. Then there is a unique matched domain $D_{\xs,\ys}$ in $D_{\beta}(\xs,\ys)$.
\end{prop}

\begin{proof}
    First, we show that $D_{\beta}(\xs,\ys) \neq \emptyset$. Let $\ell_{\xs,\ys}$ be a 1-cycle obtained by connecting $\xs \cap a_i'$ to $\ys \cap a_i'$ along $a_i'$ and $\ys \cap b_i$ to $\xs \cap b_i$ along $b_i$ for every $i \in \{1,\dots,n\}$. Since $H_1(D^2) = 0$, the 1-cycle $\ell_{\xs,\ys}$ is 0-homologous; hence, it bounds a 2-chain $\zeta$ in $D^2$. The multiplicities of $\zeta$ give a domain $D \in D_{\beta}(\xs,\ys)$. The ambiguity in the construction of $\ell_{\xs,\ys}$ arises from the choice of how we connect $\xs \cap b_i$ and $\ys \cap b_i$ along $b_i$, and different choices give rise to domains that differ by multiples of $O_i$. We now describe a deterministic procedure to obtain a matched domain $D_{\xs,\ys}$ from $D$.

    Since the closure of $\beta$ is a knot, the permutation $\pi$ it induces on the indices of the punctures $z_0,\dots,z_n$ is transitive. So, $\pi$ is an $(n+1)$-cycle. For $i \in \bZ$, let $\pi_i := \pi^i(0)$ and write $Z_i := z_{\pi_i}$ and $W_i := w_{\pi_i}$. In particular, $W_0 = w_0$ and $Z_{n+1} = z_0$. We can visit these base points by traversing the arcs $a_i'$ backwards and the discs $O_i$ in the following order:
    \[
    W_0 \stackrel{a_0'}{\longrightarrow} Z_1 \stackrel{O_{\pi_1}}{\longrightarrow} W_1 \stackrel{a_{\pi_1}'}{\longrightarrow} Z_2 \stackrel{O_{\pi_2}} {\longrightarrow} W_2 \stackrel{a_{\pi_2}'}{\longrightarrow} \dots \stackrel{a_{\pi_{n-1}}'}{\longrightarrow} Z_n \stackrel{O_{\pi_n}}{\longrightarrow} W_n \stackrel{a_{\pi_n}'}{\longrightarrow} Z_{n+1} = Z_0.
    \]
    A domain $D_{\xs,\ys}$ is matched if and only if $n_{W_0}(D_{\xs,\ys}) = 0$ and $n_{W_i}(D_{\xs,\ys}) = n_{Z_{i+1}}(D_{\xs,\ys})$ for $i \in \{1,\dots,n\}$. As $Z_{n+1} = Z_0$ does not lie in any of the $O_i$, the multiplicity here is fixed, and we work backwards: We recursively construct a sequence of domains $D_0 := D$, $D_1, \dots, D_n$, such that $n_{W_0}(D_k) = 0$ and $n_{W_i}(D_k) = n_{Z_{i+1}}(D_k)$ for $i \in \{n-k+1,\dots,n\}$. In particular, $D_{\xs,\ys} := D_n$ is matched. Given $D_k$ for $k \in \{0,\dots,n-1\}$, let
    \[
    D_{k+1} := D_k + (n_{Z_{n-k+1}}(D_k) - n_{W_{n-k}}(D_k)) O_{\pi_{n-k}}.
    \]
    Then $n_p(D_{k+1}) = n_p(D_k)$ for $p \not\in O_{\pi_{n-k}}$ and
    \[
    n_{W_{n-k}}(D_{k+1}) = n_{W_{n-k}}(D_k) + (n_{Z_{n-k+1}}(D_k) - n_{W_{n-k}}(D_k)) = n_{Z_{n-k+1}}(D_k),
    \]
    so $D_{k+1}$ satisfies the required properties.

    Finally, we show uniqueness. Suppose that $D_{\xs,\ys}$ and $D_{\xs,\ys}' \in D_{\beta}(\xs,\ys)$ are both matched. Then $D_0 := D_{\xs,\ys} - D_{\xs,\ys}'$ is a matched periodic domain. So, $D_0 = \sum_{i=1}^n \lambda_i O_i$ by Lemma~\ref{lem:periodic}, and hence $n_{W_i}(D_0) = n_{Z_i}(D_0) = \lambda_{\pi_i}$ for $i \in \{1,\dots,n\}$. We have $n_{Z_{n+1}}(D_0) = n_{W_0}(D_0) = 0$, since $Z_{n+1} = Z_0$. The matched condition $n_{W_i}(D_0) = n_{Z_{i+1}}(D_0)$ for $i \in \{1,\dots,n\}$ now implies that $\lambda_{\pi_1} = \dots = \lambda_{\pi_n} = 0$, so $D_0 = 0$.
\end{proof}

We identify the intersection points in the disc model and the Heegaard diagram.

\begin{prop}
   Let $\beta \in B_{n+1}$ be a braid, $\bD_\beta = (\cD_n, \ba', \bb)$ the corresponding disc model, and $\cH_\beta = (\Sigma, \as, \bs, w, z)$ the associated doubly pointed Heegaard diagram. The natural embedding of $D^2 \setminus N(\ws \cup \zs)$ into $\Sigma$ induces a bijection $b_\beta \colon \cI_{\bD_\beta} \to \cI_{\cH_\beta}$.
\end{prop}

\begin{proof}
When we close up the arcs $a_1', \dots, a_n'$ along the handles added to $D^2$ to obtain $\Sigma$, we do not create any new intersection points with $b_1, \dots, b_n$, since they lie in $D^2$.
\end{proof}

\begin{notation}
From now on, we will denote this set of intersection points as
\begin{equation}
    \cI :=\cI_{\bD_\beta} \cong \cI_{\cH_\beta}.
\end{equation}
\end{notation}

\begin{prop} \label{prop:domain-correspondence}
    Let $\xs$, $\ys \in \cI$. Then there is a unique domain $D^\Sigma_{\xs,\ys} \in D_{\cH}(\xs,\ys)$ such that $n_w(D^\Sigma_{\xs,\ys}) = 0$. It can be obtained from the matched domain $D_{\xs,\ys} \in D_{\beta}(\xs,\ys)$ by performing surgery along the 0-spheres $\{w_i, z_{\pi(i)}\}$ for $i \in \{1,\dots,n\}$. In other words, we attach a tube $T_i$ to $S^2 \setminus N(\ws \cup \zs)$ along $\partial N(\{w_i, z_{\pi(i)}\})$, where the multiplicity of $D^\Sigma_{\xs,\ys}$ is $n_{w_i}(D_{\xs,\ys}) = n_{z_{\pi(i)}}(D_{\xs,\ys})$.
\end{prop}

\begin{proof}
    Given the matched domain $D_{\xs,\ys}$, attaching tubes results in a domain in $D_{\cH}(\xs,\ys)$ such that $n_w(D^\Sigma_{\xs,\ys}) = 0$. Conversely, any domain $D^\Sigma_{\xs,\ys} \in D_{\cH}(\xs,\ys)$ satisfying $n_w(D^\Sigma_{\xs,\ys}) = 0$ must have constant multiplicity along each $T_i$. Indeed, it contains only a single arc of $\alpha_i$ that intersects both components of $\partial T_i$, but no other curves from $\as \cup \bs$, and so $T_i \setminus (\as \cup \bs)$ is connected. Surgering $D^\Sigma_{\xs,\ys}$ along the cores of the $T_i$ gives a matched domain in $D_{\beta}(\xs,\ys)$.
\end{proof}

\begin{rmk}\label{rmk:periodic}
    Uniqueness also follows from Theorem~\ref{THD}, since $(\Sigma,\as,\bs, w)$ is a pointed Heegaard diagram for $S^3$ and periodic domains $P$ satisfying $n_w(P) = 0$ correspond to $H_2(S^3) = 0$. Hence, if $P$ is a periodic domain, then $P_0 := P - n_w(P)[\Sigma]$ is a periodic domain with $n_w(P_0) = 0$, so $P_0 = 0$. Consequently, if $D$, $D' \in D_{\cH}(\xs,\ys)$, then $D-D' = (n_w(D)-n_w(D'))[\Sigma]$.
\end{rmk}

\subsection{The classical Alexander grading}\label{sec:gradings}
In Definition~\ref{def:grdisc}, we assigned an Alexander grading $A^{\Pd}(\xs)$ to an intersection point $\xs \in \cI$. We always have a canonical intersection point $\xs^0 \in \cI$, as in Definition~\ref{def:canonical-intersection-point}.

\begin{lem}\label{lem:canonical-grading}
    We have 
    \[
    A^{\Pd}(\xs^0)= \frac{w(\beta)+n}{2} \text{ and } A^{\qPd}(\xs^0) = n + w(\beta).
    \]
\end{lem}

\begin{proof}
    By definition, 
    \[
    A^{\Pd}(\xs_0) = \deg_x \varphi(l_{\xs_0}) + \frac{w(\beta)+n}{2} \text{ and } A^{\qPd}(\xs_0) = \deg_d \phi(l_{\xs_0}) + w(\beta).
    \]
    As the braid action on the right-hand side of $\cD_n$ is trivial, $l_{\xs_0}$ is the product of loops $l_1,\dots,l_n$ on $\cD_n$, such that $l_i$ is obtained by following $\eta_i'$ from $d_i$ to $\eta_i'(1)$, then $\eta_i'(1)$ to $x_i^0$ along a short arc on $b_i$, and finally $a_i'$ back to $d_i$. The arcs $l_i$ are not nested, so there is no relative winding. Each $l_i$ goes around a single point of $\qs$ anticlockwise once and does not wind around any point in $\zs \cup \ws$. Hence, $\deg_x \varphi(l_{\xs_0}) = \deg_x \phi(l_{\xs_0}) = 0$. Furthermore, $\deg_d \phi(l_{\xs_0}) = n$.
\end{proof}

The \emph{classical Alexander grading} $A^{\HF}$ arising from knot Floer homology was defined by Ozsv\'ath and Szab\'o~\cite{OS2}. The induced relative grading is easy to describe. 

\begin{defn}\label{def:Alex-HF}
Let $(\Sigma,\as,\bs, w, z)$ be a doubly pointed Heegaard diagram for a knot $K$ in $S^3$. For $\xs$, $\ys \in \cI$, let $D \in D_{\cH}(\xs, \ys)$. We define
\begin{equation}
A^{\HF}(D):= n_{z}(D)-n_{w}(D).
\end{equation}
\end{defn}

This is independent of the choice of domain $D$ connecting $\xs$ to $\ys$, since any two domains in $D_\cH(\xs,\ys)$ differ by a multiple of $[\Sigma]$ by Remark~\ref{rmk:periodic}, and $n_z([\Sigma]) = n_w([\Sigma])$.

The absolute classical Alexander grading of $\xs \in \cI$ is
\[
A^{\HF}(\xs) = A^{\HF}(\xs^0) + A^{\HF}(D),
\]
where $D \in D_{\cH}(\xs, \xs^0)$. We will show in Corollary~\ref{coro:canonical-HF-grading} that
\[
A^{\HF}(\xs^0) = \frac{w(\beta)+n}{2}. 
\]

\section{Knot Floer homology}\label{sec:HFK}

Let $K$ be a knot in $S^3$, and $\cH = (\Sigma,\as,\bs,w,z)$ a doubly pointed Heegaard diagram for $K$. Knot Floer homology $\HFKh(K)$ is a finitely generated $\bF_2$-module, bigraded by the Alexander grading $A^{\HF}$ (Definition~\ref{def:Alex-HF}) and the Maslov (homological) grading $M$. 

Fix an almost complex structure $J$ on $\Sym^n(\Sigma)$ satisfying the conditions described by Ozsv\'ath and Szab\'o~\cite{OS2, OS}. For $u \in \pi_2(\xs,\ys)$, let $\cM(u)$ be the moduli space of pseudo-holomorphic representatives of $u$. For generic $J$, the space $\cM(u)$ is a smooth manifold of dimension $\mu(u)$, where $\mu(u)$ is the Maslov index of $u$. Under the usual admissibility and transversality assumptions, the codimension-one boundary strata in the compactification of the relevant one-dimensional moduli spaces consist of broken flow lines. A combinatorial formula for the Maslov index was given by Lipshitz~\cite{Lipshitz}; see equation~\eqref{eqn:Maslov}. 
As $D^2 \setminus \{i,-i\}$ is conformally equivalent to $\bR \times [0,1]$, there is an $\bR$-action on $\cM(u)$, and we write $\cMh(u) := \cM(u)/\bR$. When $\mu(u) = 1$, the space $\cMh(u)$ is a compact 0-dimensional manifold.

The chain group $\CFKh(\mathscr H)$ is the $\bF_2$-vector space generated by $\cI$. For a generator $\xs \in \cI$, we define
\begin{equation}\label{eqn:HFK-differential}
\ph \xs := \sum_{\ys \in \cI} \sum_{\substack{u \in \pi_2(\xs,\ys):\\ \mu(u)=1 \\ n_z(u) = n_w(u) = 0}} \# \cMh(u) \cdot \ys,
\end{equation}
where $\# \cMh(u) \in \bF_2$ is the number of elements of $\cMh(u)$ modulo two. The sum is finite since any diagram of a knot in $S^3$ is weakly admissible. Note that 
\[
\{u \in \pi_2(\xs,\ys) : n_z(u) = n_w(u) = 0\} \neq \emptyset 
\]
if and only if $A^{\HF}(\xs) = A^{\HF}(\ys)$. Ozsv\'ath--Szab\'o and Rasmussen showed that $\ph \circ \ph = 0$. Hence, $\ph$ is a differential, and we can define
\[
\HFKh(K) := H_*\left(\CFKh(\cH), \ph \right).
\]

The Maslov grading $M \colon I \rightarrow \mathbb Z$ is the standard absolute Maslov grading obtained from the singly pointed diagram of $S^3$ after forgetting $z$. It is characterised by the following properties:
\begin{itemize}
\item For $\xs$, $\ys \in \cI$, and $u \in \pi_2(\xs,\ys)$, we have 
\[
M(\xs)-M(\ys)=\mu(u)-2 n_{w}(u).
\]
\item If $\xs \in \cI$ is a cycle in the singly pointed complex $\CFh(\Sigma,\as,\bs,w)$ and represents a generator of $\HFh(S^3) \cong \bF_2$, then $M(\xs) = 0$. 
\end{itemize}
The differential $\ph$ preserves the Alexander grading and decreases the Maslov grading by one. The pair $(A^{\HF}, M)$ is a bigrading on the set of intersection points $\cI$. 

\begin{thm}[Ozsv\'ath--Szab\'o~\cite{OS2}, Rasmussen~\cite{Rass}]\label{thm:HFK}
Let $K$ be a knot in $S^3$ and $\mathcal H$ a doubly pointed Heegaard diagram for $K$. Then, the bigrading on the complex descends to a bigrading on the homology:
\[
\HFKh(\mathcal H) = \bigoplus_{(a, m) \in \bZ^2} \HFKh_{a, m}(\mathcal H)
\]
that is an invariant of $K$ up to bigrading-preserving isomorphism.
This is a categorification of the symmetrised Alexander polynomial $\Delta_K(x)$ of $K$ in the sense that
\begin{equation}
\Delta_K(x) = \sum_{(a, m) \in \bZ^2} (-1)^{m} \rank \left( \HFKh_{a, m}(K)\right) \cdot x^{a}.
\end{equation}
\end{thm}

In light of Theorem~\ref{thm:HFK}, we write $\HFKh(K)$ for the isomorphism class of the bigraded group $\HFKh(\mathcal H)$.

\section{The Alexander gradings in the disc model and the Heegaard diagram agree}\label{S:identif}

Let $\beta$ be an $(n+1)$-braid with closure $K := \hat{\beta}$, and consider the chain complex $\CFKh(\mathscr H_{\beta})$ associated with the doubly pointed Heegaard diagram $\mathscr H_{\beta}$ from Theorem~\ref{THD}. Let 
\[
\HFKh(K) = \HFKh(\mathscr H_{\beta}) = H_*\left(\CFKh(\cH_\beta)\right)
\]
be the corresponding knot Floer homology. We denote by $(A^{\HF},M)$ the bigrading defined as in Section~\ref{sec:HFK}, given by the Alexander grading and the Maslov grading. We restate and prove Theorem~\ref{thm:Alexander-agree} from the introduction. We also simultaneously show Proposition~\ref{prop:signs}.

\thmidgr*

\begin{proof}[Proof of Theorem~\ref{thm:Alexander-agree} and Proposition~\ref{prop:signs}.] First, we show that for intersection points $\xs$, $\ys \in \cI$, the relative gradings agree:
\begin{equation}\label{eqn:Alexander-agree}
A^{\Pd}(\xs)-A^{\Pd}(\ys) = A^{\HF}(\xs) - A^{\HF}(\ys).
\end{equation}
By definition,
\[
A^{\HF}(\xs)-A^{\HF}(\ys) := A^{\HF}(D) = n_z(D) - n_w(D),
\]
where $D \in D_\cH(\xs, \ys)$.
We show that this agrees with the local system grading. Recall that
\begin{equation}
A^{\Pd}(\xs) - A^{\Pd}(\ys)= \deg_x \varphi (l_{\xs}) - \deg_x \varphi(l_{\ys}) = \deg_x \varphi(l_{\xs}l_{\ys}^{-1}),
\end{equation}
where the loops $l_{\xs}$ and $l_{\ys}$ in the configuration space are constructed using paths in $\cD_n$ that follow the paths $\eta_1,\dots,\eta_n$ or $\eta_1',\dots,\eta_n'$, the simple closed curves $b_1, \dots, b_n$, and the arcs $a_1',\dots,a_n'$.

For each $i \in \{1,\dots,n\}$, there are four possibilities:
\begin{enumerate}
    \item\label{it:1} both $l_{\xs}$ and $l_{\ys}$ use $\eta_i$,
    \item\label{it:2} both $l_{\xs}$ and $l_{\ys}$ use $\eta_i'$,
    \item\label{it:3} $l_{\xs}$ uses $\eta_i'$ and $l_{\ys}$ uses $\eta_i$,
    \item\label{it:4} $l_{\xs}$ uses $\eta_i$ and $l_{\ys}$ uses $\eta_i'$.
\end{enumerate}

We modify the loop $l_{\xs} l_{\ys}^{-1}$ to a new loop $l_{\xs,\ys}$ such that $l_{\xs,\ys}$ is supported on $\ba' \cup \bb$, as follows. In cases~\ref{it:1} and~\ref{it:2}, we have the segment $\eta_i \eta_i^{-1}$ or $\eta_i' (\eta_i')^{-1}$ in $l_{\xs} l_{\ys}^{-1}$, which we remove when constructing $l_{\xs,\ys}$. In case~\ref{it:3}, we have the segment $(\eta_i)^{-1} \eta_i'$ in $l_{\xs} l_{\ys}^{-1}$. Let $\eta''_i$ be a path connecting $\eta_i(1)$ and $\eta'_i(1)$ along the bottom of $b_i$, as in Figure~\ref{Ex2}. We replace $(\eta_i)^{-1} \eta_i'$ with $\eta_i''$ in the construction of $l_{\xs,\ys}$. Finally, in case~\ref{it:4}, we replace $(\eta_i')^{-1} \eta_i$ with $(\eta_i'')^{-1}$ when constructing $l_{\xs,\ys}$.

\begin{figure}
\centering
\def\svgwidth{0.9\textwidth}
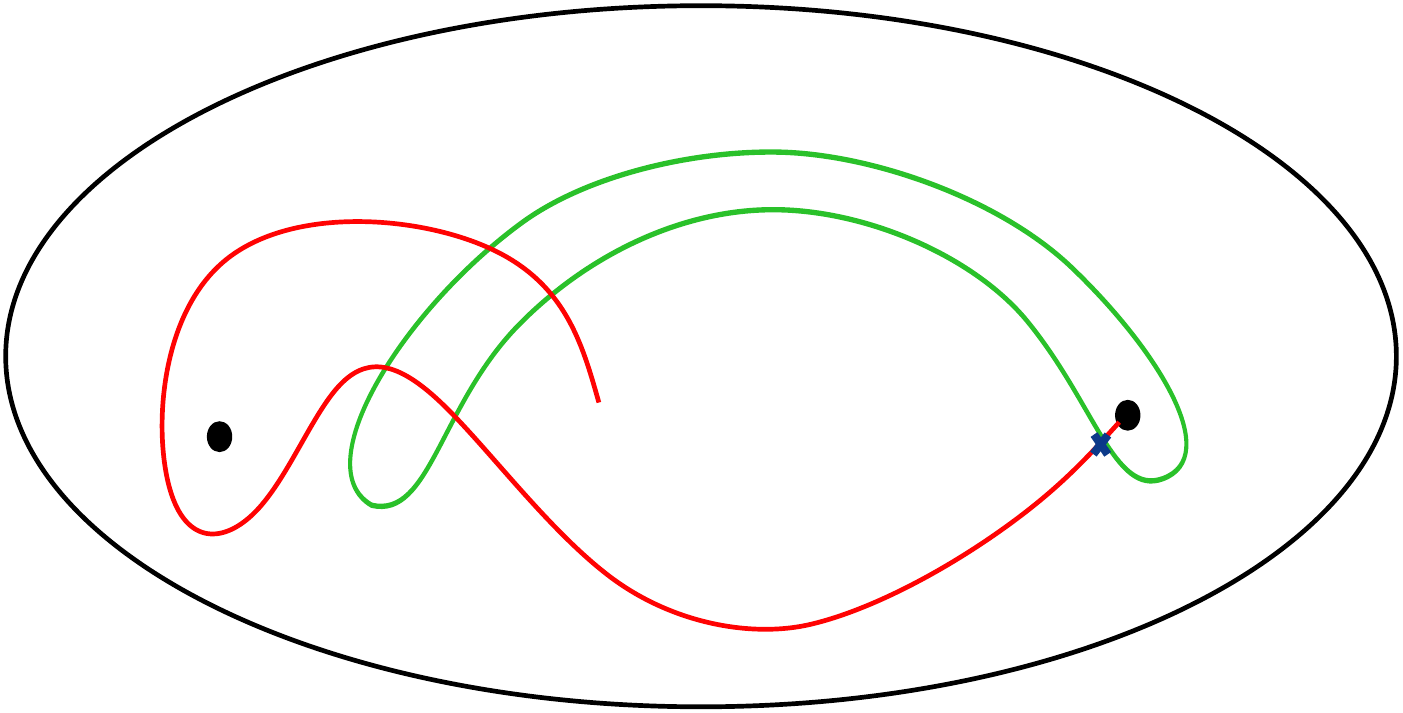
\caption{Changing the loop $l_{\xs}l_{\ys}^{-1}$ to $l_{\xs,\ys}$.}
\label{Ex2}
\end{figure} 
   
Note that the loop $\eta''_i(\eta_i')^{-1} \eta_i$ has the same winding number around $\zs$ and $\ws$. As $\varphi$ counts winding around points of $\zs$ with $x$ and around points of $\ws$ with $x^{-1}$, we have
\begin{equation}
  \deg_x \varphi\left(l_{\xs} l_{\ys}^{-1} \right)= \deg_x \varphi(l_{\xs,\ys}).
\end{equation}
Hence,
\begin{equation}\label{pdgrad}
A^{\Pd}(\xs)-A^{\Pd}(\ys)= \deg_x \varphi(l_{\xs,\ys}).
\end{equation}

Given a domain $D \in D_{\beta}(\xs,\ys)$, let 
\[
n_{\zs}(D) := \sum_{i=0}^n n_{z_i}(D) \text{ and } n_{\ws}(D) := \sum_{i=0}^n n_{w_i}(D).
\]

\begin{prop}
If $D \in D_{\beta}(\xs,\ys)$, then
\begin{equation}\label{r1}
  \deg_x \varphi\left( l_{\xs,\ys} \right)=n_{\zs}(D) - n_{\ws}(D).
\end{equation}
\end{prop}

\begin{proof}
If $\partial D = l_{\xs,\ys}$, then this follows from the fact that the local system $\varphi$ counts monodromies around the points of $\zs$ and $\ws$ with values $x$ and $x^{-1}$, respectively. If $D' \in D_{\beta}(\xs,\ys)$ is another domain, then $P := D-D'$ is a periodic domain that is a $\bZ$-linear combination of $O_1,\dots,O_n$ by Lemma~\ref{lem:periodic}. As each $O_i$ contains exactly one point from $\zs$ and one from $\ws$, we have $n_{\zs}(P) - n_{\ws}(P) = 0$; so
\[
n_{\zs}(D) - n_{\ws}(D) = n_{\zs}(D') - n_{\ws}(D'). \qedhere
\]
\end{proof}

\begin{prop}\label{chdom}
    Let $\xs$, $\ys \in \cI$. If $D_{\xs,\ys} \in D_{\beta}(\xs,\ys)$ is the matched domain and $D^{\Sigma}_{\xs,\ys} \in D_{\cH}(\xs,\ys)$ is the unique domain satisfying $n_w(D^{\Sigma}_{\xs,\ys}) = 0$, then
    \begin{equation}\label{r2}
        n_{\zs}(D_{\xs,\ys})-n_{\ws}(D_{\xs,\ys}) = n_z \left(D^{\Sigma}_{\xs,\ys}\right)-n_w \left( D^{\Sigma}_{\xs,\ys} \right).
    \end{equation}
\end{prop}

\begin{proof}
    By Proposition~\ref{prop:domain-correspondence}, we can obtain $D_{\xs,\ys}^\Sigma$ from $D_{\xs,\ys}$ by performing surgery along the 0-spheres $\{w_i,z_{\pi(i)}\}$ for $i \in \{1,\dots,n\}$, and the matched condition states that $n_{w_i}(D_{\xs,\ys}) = n_{z_{\pi(i)}}(D_{\xs,\ys})$. Hence,
    \[
    n_{\zs}(D_{\xs,\ys})-n_{\ws}(D_{\xs,\ys}) = n_{z_{\pi(0)}}(D_{\xs,\ys}) - n_{w_0}(D_{\xs,\ys}) = n_z\left(D^{\Sigma}_{\xs,\ys}\right)-n_w \left( D^{\Sigma}_{\xs,\ys} \right). \qedhere
    \]
\end{proof}

We use the domain $D^{\Sigma}_{\xs,\ys}$ in $\mathcal H_{\beta}$ to compute the relative Alexander grading:
\[
A^{\HF}(\xs)-A^{\HF}(\ys) = n_{z}(D^{\Sigma}_{\xs,\ys})-n_{w}(D^{\Sigma}_{\xs,\ys}).
\]
Combining this relation with equations \eqref{pdgrad}, \eqref{r1}, and \eqref{r2}, we conclude that
\begin{equation}
A^{\Pd}(\xs)-A^{\Pd}(\ys)=A^{\HF}(\xs)-A^{\HF}(\ys).
\end{equation}
Hence, there is a constant $c$ such that 
\begin{equation}\label{eqn:Alexander-shifted}
A^{\Pd}(\xs) = A^{\HF}(\xs) + c
\end{equation}
for every $\xs \in \cI$. We now prove a weaker version of Proposition~\ref{prop:signs}.

\begin{lem}\label{lem:signs-weak} 
Let $\beta$ be a braid whose closure is a knot. For $\xs \in \cI$, let $\varepsilon_{\xs}$ be the intersection sign of $(\beta \cup \bI)\cs$ and $\ct$ at $\xs$. Then there is $\varepsilon \in \{1,-1\}$ such that
\[
\varepsilon_{\xs} = (-1)^{M(\xs)} \varepsilon
\]
for every $\xs \in \cI$.
\end{lem}

\begin{proof}
    For an intersection point $\xs \in \cI$, the intersection sign $\varepsilon_{\xs}$ in $(\beta \cup \bI)\cs \cap \ct$ tautologically agrees with the intersection sign in $\bT_{\as} \cap \bT_{\bs}$. In Heegaard Floer homology, we have $\varepsilon_{\xs} = (-1)^{M(\xs)}$, up to an overall sign, where $M(\xs)$ is the homological grading of $\xs$; see Lipshitz--Ozsv\'ath--Thurston~\cite[Section~4.1]{bordered-notes} and Ozsv\'ath--Szab\'o~\cite[Section~1.2]{HF-lectures}.
\end{proof}

As knot Floer homology categorifies the Alexander polynomial (Theorem~\ref{thm:HFK}), it follows from Theorem~\ref{thm:Alexander-disc}, Lemma~\ref{lem:signs-weak}, and equation~\eqref{eqn:Alexander-shifted} that 
\[
(-1)^{\frac{n(n+1)}{2}}\Delta_K(x) = \varepsilon x^c \Delta_K(x).
\]
Since $\Delta_K(x)$ is a non-zero Laurent polynomial, $c = 0$. It follows that $A^{\Pd}(\xs) = A^{\HF}(\xs)$ for every $\xs \in \cI$ and $\varepsilon = (-1)^{\frac{n(n+1)}{2}}$. This concludes the proofs of Theorem~\ref{thm:Alexander-agree} and Proposition~\ref{prop:signs}.
\end{proof}
 
\begin{coro}\label{coro:canonical-HF-grading}
Let $\beta \in B_{n+1}$ be a braid whose closure is a knot, and let $\xs^0 \in \cI_{\cH_\beta}$ be the canonical intersection point. Then    
\[
A^{\HF}(\xs^0) = \frac{w(\beta)+n}{2}. 
\]
\end{coro}

\begin{proof}
    This follows from Theorem~\ref{thm:Alexander-agree} and Lemma~\ref{lem:canonical-grading}. 
\end{proof}

\section{The Jones polynomial via quantum Heegaard diagrams}\label{sec:unification} 

In Section~\ref{S:unifmodel}, we presented a unified model for the Alexander and Jones polynomials defined via a two-variable local system $\phi$ on the configuration space $C^q_n$ of the punctured disc $\cD^q_n = D^2 \setminus (\zs \cup \ws \cup \qs)$. Given a braid $\beta \in B_{n+1}$, we defined a two-variable polynomial $\Omega'_\beta(x,d)$ such that the Alexander polynomial $\Delta_K(x) = \Omega'_\beta(x,-1)$ and the Jones polynomial $V_{\mK}(x) = \Omega'_\beta(x,x^{-1})$; see Theorem~\ref{Topmod}.

In Definition~\ref{def:quantum-Heegaard-diagram}, we associated a quantum Heegaard diagram $\cH_\beta^q$ to $\beta$. In this section, we translate Theorem~\ref{Topmod} from $\bD_\beta^q$ to $\cH_\beta^q$. We show that the quantum intersection polynomial $\Omega^q_\beta$ from Definition~\ref{def:quantum-intersection} corresponds to $\Omega'_\beta$ and hence recovers the Alexander and Jones polynomials.

\begin{defn}
Let $\beta \in B_{n+1}$ be a braid and 
\[
\cH_\beta^q = (\Sigma,\as,\bs, w, z, \qs, \hs, \xs^0) 
\]
the associated quantum Heegaard diagram. Consider the \emph{fat diagonal}
\[
\Delta :=\{ (x_1,\dots,x_n) \in \Sym^n(\Sigma) \mid \exists \ 1\leq i < j \leq n : x_i=x_j \}
\]
of $\Sym^n(\Sigma)$. For $n=1$, we have $\Delta = \emptyset$. Let $\xs$, $\ys \in \bT_{\as} \cap \bT_{\bs}$, and $u \in \pi_2(\xs,\ys)$ be a Whitney disc. Then let
\[
n_\Delta(u) := u \cdot \Delta
\]
be the algebraic intersection number between $u$ and $\Delta$. (For $n=1$, we have $n_\Delta(u)=0$.) Furthermore, we write
\[
n_{\qs}(u) := \sum_{i=0}^n n_{q_i}(u).
\]
For $i \in \{1,\dots,n\}$, let $h_i \in \Sigma$ be a base point on the handle that is attached at $a_i'(0) = z_{\pi(i)}$ and $a_i'(1) = w_i$, and write $\hs := \{h_1, \dots, h_n\}$. We let
\[
n_{\hs}(u) := \sum_{i=1}^n n_{h_i}(u).
\]
For a domain $D \in D(\cH_\beta^q)$, we define $n_{\qs}(D)$ and $n_{\hs}(D)$ as the sum of the multiplicities of $D$ at $\qs$ and $\hs$, respectively.
\end{defn}

According to Rasmussen~\cite[Theorem~9.1]{Rass} and Lipshitz~\cite[Corollary~4.10]{Lipshitz}, for $\xs$, $\ys \in \cI$, and $u \in \pi_2(\xs,\ys)$, we have
\begin{equation}\label{eqn:diagonal}
n_\Delta(u) = n_{\xs}(\cD(u)) + n_{\ys}(\cD(u)) - e(\cD(u)),
\end{equation}
where the \emph{point measures} $n_{\xs}(\cD(u))$ and  $n_{\ys}(\cD(u))$ are the sum of the averages of the multiplicities of $\cD(u)$ at the points of $\xs$ and $\ys$, respectively, and $e(\cD(u))$ is the \emph{Euler measure} of $\cD(u)$. Furthermore, the Maslov index of $u$ can be computed as
\begin{equation}\label{eqn:Maslov}
\mu(u) = n_{\xs}(\cD(u)) + n_{\ys}(\cD(u)) + e(\cD(u)).
\end{equation}

\begin{defn}\label{def:diagonal-intersection}
Let $\cH = (\Sigma,\as,\bs)$ be a Heegaard diagram and $\xs$, $\ys \in \bT_{\as} \cap \bT_{\bs}$. For a domain $D \in D_{\cH}(\xs,\ys)$, we define 
\[
n_\Delta(D) := n_{\xs}(D) + n_{\ys}(D) - e(D).
\]
\end{defn}

For $\xs$, $\ys \in \cI_{\bD_\beta}$ and $D \in D_\beta(\xs,\ys)$, we define $n_\Delta(D)$ analogously. If $D = \cD(u)$ for a topological Whitney disc, then $n_\Delta(u) = u \cdot \Delta$.

\begin{defn}\label{def:W}
Let $(\Sigma,\as,\bs, w, z, \qs, \hs, \xs^0)$ be a quantum Heegaard diagram and 
\[
\xs = (x_1,\dots,x_n) \in \bT_{\as} \cap \bT_{\bs}
\]
an intersection point such that $x_i \in \beta_i$ for $i \in \{1,\dots,n\}$. We write
\[
W_i(\xs) := 
\begin{cases}
    3(n-i)+1 &\text{ if } x_i = x_i^0, \\
    0 & \text{ otherwise,}
\end{cases}
\]
and let $W(\xs) := \sum_{i=1}^n W_i(\xs)$ be the \emph{weight} of $\xs$.
\end{defn}

Given the Heegaard diagram $\mathscr H_{\beta}$, we defined the classical Alexander grading $A^{\HF}$ in Definition~\ref{def:Alex-HF} via
\[
A^{\HF}(D):=n_{z}(D)-n_{w}(D).
\]
This naturally extends to quantum Heegaard diagrams. 

\begin{lem}\label{lem:additivity}
    Let $\cH^q_\beta$ be a quantum Heegaard diagram and $\xs$, $\ys$, $\zs \in \cI$ intersection points. Given $D \in D_\cH(\xs,\ys)$ and $D' \in D_\cH(\ys, \zs)$, we have
    \[
    \begin{split}
    n_\Delta(D+D') &= n_\Delta(D) + n_\Delta(D'),\\
    n_{\hs}(D+D') &= n_{\hs}(D) + n_{\hs}(D'),\\
    n_{\qs}(D+D') &= n_{\qs}(D) + n_{\qs}(D')
    \end{split}
    \]
    for $D+D' \in D_\cH(\xs,\zs)$.
\end{lem}

\begin{proof}
    The additivity of $n_{\hs}$ and $n_{\qs}$ is straightforward. We now show the additivity of $n_\Delta$. For $n=1$, this is straightforward as $n_\Delta \equiv 0$. For $n \ge 2$, there exist $u \in \pi_2(\xs,\ys)$ and $u' \in \pi_2(\ys,\zs)$ such that $\cD(u) = D$ and $\cD(u') = D'$. Then $D+D' = \cD(u u')$, where $uu'$ is the concatenation of $u$ and $u'$. The result then follows from
    \[
    (u u') \cdot \Delta = u \cdot \Delta + u' \cdot \Delta
    \]
    and equation~\eqref{eqn:diagonal}.
\end{proof}

We now define a new grading using the additional base points $\qs$ and $\hs$, and the distinguished intersection point $\xs^0$. We will use this to provide a formula for the Jones polynomial in terms of $\cH^q_\beta$.

\begin{defn}\label{def:quantum-grading-disc}
Given a quantum Heegaard diagram $\mathcal H^q_{\beta}$, intersection points $\xs$, $\ys \in \cI$, and a domain $D \in D_{\cH}(\xs, \ys)$ such that $n_w(D)=0$, we define
\[
A^{\qHF}(D):= n_\Delta(D)-2n_{\hs}(D) + n_{\qs}(D) + W(\xs) - W(\ys).
\]
\end{defn}

The following is motivated by Lemma~\ref{lem:canonical-grading}.

\begin{defn}\label{def:canonical-quantum-grading}
We fix the quantum Alexander grading of the canonical intersection point $\xs^0 \in \cI$ to be 
\begin{equation}\label{eqgr0}
A^{\qHF}(\xs^0) := n + w(\beta).
\end{equation}
\end{defn}

Since $(\Sigma, \as, \bs, w)$ is a Heegaard diagram for $\mathbb S^3$, there is a unique domain $D$ with $n_w(D)=0$ connecting any two intersection points, and we can define absolute gradings:

\begin{defn}\label{qAlex}
Let $\xs\in \cI$ be an intersection point and $D \in D_{\cH^q_\beta}(\xs, \xs^0)$ the unique domain such that $n_w(D) = 0$. We define the \emph{quantum Alexander grading} of $\xs$ as
\[
A^{\qHF}(\xs):=A^{\qHF}(\xs^0)+A^{\qHF}(D).
\]
\end{defn}

\begin{prop}\label{prop:transitivity}
    Let $\cH^q_\beta$ be a quantum Heegaard diagram and $\xs$, $\ys \in \cI$ be intersection points. If $D \in D_\cH(\xs,\ys)$ satisfies $n_w(D) = 0$, then
    \[
    A^{\qHF}(\xs) - A^{\qHF}(\ys) = A^{\qHF}(D).
    \]
\end{prop}

\begin{proof}
    By Proposition~\ref{prop:domain-correspondence}, there is a unique domain connecting any pair of intersection points with multiplicity zero at $w$. Hence,
    \[
    \begin{split}
    A^{\qHF}(\xs) - A^{\qHF}(\ys) &= \left(A^{\qHF}(\xs) - A^{\qHF}(\xs_0)\right) - \left(A^{\qHF}(\ys) - A^{\qHF}(\xs_0)\right) \\
    &= A^{\qHF}\left(D^\Sigma_{\xs,\xs_0}\right) - A^{\qHF}\left(D^\Sigma_{\ys,\xs_0}\right) \\
    &= A^{\qHF}\left(D^\Sigma_{\xs,\xs_0} - D^\Sigma_{\ys,\xs_0}\right) \\
    &= A^{\qHF}\left(D^\Sigma_{\xs,\ys}\right),
    \end{split}
    \]
    where the second equality follows from the definition of $A^{\qHF}$ of a domain and the third equality from the definition of $A^{\qHF}$ together with Lemma~\ref{lem:additivity}.
\end{proof}

We introduced the quantum intersection polynomial $\Omega^q_\beta(x,d)$ in Definition~\ref{def:quantum-intersection}. In Theorem~\ref{thm:Alexander-agree}, we showed that $A^{\Pd}(\xs) = A^{\HF}(\xs)$ for every $\xs \in \cI$. We now restate and prove the quantum version of this result.

\thmidqgr*

\begin{proof}
This proof is analogous to that of Theorem~\ref{thm:Alexander-agree}. However, as we use the local system $\phi$ instead of $\varphi$, we also have to deal with the relative winding and the marked points $\qs$. We first prove that the relative gradings agree. That is, for any $\xs$, $\ys \in \cI$, we have
\begin{equation}\label{eqn:relative-quantum-agree}
A^{qPd}(\xs)-A^{qPd}(\ys) = A^{\qHF}(\xs) - A^{\qHF}(\ys).
\end{equation}
    
{\it Step I: Replacing the loop $l_{\xs} l_{\ys}^{-1}$ with $l_{\xs,\ys}$.}
We recall that the loops $l_{\xs}$ and $l_{\ys}$ are obtained by following the paths $\eta_1,\dots,\eta_n$ or $\eta_1',\dots,\eta_n'$, then the curves $b_1, \dots, b_n$, and finally the arcs $a_1',\dots,a_n'$; see Definition~\ref{def:loop}.
As in the proof of Theorem~\ref{thm:Alexander-agree}, for each $i \in \{1,\dots,n\}$, there are four possibilities:
\begin{enumerate}
    \item\label{it:1q} both $l_{\xs}$ and $l_{\ys}$ use $\eta_i$; i.e., $x_i \neq x_i^0$ and $y_i \neq x_i^0$;
    \item\label{it:2q} both $l_{\xs}$ and $l_{\ys}$ use $\eta_i'$; i.e., $x_i = x_i^0$ and $y_i = x_i^0$;
    \item\label{it:3q} $l_{\xs}$ uses $\eta_i'$ and $l_{\ys}$ uses $\eta_i$; i.e., $x_i = x_i^0$ and $y_i \neq x_i^0$;
    \item\label{it:4q} $l_{\xs}$ uses $\eta_i$ and $l_{\ys}$ uses $\eta_i'$; i.e., $x_i \neq x_i^0$ and $y_i = x_i^0$.
\end{enumerate}
Here, we used the fact that the only point of $\ba' \cap \bb$ on the right-hand side of $b_i$ is $x_i^0$, since we obtained $\ba'$ from $\ba$ using a representative of the action of $\beta \cup \bI$ that is supported in the half-plane $\{x<0\}$. We let
\[
\omega_i(\xs,\ys) :=
\begin{cases}
    1 &\text{ if } x_i = x_i^0 \text{ and } y_i \neq x^0_i \\
    -1 &\text{ if } x_i \neq x_i^0 \text{ and } y_i = x^0_i, \\
    0 &\text{ otherwise.}
\end{cases}
\]

We replace $l_{\xs} l_{\ys}^{-1}$ with a loop $l_{\xs,\ys}$ supported on $\ba' \cup \bb$. The path $\eta''_i$ connects $\eta_i(1)$ and $\eta'_i(1)$ along the bottom of $b_i$, as shown in Figure~\ref{Ex2}. In cases~\ref{it:1q} and~\ref{it:2q}, we have the segment $\eta_i \eta_i^{-1}$ or $\eta_i' (\eta_i')^{-1}$ in $l_{\xs} l_{\ys}^{-1}$, which we remove when constructing $l_{\xs,\ys}$. In case~\ref{it:3q}, we replace $(\eta_i)^{-1} \eta_i'$ with $\eta_i''$ in the construction of $l_{\xs,\ys}$. In case~\ref{it:4q}, we replace $(\eta_i')^{-1} \eta_i$ with $(\eta_i'')^{-1}$ when constructing $l_{\xs,\ys}$.

We recall the definition of the quantum Alexander grading via the punctured disc from Definition~\ref{def:grdisc''}:
\begin{equation}
A^{\qPd}(\xs):=\deg_d \phi(l_{\xs}) + w(\beta).
\end{equation}
So, we have
\begin{equation}
A^{\qPd}(\xs)-A^{\qPd}(\ys)= \deg_d \phi(l_{\xs} l_{\ys}^{-1}).
\end{equation}
We introduce the three-variable local system 
\[
\psi := r \circ \rho \colon \pi_1(C^q_n) \to \bZ[x^{\pm 1}, y^{\pm 1}, d^{\pm 1}]^\times, 
\]
where
\[
r\left([\tilde{\zeta}_i]\right)= x \text{, } r\left([\tilde{\omega}_i]\right)= x^{-1} \text{, } r\left([\tilde{\gamma}_i]\right)= y \text{, and } r\left([\tilde{\delta}]\right)= d;
\]
cf.~Definition~\ref{localsystemg}. Then $\phi = \psi|_{y=-d}$, and hence
\[
\deg_d \phi = \deg_d \psi + \deg_y \psi.
\]

Consider the change in $\deg_d \psi$ when we replace the loop $l_{\xs} l_{\ys}^{-1}$ with $l_{\xs,\ys}$. There is no change in cases~\ref{it:1q} and~\ref{it:2q}. In case~\ref{it:3q}, 
\[
\left(\eta''_i(\eta_i')^{-1} \eta_i \right) \left((\eta_i)^{-1} \eta_i'\right) = \eta_i''.
\]
Hence, there is a change in $\deg_d \psi$ due to the fact that the closed loop $\eta''_i(\eta_i')^{-1} \eta_i$ bounds a disc $D$ containing precisely $n-i$ particles on the loop $l_{\xs} l_{\ys}^{-1}$. Then, the relative winding of the $i$th particle following $(\eta_i)^{-1} \eta_i'$ will change by a full clockwise rotation around each of the other $n-i$ particles belonging to $D$. Each full clockwise rotation contributes a $d^{-2}$ factor to $\psi$. Overall, the relative rotation of the loop will change by $n-i$ full clockwise rotations in the configuration space. In case~\ref{it:4q}, we add $n-i$ full anticlockwise rotations. So, we obtain
\begin{equation}
\deg_d \psi(l_{\xs} l_{\ys}^{-1})=\deg_d \psi(l_{\xs,\ys})+ \sum_{i=1}^{n} 2(n-i) \omega_i(\xs,\ys).
\end{equation}

Now consider the change in $\deg_y \psi$ when we replace the loop $l_{\xs} l_{\ys}^{-1}$ with $l_{\xs,\ys}$. In case~\ref{it:3q}, the loop $\eta''_i(\eta_i')^{-1} \eta_i$ that bounds $D$ contains $n-i+1$ $q$-punctures and rotates clockwise around each once. In case~\ref{it:4q}, the loop $(\eta_i'')^{-1}(\eta_i)^{-1} \eta_i'$ rotates around $n-i+1$ $q$-punctures anticlockwise once. Hence, 
\begin{equation}
\deg_y \psi(l_{\xs} l_{\ys}^{-1})=\deg_y \psi(l_{\xs,\ys})+\sum_{i=1}^n (n-i+1) \omega_i(\xs,\ys).
\end{equation}
This implies that
\begin{equation}\label{eqn:arc-change}
\begin{split}
A^{\qPd}(\xs)-A^{\qPd}(\ys) &= \deg_d \psi(l_{\xs,\ys}) + \deg_y \psi(l_{\xs,\ys}) + \sum_{i=1}^n (3(n-i)+1) \omega_i(\xs,\ys) \\
&= \deg_d \psi(l_{\xs,\ys}) + \deg_y \psi(l_{\xs,\ys}) + W(\xs) - W(\ys). 
\end{split}
\end{equation}

{\it Step II: From the loop $l_{\xs,\ys}$ to a domain on the quantum Heegaard diagram.} Let $D \in D_{\beta}(\xs,\ys)$ be a domain with $\partial D = l_{\xs,\ys}$ and $n_{w_0}(D) = 0$. Then
\begin{equation}\label{eqn:domain}
\deg_d \psi(l_{\xs,\ys}) = n_\Delta(D) \text{ and }
\deg_y \psi(l_{\xs,\ys}) = n_{\qs}(D),
\end{equation}
as $\deg_d \psi(l_{\xs,\ys})$ and $\deg_y \psi(l_{\xs,\ys})$ are the linking numbers of $l_{\xs,\ys}$ with the fat diagonal and with the divisors $V_{q_i}$, respectively. Since a Whitney disc $u \in \pi_2(\xs,\ys)$ with $\cD(u) = D$ bounds $l_{\xs,\ys}$, these linking numbers are $u \cdot \Delta = n_\Delta(D)$ and $\sum_{i=0}^n u \cdot V_{q_i} = n_{\qs}(D)$.

By Proposition~\ref{prop:matched}, there is a unique matched domain $D_{\xs,\ys} \in D_{\beta}(\xs,\ys)$. The difference $D_{\xs,\ys} - D$ is a periodic domain and, therefore, a linear combination of $O_1,\dots,O_n$ by Lemma~\ref{lem:periodic}. For $i \in \{1,\dots,n\}$, we have
\[
n_\Delta(O_i) = n_{\xs}(O_i) + n_{\ys}(O_i) - e(O_i) = \frac{1}{2} + \frac{1}{2} - 1 = 0,
\]
so $n_\Delta(D) = n_\Delta(D_{\xs,\ys})$. As $\qs \cap O_i = \emptyset$, we also have 
$n_{\qs}(D_{\xs,\ys}) = n_{\qs}(D)$.

By Proposition~\ref{prop:domain-correspondence}, we can obtain the unique domain $D^\Sigma_{\xs,\ys} \in D_{\cH}(\xs,\ys)$ satisfying $n_w(D^\Sigma_{\xs,\ys})=0$ from $D_{\xs,\ys}$ by performing surgery along the 0-spheres $\{w_i,z_{\pi(i)}\}$ for $i \in \{1,\dots,n\}$. This leaves $n_{\qs}$ unchanged as $\qs$ is disjoint from the tubes $T_1,\dots,T_n$ added, so
\begin{equation}\label{eqn:Q}
    n_{\qs}(D) = n_{\qs}(D^\Sigma_{\xs,\ys}).
\end{equation}
Furthermore,
\[
e(D_{\xs,\ys}) = e\left(D^\Sigma_{\xs,\ys}\right) + 2n_{\hs}\left(D^\Sigma_{\xs,\ys}\right),
\]
as surgery along $\{z_{\pi(i)}, w_i\}$ decreases the Euler measure by $2n_{h_i}(D^\Sigma_{\xs,\ys})$, where $h_i \in T_i$. In addition,
\[
n_{\xs}(D_{\xs,\ys}) = n_{\xs}\left(D^\Sigma_{\xs,\ys}\right) \text{ and } n_{\ys}(D_{\xs,\ys}) = n_{\ys}\left(D^\Sigma_{\xs,\ys}\right).
\]
It follows that
\begin{equation}\label{eqn:delta}
n_\Delta(D) = n_\Delta(D_{\xs,\ys}) = n_\Delta\left(D^\Sigma_{\xs,\ys}\right) - 2n_{\hs}\left(D^\Sigma_{\xs,\ys}\right).
\end{equation}
If we combine equation~\eqref{eqn:arc-change}, \eqref{eqn:domain}, \eqref{eqn:Q}, and~\eqref{eqn:delta} and Proposition~\ref{prop:transitivity}, we obtain that
\[
\begin{split}
A^{\qPd}(\xs)-A^{\qPd}(\ys) &= n_\Delta\left(D^\Sigma_{\xs,\ys}\right) - 2n_{\hs}\left(D^\Sigma_{\xs,\ys}\right) + n_{\qs}(D^\Sigma_{\xs,\ys}) + W(\xs) - W(\ys) \\ 
&= A^{\qHF}(D^\Sigma_{\xs,\ys}) \\ 
&= A^{\qHF}(\xs) - A^{\qHF}(\ys).
\end{split}
\]

Having shown that the relative gradings agree, we now turn our attention to the absolute gradings. By Definition~\ref{def:canonical-quantum-grading} and Lemma~\ref{lem:canonical-grading}, we have
\[
A^{\qHF}(\xs^0) := n+w(\beta) = A^{\qPd}(\xs^0).
\]
By substituting $\ys = \xs^0$ into equation~\eqref{eqn:relative-quantum-agree}, we obtain that
\[
A^{qPd}(\xs)-A^{qPd}(\xs^0) = A^{\qHF}(\xs) - A^{\qHF}(\xs^0),
\]
which concludes the proof of Theorem~\ref{thm:Jones-agree}.
\end{proof}

We are now ready to prove Theorem~\ref{qINT}, which we restate.

\thmAlexanderJones*

\begin{proof}
By Theorem~\ref{Topmod}, the two-variable polynomial $\Omega'_\beta(x,d)$ unifies the Alexander and Jones polynomials via
\begin{equation}
\begin{aligned}
&\Omega'_\beta(x,-1)=\Delta_K(x),\\
&\Omega'_\beta(x,x^{-1})=V_{\mK}(x). 
\end{aligned}
\end{equation}
Hence, the result follows once we show that
\begin{equation}\label{forms}
    \Omega'_\beta(x,-d) = \Omega^q_\beta(x,d).
\end{equation}
Since $\beta \in B_{n+1}$ and $\hat{\beta}$ is a knot, $w(\beta) \equiv n \pmod 2$. Hence, by Definition~\ref{defn}, we have
\begin{equation}\label{formulao}
\Omega'_\beta(x,-d) = (-1)^n x^{\frac{w(\beta)+n}{2}} d^{w(\beta)} \langle  (\beta \cup \mathbb I) \cs, \ct \rangle(x,-d) \in \bZ[x^{\pm 1}, d^{\pm 1}].
\end{equation}
By Definition~\ref{def:graded-intersection}, we have
\begin{equation}
\langle (\beta\cup \bI) \cs,\ct \rangle(x,-d) = \sum_{\xs = (x_1,\dots,x_n) \in \cI_{\bD_\beta}} \varepsilon_{x_1} \dots \varepsilon_{x_n} \cdot \phi(l_{\xs})(x,-d),
\end{equation}
where $\phi$ is the two-variable local system with values in $\bZ[x^{\pm 1},d^{\pm 1}]^\times$ from Definition~\ref{localsystemg}. Lemma~\ref{lem:signs} states that
\[
\varepsilon_{\xs} = (-1)^{\frac{n(n-1)}{2}} \sgn(\pi_{\xs}) \varepsilon_{x_1} \dots \varepsilon_{x_n}.
\]
It follows that
\[
\langle (\beta\cup \bI) \cs,\ct \rangle(x,-d) = (-1)^{\frac{n(n-1)}{2}} \sum_{\xs \in \cI_{\bD_\beta}} \varepsilon_{\xs} \cdot \sgn(\pi_{\xs}) \cdot \phi(l_{\xs})(x,-d).
\]
By Definition~\ref{def:grdisc''} and Theorem~\ref{thm:Alexander-agree}, we have
\begin{equation}
\deg_x \phi(l_{\xs}) + \frac{w(\beta)+n}{2} = A^{\Pd}(\xs) = A^{\HF}(\xs).
\end{equation}
Furthermore, by Definition~\ref{def:grdisc''} and Theorem~\ref{thm:Jones-agree}, we have
\begin{equation}
\deg_d \phi(l_{\xs}) + w(\beta) = A^{qPd}(\xs) = A^{\qHF}(\xs).
\end{equation}
According to Lemma~\ref{lem:delta}, we have
\[
\sgn(\varphi(l_{\xs})) = \sgn(\pi_{\xs}).
\]
However, $\varphi(l_{\xs}) = \phi(l_{\xs})(x,-1)$; hence
\[
\sgn(\phi(l_{\xs})(x,-d)) = \sgn(\phi(l_{\xs})(x,-1)) = \sgn(\pi_{\xs}).
\]
It follows that
\[
x^{\frac{w(\beta)+n}{2}} d^{w(\beta)} \phi(l_{\xs})(x,-d) = \sgn(\pi_{\xs}) \cdot x^{A^{\HF}(\xs)} d^{A^{\qHF}(\xs)}.
\]
We conclude that
\[
\Omega'_\beta(x,-d) = (-1)^{\frac{n(n+1)}{2}}\sum_{\xs \in \cI} \varepsilon_{\xs} \, x^{A^{\HF}(\xs)} d^{A^{\qHF}(\xs)} = \Omega^q_\beta(x,d),
\]
where the second equality is Definition~\ref{def:quantum-intersection}.
\end{proof}

\section{Quantum braid Floer homology}

We define a Floer-theoretic invariant of quantum Heegaard diagrams that categorifies the quantum intersection polynomial $\Omega^q_\beta(x,d)$ and admits a spectral sequence to knot Floer homology. Throughout this section, let $\beta \in B_{n+1}$ be a braid whose closure is a knot, and let
\[
\cH^q_\beta = (\Sigma,\as,\bs, w, z, \qs, \hs, \xs^0)
\]
be an associated quantum Heegaard diagram.

\begin{rmk}
    Neither $A^{\qHF}$ nor $-A^{\qHF}$ is a filtration on the knot Floer complex 
    \[
    \left(\CFKh(\Sigma,\as,\bs,w,z), \ph\right).
    \]
    In Figure~\ref{fig:trefoil-Alexander-Jones}, the shaded disc $u$ from $\xs$ to $\ys$ satisfies $A^{\qHF}(\ys) - A^{\qHF}(\xs) = 1$, as $n_\Delta(u) = n_{\hs}(u) = n_{\qs}(u) = 0$, $W(\ys) = 1$, and $W(\xs) = 0$. On the other hand, consider the rectangle $v$ from $\xs' := (x_3,z_7)$ to $\ys' := (x_5,z_5)$ in the diagram of $8_{19}$ in \cite[Figure~7.2]{Anghel2023TAMS}. This lies entirely in $\{x<0\}$, hence 
    \[
    n_{\hs}(v) = n_{\qs}(v) = W(\xs') = W(\ys') = 0. 
    \]
    As $n_\Delta(v) = 1$, we have $A^{\qHF}(\ys') - A^{\qHF}(\xs') = -1$. Since $\ph$ preserves $A^{\HF}$, it also follows that $A^{\HF} - A^{\qHF}$ is not a filtration, even though it is a natural candidate filtration for obtaining a categorification of the Jones polynomial.
\end{rmk}

\begin{defn}\label{def:P}
    For $\xs \in \cI$, let
    \[
    P(\xs) := A^{\qHF}(\xs) - W(\xs).
    \]
\end{defn}

\begin{defn}
    For $\xs$, $\ys \in \cI$, a domain $D \in D_{\cH}(\xs,\ys)$, and $u \in \pi_2(\xs,\ys)$, we set
    \[
    \kappa(D) := n_{\Delta}(D) - 2n_{\hs}(D) + n_{\qs}(D), 
    \]
    and write $\kappa(u) := \kappa(\cD(u))$.
\end{defn}

By Definition~\ref{def:quantum-grading-disc}, for $\xs$, $\ys \in \cI$, and $D \in D_\cH(\xs,\ys)$ such that $n_w(D) = 0$, we have
\begin{equation}\label{eqn:P}
    P(\xs) - P(\ys) = \kappa(D).
\end{equation}

\begin{defn}
    For $\xs \in \cI$, we let
    \[
    S(\xs) := \{i \in \{1,\dots,n\} : x_i = x^0_i\}.
    \]
\end{defn}

Let $q_{n+1}$ be a point in the region neighbouring the one containing $q_n$ along $a_n$. In the following result, we use the cylindrical reformulation of Heegaard Floer homology due to Lipshitz~\cite{Lipshitz}.

\begin{prop}\label{prop:positivity}
    Let $\beta \in B_{n+1}$ be a braid whose closure is a knot and $\cH^q_\beta$ the associated quantum Heegaard diagram. There is a non-empty class $\cJ(\cH^q_\beta)$ of generic almost complex structures on $\Sigma \times [0,1] \times \bR$ such that the following holds for every $J \in \cJ(\cH^q_\beta)$: If $\xs$, $\ys \in \cI$, and $D \in D_{\cH}(\xs,\ys)$ is such that $n_w(D) = 0$, and $D$ admits a $J$-holomorphic representative, then
    \[
    n_{\Delta}(D) - 2n_{\hs}(D) \ge 0, \quad n_{\qs}(D) \ge 0, \quad n_{q_{n+1}}(D) \ge 0.
    \]
    Furthermore, if $n_{\Delta}(D) - 2n_{\hs}(D) = 0$, then
    \[
    S(\xs) \subseteq S(\ys).
    \]
\end{prop}

\begin{proof}
    For $L \in \bR_{>0}$, let $\cJ_L(\cH^q_\beta)$ be a family of admissible almost complex structures on $\Sigma \times [0,1] \times \bR$ obtained by inserting standard necks of length $L$ along the tubes $T_1,\dots,T_n$ into a fixed $C^\infty$-precompact family of admissible almost complex structures. We require this family to be uniformly controlled away from the necks and fixed near $\{q_i\} \times [0,1] \times \bR$ for $i \in \{0,\dots,n+1\}$, and such that $\Delta$ and the divisors $V_{q_0},\dots,V_{q_{n+1}}$ are holomorphic with respect to the corresponding family of almost complex structures on $\Sym^n(\Sigma)$.

    We claim that there is $L_0 \in \bR_+$ such that, for every $L \ge L_0$, every generic $J \in \cJ_L(\cH^q_\beta)$, every $\xs$, $\ys \in \cI$, and $D \in D_{\cH}(\xs,\ys)$ with $n_w(D) = 0$ that admits a $J$-holomorphic representative, we have $n_{\Delta}(D) - 2n_{\hs}(D) \ge 0$ and $S(\xs) \subseteq S(\ys)$ whenever $n_{\Delta}(D) - 2n_{\hs}(D) = 0$. Otherwise, there is a sequence $L_k \to \infty$ for $k \in \bN$, almost complex structures $J_k \in \cJ_{L_k}(\cH^q_\beta)$, and $J_k$-holomorphic representatives $f_k$ of domains $D_k$ from $\xs_k$ to $\ys_k$, with $n_w(D_k) = 0$, such that either
    \[
    n_\Delta(D_k) - 2n_{\hs}(D_k) < 0,
    \]
    or $n_{\Delta}(D_k) - 2n_{\hs}(D_k) = 0$ and $S(\xs_k) \nsubseteq S(\ys_k)$. Since $\cI$ is finite and 
    \[
    |\{D \in D_{\cH}(\xs,\ys) : n_w(D) = 0\}| = 1, 
    \]
    after passing to a subsequence, we may assume that $\xs_k = \xs$, $\ys_k = \ys$, and $D_k = D$ for every $k \in \bN$. As $n_w(D) = 0$, we have $D = D^\Sigma_{\xs,\ys}$. Let $D_{\xs,\ys}$ be the corresponding matched domain in the quantum disc model $\bD^q_\beta$, given by Proposition~\ref{prop:domain-correspondence}.

    By SFT neck-stretching compactness~\cite[Theorems~10.3 and 10.6]{SFT}, in the setting of the cylindrical reformulation of Heegaard Floer homology as in the proof of \cite[Proposition~12.4]{Lipshitz}, after passing to a subsequence, the $J_k$-holomorphic curves $f_k$ converge to a holomorphic building whose planar part $f_{\bD}$, possibly consisting of several stories, has total domain $D_{\xs,\ys}$. Here, after passing to a subsequence, the asymptotic generators and the domain are fixed, and the corresponding curves have uniformly bounded energy.

    Let 
    \[
    \pi_{[0,1] \times \bR} \colon \Sigma \times [0,1] \times \bR \to [0,1] \times \bR 
    \]
    be the projection. The map $\pi_{[0,1] \times \bR}$ is holomorphic. The number $n_\Delta(D_{\xs,\ys}) = n_\Delta(f_{\bD})$ is the total branching order of $\pi_{[0,1] \times \bR} \circ f_{\bD}$, and is therefore non-negative. If $n_\Delta(f_{\bD}) = 0$, then every planar story is unbranched.
    
    By Definition~\ref{def:diagonal-intersection} and since $e(D) + 2n_{\hs}(D) = e(D_{\xs,\ys})$, we have
    \begin{equation}\label{eqn:diagonal-surgery}
    n_{\Delta}(D) - 2n_{\hs}(D) = n_{\Delta}(D_{\xs,\ys}) \ge 0.
    \end{equation}
    This rules out the possibility that $n_\Delta(D) - 2n_{\hs}(D) < 0$.

    Now suppose that $n_{\Delta}(D) - 2n_{\hs}(D) = 0$. Let $f^1,\dots,f^r$ be the planar stories of $f_{\bD}$, with intermediate generators $\xs = \zs^0, \zs^1, \dots, \zs^r = \ys$ (i.e., $f^k$ connects $\zs^{k-1}$ and $\zs^k$). Since $\qs \cap (T_1 \cup \dots \cup T_n) = \emptyset$, we have $n_{\qs}(D_{\xs,\ys}) = n_{\qs}(D)$. By equation~\eqref{eqn:diagonal-surgery}, we have 
    \begin{equation}\label{eqn:unbranched}
    n_\Delta(D_{\xs,\ys}) = n_\Delta(D) - 2n_{\hs}(D) = 0. 
    \end{equation}
    Fix $k \in \{1,\dots,r\}$. By equation~\eqref{eqn:unbranched}, the map $\pi_{[0,1] \times \bR} \circ f^k$ is unbranched: The domain of $f^k$ is a disjoint union $\coprod_{i=1}^n [0,1] \times \bR \times \{i\}$, and $\pi_{[0,1] \times \bR} \circ f^k$ is the trivial $n$-fold covering of $[0,1] \times \bR$. For $i \in \{1,\dots,n\}$, let
    \[
    g_i := f^k|_{[0,1] \times \bR \times \{i\}}.
    \]
    Let $\zs^{k-1} = (x_1,\dots,x_n)$ and $\zs^k = (y_1,\dots,y_n)$ be ordered such that $x_i$, $y_i \in \beta_i$. We can label the $g_i$ such that $\lim_{t \to -\infty} g_i(s,t,i) = x_i$ and $\lim_{t \to \infty} g_i(s,t,i) = y_i$.

    For $i \in \{1,\dots,n\}$, let $B_i$ be the bigon with corners $a_i \cap b_i$ and boundary contained in $a_i$ and the lower half of $b_i$; see Figure~\ref{fig:disc-model}. Since $\beta \cup \bI$ is the identity in $\{x > 0\}$ and $n_w(g_i) = 0$, the domain $\cD(g_i) \cap \{x>0\}$ is either $0$, $B_i \cap \{x>0\}$, or $(B_i \cup N(c_i)) \cap \{x>0\}$. Hence, $x_i = x_i^0$ and $y_i \neq x_i^0$ are impossible, as in this case, $\partial\cD(g_i) \cap a_i$ would be a 1-chain with boundary $y_i - x_i^0$. It follows that $S(\zs^{k-1}) \subseteq S(\zs^k)$ for $k \in \{1,\dots,r\}$, and so $S(\xs) \subseteq S(\ys)$. This rules out the possibility that $n_{\Delta}(D) - 2n_{\hs}(D) = 0$ and $S(\xs) \nsubseteq S(\ys)$.

    The inequalities $n_{\qs}(D) \ge 0$ and $n_{q_{n+1}}(D) \ge 0$ follow from the positivity of the intersection with the divisors $V_{q_i}$ for $i \in \{0,\dots,n+1\}$. Hence, the result holds if we set $\cJ(\cH^q_\beta)$ to be the set of generic almost complex structures in $\bigcup_{L \ge L_0} \cJ_L(\cH^q_\beta)$.
\end{proof}

\begin{coro}
Under the assumptions of Proposition~\ref{prop:positivity}, we have 
\begin{equation}\label{eqn:kappa}
(n_{\Delta}(D) - 2n_{\hs}(D) = 0) \wedge (n_{\qs}(D) = 0) \iff \kappa(D) = 0.
\end{equation}
\end{coro}

\begin{rmk}\label{rem:parametrised-positivity}
    For any compact one- or two-parameter family $\{J_\lambda : \lambda \in K\}$ of admissible almost complex structures satisfying the preceding conditions and $L \in \bR_+$, let $J_{\lambda, L}$ be obtained by inserting standard necks of length $L$. Then there is a single $L_0 \in \bR_+$ for which the conclusions of Proposition~\ref{prop:positivity} hold for every $\lambda \in K$ and $L \ge L_0$. Indeed, otherwise one could choose parameters $\lambda_k$, neck lengths $L_k \to \infty$, and counterexamples as in the proof of Proposition~\ref{prop:positivity}. After passing to a subsequence $\lambda_k \to \lambda_\infty$, parametrised neck-stretching compactness gives the same contradiction. We choose such families generically so that the corresponding parametrised moduli spaces are regular. For every compact one- or two-parameter family of admissible almost complex structures used in continuation equations with moving $\alpha$-boundary conditions, there is an $L_0 \in \bR_+$, independent of the parameter, such that the conclusions of Proposition~\ref{prop:positivity} hold for every $L \geq L_0$, provided the moving boundary conditions are constant along the tubes and the almost complex structures satisfy the preceding conditions near the marked cylinders.
\end{rmk}

The following local lemma about domains will play a central role.

\begin{lem}\label{lem:S-monotonicity}
    Let $\xs = (x_1,\dots,x_n)$, $\ys = (y_1,\dots,y_n) \in \cI$, labelled such that $x_i$, $y_i \in b_i$ for $i \in \{1,\dots,n\}$. Furthermore, let $D \in D_\cH(\xs,\ys)$ be a positive domain. If $n_{q_i}(D) = 0$ for some $i \in \{1,\dots,n\}$, then $x_i = x_i^0$ implies that $y_i = x_i^0$. In particular, if $n_{\qs}(D) = 0$, then $S(\xs) \subseteq S(\ys)$.
\end{lem}

\begin{proof}
    Suppose that $x_i = x_i^0$. Let $a$, $b$, $c$, and $d$ be the multiplicities of $D$ in the quadrants around $x_i^0$ labelled anticlockwise, starting from the quadrant containing $q_i$. Since $n_{q_i}(D)=0$, we have $a = 0$. As the two quadrants in $N(c_i) \setminus N(\{z_i,w_i\})$ lie in the same region (where $N(c_i)$ is the disc bounded by $b_i$ in $D^2$), we have $b = c$. Since $D \ge 0$, we have $d \ge 0$. By Definition~\ref{def:xy-domains}, we have 
    \[
    \partial(\partial D \cap \as) = y_i - x_i.
    \]
    The coefficient of $x_i^0$ in $\partial(\partial D \cap \as)$ is
    \[
    -a + b - c + d = d.
    \]
    As $d \ge 0$ and $x_i = x_i^0$, this is only possible if $d = 0$ and $y_i = x_i^0$.
\end{proof}

\begin{prop}\label{prop:W-S}
    Let $\beta \in B_{n+1}$ be a braid whose closure is a knot and $\cH^q_\beta$ an associated quantum Heegaard diagram. Suppose that $\xs$, $\ys \in \cI$, and $D \in D_{\cH}(\xs,\ys)$ is such that $n_w(D) = 0$. Then
    \begin{equation}\label{eqn:W-S}
    W(\ys) - W(\xs) = 3n_{\qs}(D) + |S(\ys)| - |S(\xs)|.
    \end{equation}
    Furthermore, if $D \ge 0$, then $|S(\xs)| \le |S(\ys)|$ and $W(\xs) \le W(\ys)$.
\end{prop}

\begin{proof}
    Given $H \subseteq \{1,\dots,n\}$, we write $\one_H$ for the indicator function of $H$. Set
    \[
    \omega_i(D) := \one_{S(\ys)}(i) - \one_{S(\xs)}(i) \in \{-1,0,1\}.
    \]
    Then
    \begin{equation}\label{eqn:omega-S}
    \sum_{i=1}^n \omega_i(D) = |S(\ys)| - |S(\xs)|.
    \end{equation}
    For $i \in \{1,\dots,n\}$, let $m_i(D)$ be the multiplicity of 
    \[
    a_i \cap \left\{x \in \bR^2 : 0 < x < \frac{1}{2(n+2)}\right\} \subseteq \alpha_i
    \]
    in $\partial_\alpha D := \partial D \cap \as$. Since $a_i \cap \bb \cap \{x>0\} = \{x_i^0\}$, the equation $\partial(\partial_\alpha D) = \ys - \xs$ gives $m_i(D) = \omega_i(D)$.

    As before, let $q_{n+1}$ be a point in the region neighbouring the one containing $q_n$ along $a_n$. For $i \in \{1,\dots,n+1\}$, join $w$ to $q_i$ by a path $\nu_i$ such that $\nu_i \pitchfork (\as \cup \bs)$ consists of $i-1$ points lying successively in $a_1,\dots,a_{i-1}$. Since $n_w(D) = 0$, the change in the coefficient of $D$ along $\nu_i$ gives
    \begin{equation}\label{eqn:omega-nq-i}
    n_{q_i}(D) = \sum_{j=1}^{i-1} \omega_j(D).
    \end{equation}
    As $q_0$ and $w$ lie in the same region, $n_{q_0}(D) = 0$. We obtain that
    \begin{equation}\label{eqn:omega-nq}
    n_{\qs}(D) = \sum_{i=1}^n n_{q_i}(D) = \sum_{i=1}^n \sum_{j=1}^{i-1} \omega_j(D) = \sum_{j=1}^n (n-j)\omega_j(D).
    \end{equation}
    By definition,
    \[
    W_i(\xs) = 
    \begin{cases}
        3(n-i)+1, & x_i = x_i^0,\\
        0, & x_i \neq x_i^0.
    \end{cases}
    \]
    Hence,
    \[
    W(\ys) - W(\xs) = \sum_{i=1}^n (3(n-i)+1) \omega_i(D) = 3\sum_{i=1}^n (n-i)\omega_i(D) + \sum_{i=1}^n \omega_i(D).
    \]
    Using equations~\eqref{eqn:omega-S} and~\eqref{eqn:omega-nq}, this becomes equation~\eqref{eqn:W-S}.

    Now assume that $D \ge 0$. By equation~\eqref{eqn:omega-nq-i}, we have 
    \[
    0 \le n_{q_{n+1}}(D) = \sum_{j=1}^n \omega_j(D) = |S(\ys)| - |S(\xs)|. 
    \]
    This, together with $n_{\qs}(D) \ge 0$ and equation~\eqref{eqn:W-S} implies that $W(\ys) - W(\xs) \ge 0$.
\end{proof}

\begin{defn}
    Consider $\bZ^2$ with the partial order $(a,b) \le (a',b')$ if and only if $a \le a'$ and $b \le b'$. We define the \emph{quantum bifiltration} of $\xs \in \cI$ to be
    \[
    \cF(\xs) := (P(\xs), -W(\xs)) \in \bZ^2.
    \]
\end{defn}

\begin{prop}\label{prop:filtration}
    Consider an admissible almost complex structure $J \in \cJ(\cH^q_\beta)$. Then the differential $\ph$ on the knot Floer complex $\CFKh(\Sigma, \as, \bs, w, z)$ defined in equation~\eqref{eqn:HFK-differential} is non-increasing with respect to $\cF$.
\end{prop}

\begin{proof}
    Let $\xs$, $\ys \in \cI$, and $u \in \pi_2(\xs,\ys)$ be such that $\cMh(u) \neq \emptyset$; i.e., $\ys$ appears in $\ph \xs$ with a non-zero coefficient. Then
    \[
    \cF(\xs) - \cF(\ys) = (\kappa(u), W(\ys) - W(\xs)). 
    \]
    We have $\kappa(u) = n_\Delta(u) - 2n_{\hs}(u) + n_{\qs}(u) \ge 0$ by Proposition~\ref{prop:positivity}. As $\cD(u) \ge 0$, we have $W(\xs) \le W(\ys)$ by Proposition~\ref{prop:W-S}. It follows that $\cF(\xs) \ge \cF(\ys)$.
\end{proof}

\begin{defn}\label{def:quantum-braid-Floer-differential}
    Let $\beta \in B_{n+1}$ be a braid, and let $J \in \cJ(\cH^q_\beta)$ be an almost complex structure. For $\xs \in \cI$, we let
    \[
    \partial^q \xs = \sum_{\substack{\ys \in \cI: \\ W(\xs) = W(\ys)}} \sum_{\substack{u \in \pi_2(\xs,\ys):\\ \mu(u)=1 \\ n_w(u) = n_z(u) = 0 \\ \kappa(u) = 0}} \#\cMh(u) \cdot \ys,
    \]
    and define the \emph{quantum braid Floer differential}
    \[
    \partial^q \colon \CFKh(\cH^q_\beta) \to \CFKh(\cH^q_\beta)
    \]
    by extending $\partial^q$ $\bF_2$-linearly.
\end{defn}

\begin{rmk}
    Suppose that $\xs$, $\ys \in \cI$ and $u \in \pi_2(\xs,\ys)$ is such that $\cMh(u) \neq \emptyset$. By Proposition~\ref{prop:positivity}, for $J \in \cJ(\cH_\beta^q)$, we have $\kappa(u) = 0$ if and only if
    \begin{equation}\label{eqn:b-q}
        n_\Delta(u) - 2n_{\hs}(u) = 0 \text{ and } n_{\qs}(u) = 0.
    \end{equation}
    Hence, we could replace the condition $\kappa(u) = 0$ with~\eqref{eqn:b-q} in Definition~\ref{def:quantum-braid-Floer-differential}.
\end{rmk}

\begin{prop}\label{prop:quantum-differential}
    The endomorphism $\partial^q$ satisfies the following:
    \begin{enumerate}
        \item\label{it:differential} $\partial^q \circ \partial^q = 0$;
        \item\label{it:triply-graded} $\partial^q$ decreases $M$ by one and preserves both $A^{\HF}$ and $A^{\qHF}$;
        \item\label{it:S} if $\xs$, $\ys \in \cI$, and $\ys$ appears in $\partial^q \xs$ with a non-zero coefficient, then $S(\xs) = S(\ys)$.
    \end{enumerate}
\end{prop}

\begin{proof}
    By Proposition~\ref{prop:filtration}, the map $\ph$ is non-increasing with respect to $\cF$. Note that $\partial^q$ is the degree-$(0,0)$ part of $\ph$ with respect to $\cF$. The property $\partial^q \circ \partial^q = 0$ now follows from $\ph \circ \ph = 0$.

    Since $\ph$ preserves $A^{\HF}$, so does $\partial^q$. If $\ys$ appears in $\partial^q \xs$ with a non-zero coefficient, then $W(\xs) = W(\ys)$, and there is some $u \in \pi_2(\xs,\ys)$ with $\kappa(u) = 0$. Hence, $P(\xs) = P(\ys)$ by equation~\eqref{eqn:P}. It follows from Definition~\ref{def:P} that $A^{\qHF}(\xs) = A^{\qHF}(\ys)$. As 
    \[
    n_\Delta(u) - 2n_{\hs}(u) + n_{\qs}(u) = \kappa(u) = 0, 
    \]
    it follows from Proposition~\ref{prop:positivity} that $n_{\qs}(u) = 0$. Hence, Lemma~\ref{lem:S-monotonicity} implies that $S(\xs) \subseteq S(\ys)$. By definition, 
    \[
    W(\xs) = \sum_{i \in S(\xs)} W_i(\xs), \quad
    W(\ys) = \sum_{i \in S(\ys)} W_i(\ys).
    \]
    We have $W_i(\ys) > 0$ for every $i \in S(\ys)$ and $W_i(\xs) = W_i(\ys)$ for $i \in S(\xs)$. As $W(\xs) = W(\ys)$, we have $S(\xs) = S(\ys)$. 
\end{proof}

\begin{defn}
    We call
    \[
    \left(\CFKh(\cH^q_\beta), \partial^q\right)
    \]
    the \emph{quantum braid Floer complex} of $\cH^q_\beta$. The \emph{quantum braid Floer homology} of $\cH^q_\beta$ is
    \[
    \qHFB(\cH^q_\beta) := H_*\left(\CFKh(\cH^q_\beta), \partial^q\right).
    \]
\end{defn}

\subsection{Quantum braid Floer homology is a braid invariant}

By the work of the first author~\cite{Anghel2024AIF}, the intersection polynomial $\Omega'_\beta(x,d)$ is not a knot invariant. It is neither invariant under Markov conjugation nor Markov stabilisation. By \cite[Theorem~1.6]{Anghel2024AIF}, we obtain a knot invariant if we pass to the quotient $\bZ[x^{\pm 1}, d^{\pm 1}]/((d+1)(xd-1))$. Hence, quantum braid Floer homology also cannot be a knot invariant, as it categorifies $\Omega^q_\beta(x,d)$ by Proposition~\ref{prop:categorification}, and since $\Omega^q_\beta(x,d) = \Omega'_\beta(x,-d)$ by Theorem~\ref{qINT}. In this section, we prove the following result, which implies that quantum braid Floer homology depends only on the braid $\beta$ and not on the choice of the action of the braid $\beta \cup \bI$ on the punctured disc $\cD_n^q$.

\begin{thm}\label{thm:braid-invariance}
    Let $\beta \in B_{n+1}$ be a braid whose closure $K$ is a knot. Let $d_0$ and $d_1$ be actions of the braid $\beta \cup \bI \in B_{3n+3}$ on $\cD_n^q$ relative to $\partial \cD_n^q$ supported in the $\{x<0\}$ half-plane. For $i \in \{0,1\}$, write $\ba_i' := d_i(\ba)$ and $\as_i$ for the closure of $\ba_i'$ in the Heegaard surface obtained by adding tubes to $\cD_n^q$. Let
    \[
    \cH^q_i = (\Sigma,\as_i,\bs, w, z, \qs, \hs, \xs^0)
    \]
    be the corresponding quantum Heegaard diagram and $\partial^q_i$ the respective quantum braid Floer differential. Choose $J_i \in \cJ(\cH^q_i)$ for $i \in \{0,1\}$. Then $\CFKh(\cH^q_0, \partial^q_0)$ and $\CFKh(\cH^q_1, \partial^q_1)$ have the same $(A^{\HF}, A^{\qHF}, M)$-graded chain homotopy type. In particular, 
    \[
    H_*\left(\CFKh(\cH^q_0, \partial^q_0)\right) \cong H_*\left(\CFKh(\cH^q_1, \partial^q_1)\right),
    \]
    and the isomorphism is $(A^{\HF}, A^{\qHF}, M)$-graded. 
\end{thm}

\begin{proof}
    Let us write $\cH_i := (\Sigma,\as_i,\bs, w, z)$ for $i \in \{0,1\}$. As we noted before, different choices of closing up the arcs $\ba_i'$ along the tubes to obtain the curves $\as_i$ lead to diffeomorphic quantum Heegaard diagrams, and hence isomorphic chain complexes.  By the work of the second author, Thurston, and Zemke~\cite[Section~9]{naturality} on the naturality of knot Floer homology and a chain-level extension of it due to Zemke~\cite[Proposition~3.5]{Zemke}, there is an $(A^{\HF}, M)$-graded chain homotopy equivalence
    \[
    \Phi_{\cH_0 \to \cH_1} \colon \CFKh(\cH_0) \to \CFKh(\cH_1)
    \]
    that is well-defined up to $(A^{\HF}, M)$-graded chain homotopy. Let $\cI_i := \bT_{\as_i} \cap \bT_{\bs}$, and let $\cF_i \colon \cI_i \to \bZ^2$ be the quantum bifiltration on $\cH^q_i$ for $i \in \{0,1\}$. We claim that $\Phi_{\cH_0 \to \cH_1}$ has a representative that preserves the quantum bifiltration; i.e.,
    \begin{equation}\label{eqn:filtered-naturality}
    \cF_1(\Phi_{\cH_0 \to \cH_1}(\xs)) \le \cF_0(\xs)
    \end{equation}
    for every $\xs \in \cI_0$.

    The diffeomorphisms $d_0$ and $d_1$ represent the same mapping class in $(\cD_n^q, \partial \cD_n^q)$ as they represent the same braid; hence, the identity is isotopic to $d_1d_0^{-1}$ relative to $\partial \cD_n^q$ through some isotopy $\{\iota_t : t \in [0,1]\}$. Furthermore, since $d_0$ and $d_1$ are supported in $\{x<0\}$, we can assume that $\iota_t$ is also supported in $\{x<0\}$ for every $t \in [0,1]$. Different choices of closing up the arcs $\ba_i'$ in the tubes $T_1,\dots,T_n$ lead to diffeomorphic quantum Heegaard diagrams (the diffeomorphism is a composition of Dehn twists and isotopies of the identity near the tubes), which tautologically have isomorphic quantum braid Floer complexes. Hence, we can assume that the differentials of $d_0$ and $d_1$ agree at $\zs$. So, we can choose $\iota_t$ for every $t \in [0,1]$ such that the differential $d_{z_i}\iota_t = \Id_{T_{z_i}D^2}$ for $i \in \{0,\dots,n\}$, and we can arrange that $\iota_t|_{\partial N(\zs)} = \Id$. Let 
    \[
    U := \{x<0\} \setminus N(\zs).
    \]
    Then the isotopy $\{\iota_t : t \in [0,1]\}$ extends to an isotopy $\{\iota_t^\Sigma : t \in [0,1]\}$ of $\Sigma$ that is supported in $U$. For $t \in [0,1]$, let us write $\as_t := \iota_t^\Sigma(\as_0)$.
    
    By the discussion preceding \cite[Definition~5.5]{OS}, we can break down the isotopy into Hamiltonian isotopies, each of which is supported in a small disc around a tangency between $\as_t$ and $\bs$, and diffeomorphisms of the Heegaard diagram supported in $\{x<0\}$. Hence, the chain map $\Phi_{\cH_0 \to \cH_1}$ is the composition of continuation maps induced by the Hamiltonian isotopies and diffeomorphism maps. Diffeomorphism maps preserve the quantum bifiltration since the diffeomorphisms preserve $\hs$, $\qs$, and $\xs^0$ (as they lie outside $U$) and the diagonal in $\Sym^n(\Sigma)$. Furthermore, each induces an invertible map on the degree-$(0,0)$ part.
    
    Hence, it suffices to prove the claim when $\{\iota_t : t \in [0,1]\}$ is a Hamiltonian isotopy supported in $U$. Let $\sI_D$ be the isotopy of $\Sym^n(D^2)$ induced by $\{\iota_t : t \in [0,1]\}$. Similarly, let $\sI$ be the isotopy of $\Sym^n(\Sigma)$ induced by $\{\iota_t^\Sigma : t \in [0,1]\}$. In this case, $\Phi_{\cH^q_0 \to \cH^q_1}$ is the continuation map $\Gamma^{\sI}$ defined by Ozsv\'ath and Szab\'o~\cite[p.~54]{OS}. Choose a smooth function $\rho \colon \bR \to [0,1]$ such that $\rho|_{[0,1]} \colon [0,1] \to [0,1]$ is a homeomorphism and 
    \[
    \rho(t) = 
    \begin{cases}
        0, & t \le 0; \\
        1, & t \ge 1.
    \end{cases}
    \] 
    For $\xs \in \bT_{\as_0} \cap \bT_{\bs}$ and $\ys \in \bT_{\as_1} \cap \bT_{\bs}$, let $\pi_2^\sI(\xs,\ys)$ be the set of homology classes of maps 
    \[
    u \colon [0,1] \times \bR \to \Sym^n(\Sigma) 
    \]
    such that $u(1,t) \in \bT_{\as_{\rho(t)}}$ and $u(0,t) \in \bT_{\bs}$ for $t \in \bR$, with asymptotics 
    \[
    \lim_{t \to -\infty} u(s,t) = \xs, \quad \lim_{t \to \infty} u(s,t) = \ys.
    \]

    Choose a generic one-parameter family of admissible almost complex structures connecting $J_0$ to $J_1$. After first increasing the neck lengths at the endpoints, if necessary, the uniform statement of Remark~\ref{rem:parametrised-positivity} allows us to choose a common sufficiently large neck length along this path. Concatenating the two continuation paths with this uniformly stretched path gives a generic family $\cK = \{J_t : t \in [0,1]\}$ of admissible almost complex structures that is constant near the endpoints and is uniformly sufficiently stretched.
    
    In the moving boundary continuation equation, we use $\bT_{\as_{\rho(t)}}$ and $J_{\rho(t)}$. Given $u \in \pi_2^{\sI}(\xs,\ys)$, we write $\cM^{\sI, \cK}(u)$ for the moduli space of pseudo-holomorphic representatives of $u$. We then let
    \[
    \Gamma^{\sI, \cK}(\xs) := \sum_{\ys \in \bT_{\as_1} \cap \bT_{\bs}} \sum_{\substack{u \in \pi_2^{\sI}(\xs,\ys): \\ \mu(u) = 0\\ n_w(u) = n_z(u) = 0}} \#\cM^{\sI, \cK}(u) \cdot \ys.
    \]

    As $V_p \cap \bT_{\as_t} = V_p \cap \bT_{\bs} = \emptyset$ for $p \in \hs \cup \qs$ and $\Delta \cap \bT_{\as_t} = \Delta \cap \bT_{\bs} = \emptyset$ for $t \in [0,1]$, we can define $n_{\hs}(u)$, $n_{\qs}(u)$, and $n_{\Delta}(u)$ for $u \in \pi_2^{\sI}(\xs,\ys)$. Set
    \[
    \kappa(u) := n_\Delta(u) - 2n_{\hs}(u) + n_{\qs}(u).
    \]

    Suppose that $u \in \pi_2^{\sI}(\xs,\ys)$ contributes to $\Gamma^{\sI, \cK}$. In particular, $n_w(u) = n_z(u) = 0$. Let $u_D \in \pi_2^{\sI_D}(\xs,\ys)$ denote the continuation class obtained by surgering $\Sigma$ along the cores of the tubes $T_1,\dots,T_n$. The diffeomorphism $\iota_t^\Sigma$ is the identity in a neighbourhood of $\xs^0$ for every $t \in [0,1]$. Let $u^0(s,t) \equiv \xs^0$ in $\pi_2^\sI(\xs^0,\xs^0)$ be the class of the constant continuation strip. Then 
    \begin{equation}\label{eqn:constant-strip}
    n_\Delta(u^0) = 0, \quad n_{\hs}(u^0) = 0, \quad n_{\qs}(u^0) = 0, \quad n_\Delta(u^0_D) = 0, 
    \end{equation}
    and hence $\kappa(u^0) = 0$.
        
    Choose $u_{\xs} \in \pi_2(\xs,\xs^0)$ in $\cH_0$ and $u_{\ys} \in \pi_2(\ys,\xs^0)$ in $\cH_1$ such that $n_w(u_{\xs}) = n_w(u_{\ys}) = 0$. Then $u_{\xs} \ast u^0$, $u \ast u_{\ys} \in \pi_2^{\sI}(\xs,\xs_0)$. Since $\pi_2^{\sI}(\xs,\xs_0) \cong \bZ$ as an affine space, where the isomorphism is given by $n_w$, and because $n_w(u_{\xs} \ast u^0) = n_w(u \ast u_{\ys}) = 0$, we have $u_{\xs} \ast u^0 = u \ast u_{\ys}$. As $n_\Delta$, $n_{\hs}$, and $n_{\qs}$ are additive with respect to concatenation,
    \[
    \kappa(u_{\xs}) + \kappa(u^0) = \kappa(u) + \kappa(u_{\ys}).
    \]
    Because $n_w(u_{\xs}) = n_w(u_{\ys}) = 0$ and $\kappa(u^0) = 0$, we have 
    \begin{equation}\label{eqn:kappa-continuation}
    \begin{split}
    P(\xs) - P(\ys) &= (P(\xs)-P(\xs_0)) - (P(\ys)-P(\xs_0))\\
    &= \kappa(u_{\xs}) - \kappa(u_{\ys})\\
    &= \kappa(u).
    \end{split}
    \end{equation}
    It follows that
    \begin{equation}\label{eqn:qHF-continuation}
        A^{\qHF}(\xs) - A^{\qHF}(\ys) = \kappa(u) + W(\xs) - W(\ys).
    \end{equation}
    We now prove a moving boundary generalisation of equation~\eqref{eqn:diagonal-surgery}.

    \begin{lem}\label{lem:continuation-surgery}
        For $\xs \in \bT_{\as_0} \cap \bT_{\bs}$, $\ys \in \bT_{\as_1} \cap \bT_{\bs}$, and $u \in \pi_2^\sI(\xs,\ys)$ with $n_w(u) = 0$, we have
        \[
        n_\Delta(u_D) = n_\Delta(u) - 2n_{\hs}(u).
        \]
    \end{lem}

    \begin{proof}
        Since $(\Sigma, \as_0, \bs, z)$ is a diagram of $S^3$, we have 
        \[
        \pi_2(\xs^0, \xs) \cong \bZ, \quad \pi_2(\xs^0, \ys) \cong \bZ, \quad \pi_2^{\sI}(\xs^0, \ys) \cong \bZ 
        \]
        as affine spaces, where the isomorphisms are given by $n_w$. Let $A \in \pi_2(\xs^0, \xs)$ and $B \in \pi_2(\xs^0, \ys)$ be the unique classes such that $n_w(A) = n_w(B) = 0$. Then 
        \[
        A \ast u = u^0 \ast B \in \pi_2^{\sI}(\xs^0, \ys)
        \]
        as $n_w(A \ast u) = 0$ and $n_w(u^0 \ast B) = 0$. Since $n_\Delta$ and $n_{\hs}$ are additive under concatenation,
        \begin{equation}\label{eqn:concatenation}
        \begin{split}
        n_\Delta(A) + n_\Delta(u) &= n_\Delta(u^0) + n_\Delta(B), \\
        n_{\hs}(A) + n_{\hs}(u) &= n_{\hs}(u^0) + n_{\hs}(B).
        \end{split}
        \end{equation}
        As the isotopy $\{\iota_t^\Sigma : t \in [0,1]\}$ is constant along the tubes,
        \[
        A_D \ast u_D = u^0_D \ast B_D,
        \]
        and so
        \begin{equation}\label{eqn:concatenation-disc}
        n_\Delta(A_D) + n_\Delta(u_D) = n_\Delta(u^0_D) + n_\Delta(B_D).
        \end{equation}
        Hence,
        \[
        \begin{split}
        n_\Delta(u_D) &= n_\Delta(u_D) - n_\Delta(u_D^0) \\
        &= n_\Delta(B_D) - n_\Delta(A_D) \\
        &= (n_\Delta(B) - 2n_{\hs}(B)) - (n_\Delta(A) - 2n_{\hs}(A)) \\
        &= n_\Delta(u) - n_\Delta(u^0) - 2(n_{\hs}(u) - n_{\hs}(u^0)) \\
        &= n_\Delta(u) - 2n_{\hs}(u),
        \end{split}
        \]
        where the first equality follows from equation~\eqref{eqn:constant-strip}, the second from equation~\eqref{eqn:concatenation-disc}, the third from equation~\eqref{eqn:diagonal-surgery}, the fourth from equation~\eqref{eqn:concatenation}, and the last from equation~\eqref{eqn:constant-strip}.
    \end{proof}
    
    Extending Proposition~\ref{prop:positivity}, we claim that, for $u \in \pi_2^\sI(\xs,\ys)$ with $\cM^{\sI, \cK}(u) \neq \emptyset$, we have
    \begin{equation}\label{eqn:continuation-positivity}
    n_\Delta(u) - 2n_{\hs}(u) \ge 0.
    \end{equation}
    Indeed, the isotopy $\{\iota^\Sigma_t : t \in [0,1]\}$ is constant along the tubes. Hence, by Remark~\ref{rem:parametrised-positivity}, every continuation strip contributing to $\Gamma^{\sI, \cK}$ satisfies $n_\Delta(u_D) \ge 0$. Consequently,
    \[
    n_\Delta(u) - 2n_{\hs}(u) = n_\Delta(u_D) \ge 0
    \]
    by Lemma~\ref{lem:continuation-surgery}. As the Hamiltonian isotopy is constant around $\qs$, positivity of intersection also implies that $n_{\qs}(u) \ge 0$ if $\cM^{\sI, K}(u) \neq \emptyset$. It follows from equation~\eqref{eqn:continuation-positivity} that $\kappa(u) \ge 0$.

    By equation~\eqref{eqn:kappa-continuation},
    \[
    \cF_0(\xs) - \cF_1(\ys) = (\kappa(u), W(\ys) - W(\xs)).
    \]
    Since the isotopy is constant around $\xs^0$, an analogue of Proposition~\ref{prop:W-S} implies that $W(\xs) \le W(\ys)$. Hence,
    \[
    \cF_1\left(\Gamma^{\sI, \cK}(\xs)\right) \le \cF_0(\xs).
    \]
    Let $\Gamma^{\sI, \cK}_q$ be the degree-$(0,0)$ part of $\Gamma^{\sI, \cK}$. For $\xs \in \cI_0$, this is given by
    \[
    \Gamma^{\sI, \cK}_q(\xs) = \sum_{\substack{\ys \in \cI_1: \\ W(\xs) = W(\ys)}} \sum_{\substack{u \in \pi_2^{\sI}(\xs,\ys):\\ \mu(u)=0 \\ n_w(u) = n_z(u) = 0 \\ \kappa(u) = 0}} \#\cM^{\sI, \cK}(u) \cdot \ys.
    \]
    Hence, if $\ys \in \cI_1$ appears in $\Gamma^{\sI, \cK}_q$ with a non-zero coefficient, then
    \[
    A^{\qHF}(\xs) = A^{\qHF}(\ys).
    \]
    Thus $\Gamma^{\sI, \cK}_q$ preserves the quantum grading. If we take the degree-$(0,0)$ part of the equation $\ph_1 \Gamma^{\sI, \cK} = \Gamma^{\sI, \cK} \ph_0$, we obtain that 
    \[
    \partial^q_1 \Gamma_q^{\sI, \cK} = \Gamma_q^{\sI, \cK} \partial^q_0;
    \]
    i.e., $\Gamma^{\sI, \cK}_q$ is a chain map.

    Let $\bar\sI$ denote the reverse of the isotopy $\sI$, and write $\bar \cK := \{J_{1-t} : t \in [0,1]\}$. This gives the filtered map $\Gamma^{\bar\sI, \bar \cK}_q$. By Ozsv\'ath and Szab\'o~\cite[pp.~55--56]{OS}, we have chain homotopies $\Gamma^{\bar\sI, \bar \cK}\Gamma^{\sI, \cK} \cong \Id$ and $\Gamma^{\sI, \cK}\Gamma^{\bar\sI, \bar \cK} \cong \Id$. We claim that these chain homotopies are filtered. Let $\{\Psi_\tau : \tau \in [0,1]\}$ be a homotopy through isotopies between the isotopies $\Psi_0 = \sI \bar \sI$ and the constant identity isotopy $\Psi_1$. Furthermore, let $\{\cK_\tau : \tau \in [0,1]\}$ be a generic 2-parameter family of uniformly sufficiently stretched almost complex structures for the continuation equations associated with $\{\Psi_\tau : \tau \in [0,1]\}$, as in Remark~\ref{rem:parametrised-positivity}, such that $\cK_0 = \cK\bar \cK$ and $\cK_1 \equiv J_0$. We can arrange that $\Psi_\tau$ is supported in $U$ and is constant around $\qs$, and the tubes for every $\tau \in [0,1]$. Let $H$ be the resulting chain homotopy. If $u \in \pi_2(\xs,\ys)$ for $\xs$, $\ys \in \bT_{\as_0} \cap \bT_{\bs}$, $\mu(u) = -1$, $n_w(u) = n_z(u) = 0$, and $f \in \cM^{\Psi_\tau, \cK_\tau}(u)$ for some $\tau \in [0,1]$, then 
    \[
    (\tau, f) \in \bigcup_{\tau \in [0,1]} \cM^{\Psi_\tau, \cK_\tau}(u)
    \]
    is counted by $H(\xs)$ with output $\ys$ (this union is zero-dimensional when $\mu(u) = -1$). The uniform two-parameter statement of Remark~\ref{rem:parametrised-positivity} applies to the family $\{\Psi_\tau : \tau \in [0,1]\}$, so every parametrised solution counted by $H$ satisfies $n_\Delta(u) - 2n_{\hs}(u) \ge 0$ and $n_{\qs}(u) \ge 0$; hence, $\kappa(u) \ge 0$. Furthermore, the argument showing equation~\eqref{eqn:kappa-continuation} at $\tau$ gives
    \[
    P(\xs) - P(\ys) = \kappa(u) \ge 0.
    \]
    Since the isotopy $\Psi_\tau$ is fixed around the points and paths used in the proof of Proposition~\ref{prop:W-S}, positivity of intersection of $u$ with $V_{q_i}$ for $i \in \{0,\dots,n+1\}$ gives the required inequalities, and the equations used in the proof remain unchanged. Hence, $W(\xs) \le W(\ys)$. We conclude that
    \[
    \cF(\xs) - \cF(\ys) = (\kappa(u), W(\ys) - W(\xs)) \in \bN^2,
    \]
    so $H$ is filtered. The same argument applies to the homotopy for $\bar \sI \sI$ and $\bar \cK \cK$. Taking the degree-$(0,0)$ components of
    \[
    \Gamma^{\bar\sI, \bar \cK} \Gamma^{\sI, \cK} - \Id = \partial H + H \partial
    \]
    therefore gives
    \[
    \Gamma^{\bar\sI, \bar \cK}_q \Gamma^{\sI, \cK}_q - \Id = \partial_q H_q + H_q \partial_q,
    \]
    and similarly for $\Gamma^{\sI, \cK}_q \Gamma^{\bar\sI, \bar \cK}_q$. Hence,
    \[
    \Gamma^{\bar\sI, \bar \cK}_q\Gamma^{\sI, \cK}_q \cong \Id, \quad \Gamma^{\sI, \cK}_q\Gamma^{\bar\sI, \bar \cK}_q \cong \Id,
    \]
    and the result follows.
    
    We remark that one can use small triangle maps instead of continuation maps. In fact, if the triple $(\Sigma, \as_1, \as_0, \bs, w, z)$ is admissible, then \cite[Lemma~9.7]{naturality} implies that $\Gamma^{\sI, \cK}$ is chain homotopic to the triangle map $F_{\as_1,\as_0,\bs}(\Theta_{\as_1,\as_0} \otimes -)$ for a suitable choice of $\cK$.
\end{proof}

\begin{defn}
    Let $\beta \in B_{n+1}$ be a braid whose closure is a knot. Then we write
    \[
    \qHFB(\beta) := \qHFB(\cH^q_\beta).
    \]
\end{defn}

\subsection{Properties of quantum braid Floer homology}

By part~\eqref{it:triply-graded} of Proposition~\ref{prop:quantum-differential}, quantum braid Floer homology is triply-graded by $(A^{\HF}, A^{\qHF}, M)$:
\[
\qHFB(\beta) = \bigoplus_{(a,q,m) \in \bZ^3} \qHFB_{a,q,m}(\beta).
\]
It categorifies the quantum intersection polynomial:

\begin{prop}\label{prop:categorification}
    Let $\beta$ be a braid whose closure is a knot. Then
    \[
    \sum_{(a,q,m) \in \bZ^3} (-1)^m \rank\left(\qHFB_{a,q,m}(\beta)\right)x^a d^q = \Omega^q_\beta(x,d).
    \]
\end{prop}

\begin{proof}
    The quantum braid Floer complex is generated by $\cI$. By part~\eqref{it:triply-graded} of Proposition~\ref{prop:quantum-differential}, the differential $\partial^q$ preserves all three gradings. Hence,
    \[
    \begin{split}
    \sum_{(a,q,m) \in \bZ^3} (-1)^m \rank\left(\qHFB_{a,q,m}(\beta)\right)x^a d^q &=
    \sum_{\xs\in \cI} (-1)^{M(\xs)} \, x^{A^{\HF}(\xs)} d^{A^{\qHF}(\xs)} \\
    &= (-1)^{\frac{n(n+1)}{2}}\sum_{\xs\in \cI} \varepsilon_{\xs} \, x^{A^{\HF}(\xs)} d^{A^{\qHF}(\xs)} \\
    &= \Omega^q_\beta(x,d).
    \end{split}
    \]
    where the second equality follows from Proposition~\ref{prop:signs} and the third is Definition~\ref{def:quantum-intersection}.
\end{proof}

By Theorem~\ref{qINT}, we have $V_{\mK}(x) = \Omega^q_\beta(x,-x^{-1})$. Hence, we obtain the following:

\begin{coro}\label{coro:Jones}
    If $\beta$ is a braid whose closure is a knot $K$, then
    \[
    V_{\mK}(x) = \sum_{(a,q,m) \in \bZ^3} (-1)^{m+q} \rank\left(\qHFB_{a,q,m}(\beta)\right)x^{a-q}.
    \]
\end{coro}

\begin{defn}
    For $\xs$, $\ys \in \cI$, and $D \in D_{\cH}(\xs,\ys)$, we write
    \[
    \begin{split}
    \fb(D) &:= n_\Delta(D) - 2n_{\hs}(D),\\ 
    \fc(D) &:= W(\ys) - W(\xs) - n_{\qs}(D),\\
    \fq(D) &:= n_{\qs}(D),\\
    \Delta W(D) &:= W(\ys) - W(\xs).
    \end{split}
    \]
    If $u \in \pi_2(\xs,\ys)$ and $D := \cD(u)$, then we put $\fb(u) := \fb(D)$, $\fc(u) := \fc(D)$, $\fq(u) := \fq(D)$, and $\Delta W(u) := \Delta W(D)$.
\end{defn}

Using this notation, if $\xs$, $\ys \in \cI$, and $D \in D_{\cH}(\xs,\ys)$ is such that $n_w(D) = 0$, then
\[
P(\xs) - P(\ys) = \kappa(D) = \fb(D) + \fq(D)
\]
by equation~\eqref{eqn:P}.

\begin{thm}\label{thm:spectral-sequence}
    Let $\beta \in B_{n+1}$ be a braid whose closure is a knot $K$. Then there is a spectral sequence whose $E^1$ page is $\qHFB(\beta)$ and which converges to $\HFKh(K)$ equipped with the filtration induced by $G := P - W$. More precisely, 
    \[
    E^\infty \cong \gr_G\left(\HFKh(K)\right).
    \]
\end{thm}

\begin{proof}
    Let $\cH^q_\beta$ be a quantum Heegaard diagram associated with $\beta$. Define
    \[
    G(\xs) := P(\xs) - W(\xs) = A^{\qHF}(\xs) - 2W(\xs). 
    \]
    If there is $u \in \pi_2(\xs,\ys)$ such that $\cMh(u) \neq \emptyset$, then $\cD(u) \ge 0$. Furthermore, $\fb(u) \ge 0$ and $\fq(u) \ge 0$ by Proposition~\ref{prop:positivity}, and $|S(\xs)| \le |S(\ys)|$ by Proposition~\ref{prop:W-S}. Hence,
    \[
    \begin{split}
        G(\xs) - G(\ys) &= \fb(u) + \fq(u) + W(\ys) - W(\xs) \\
        &= \fb(u) + 4\fq(u) + |S(\ys)| - |S(\xs)| \\
        &\ge 0, 
    \end{split}
    \]
    where the first equality follows from equation~\eqref{eqn:P} and the second from Proposition~\ref{prop:W-S}. By Lemma~\ref{lem:S-monotonicity}, equality holds exactly when $\fb(u) = 0$, $\fq(u) = 0$, and $S(\xs) = S(\ys)$, so the degree-zero part of $\ph$ with respect to $G$ is $\partial^q$. Hence, the spectral sequence associated with the knot Floer complex with the filtration $G$ has $E^1$ page $\qHFB(\beta)$ and $E^\infty$ page $\gr_G \HFKh(K)$.
\end{proof}

\begin{defn}
    For $\xs \in \cI$, write $s(\xs) := |S(\xs)|$ and
    \[
    R(\xs) := \sum_{i \in S(\xs)} (n-i) = \frac{W(\xs) - s(\xs)}{3}.
    \]
\end{defn}

By Proposition~\ref{prop:W-S}, if $\xs$, $\ys \in \cI$, and $D \in D_{\cH}(\xs,\ys)$ satisfies $n_w(D) = 0$, then
\[
W(\ys) - W(\xs) = 3\fq(D) + s(\ys) - s(\xs),
\]
which is equivalent to
\[
\fq(D) = R(\ys) - R(\xs).
\]

\begin{defn}
    For $\xs \in \cI$, let
    \[
    \begin{split}
        \cB(\xs) &:= P(\xs) + R(\xs), \\
        \cC(\xs) &:= R(\xs) - W(\xs).
    \end{split}
    \]
\end{defn}

We have
\[
\cB(\xs) - \cB(\ys) = P(\xs) - P(\ys) + R(\xs) - R(\ys) = \fb(D) + \fq(D) - \fq(D) = \fb(D).
\]
Furthermore,
\[
\cC(\xs) - \cC(\ys) = R(\xs) - R(\ys) + W(\ys) - W(\xs) = -\fq(D) + W(\ys) - W(\xs) = \fc(D).
\]

\begin{prop}\label{prop:bifiltration-BC}
    The pair $(\cB, \cC)$ is a bifiltration on $(\CFKh(\cH^q_\beta), \ph)$. The degree-$(0,0)$ part of $\ph$ is $\partial^q$.
\end{prop}

\begin{proof}
By Proposition~\ref{prop:W-S}, we have
\begin{equation}\label{eqn:fc}
    \fc(D) = 2\fq(D) + s(\ys) - s(\xs).
\end{equation}
Furthermore, if $D \ge 0$ and $n_w(D) = 0$, then $s(\xs) \le s(\ys)$. Hence, $\fc(D) \ge 0$. By Proposition~\ref{prop:positivity}, for a suitably stretched almost complex structure, $\fb(D) \ge 0$ if $D$ has a pseudo-holomorphic representative. So, if $u$ contributes to $\ph$, then $(\fb(u), \fc(u)) \in \bN^2$. As both $\fb$ and $\fc$ are additive under concatenation, $\ph$ is bifiltered with respect to $(\cB,\cC)$. 

Suppose that $\fb(u) = \fc(u) = 0$. Then $\fq(u) = 0$ and $|S(\xs)| = |S(\ys)|$. It now follows from Proposition~\ref{prop:W-S} and Lemma~\ref{lem:S-monotonicity} that $S(\xs) = S(\ys)$ and $W(\xs) = W(\ys)$. Hence, $\fb(u) = \fc(u) = 0$ is equivalent to $\kappa(u) = \fb(u) + \fq(u) = 0$ and $W(\xs) = W(\ys)$, which are the requirements for $u$ to contribute to $\partial^q$. Hence, $\partial^q$ is the degree-$(0,0)$ part of $\ph$ with respect to the bifiltration $(\cB,\cC)$.
\end{proof}

By Corollary~\ref{coro:Jones}, the Jones grading of $\xs \in \cI$ is
\[
J(\xs) := A^{\HF}(\xs) - A^{\qHF}(\xs).
\]
If $u \in \pi_2(\xs,\ys)$ and $n_w(u) = n_z(u) = 0$, then $A^{\HF}(\xs) = A^{\HF}(\ys)$. Furthermore, by Proposition~\ref{prop:W-S},
\[
A^{\qHF}(\xs) - A^{\qHF}(\ys) = \fb(u) + \fq(u) - (W(\ys) - W(\xs)) = \fb(u) - 2\fq(u) - (s(\ys) - s(\xs)).
\]
Hence, by equation~\eqref{eqn:fc},
\begin{equation}\label{eqn:Jones-bc}
J(\xs) - J(\ys) = 2\fq(u) + s(\ys) - s(\xs) - \fb(u) = \fc(u) - \fb(u).
\end{equation}
This motivates the following definition.

\begin{defn}
    Consider the two-variable polynomial ring $\bF_2[U,V]$, and let
    \[
    \qCFB^+(\cH^q_\beta) := \CFKh(\cH^q_\beta) \otimes_{\bF_2} \bF_2[U,V].
    \]
    For $\xs \in \cI$, we define
    \[
    \partial^q_+ \xs := \sum_{\ys \in \cI} \sum_{\substack{u \in \pi_2(\xs,\ys):\\ \mu(u)=1,\\ n_w(u) = n_z(u) = 0}} \#\cMh(u) U^{\fb(u)} V^{\fc(u)} \ys,
    \]
    and extend it to $\qCFB^+(\cH^q_\beta)$ such that it is $\bF_2[U,V]$-linear.
\end{defn}

That $\partial^q_+ \circ \partial^q_+ = 0$ follows from the additivity of $\fb$ and $\fc$ under gluing and analysing the ends of 1-dimensional moduli spaces. If we set $\deg_J(U) = -1$ and $\deg_J(V) = 1$, then the complex becomes Jones-homogeneous by equation~\eqref{eqn:Jones-bc}. When $U = V = 1$, we obtain $\ph$ and knot Floer homology. For $U=V=0$, we get $\partial^q$ and quantum braid Floer homology. 

With the notation of Theorem~\ref{thm:braid-invariance}, for $i \in \{0,1\}$, let $S_i$, $W_i$, $P_i$, $R_i$, $\cB_i$, and $\cC_i$ denote the corresponding functions on $\cI_i$. Thus,
\[
R_i(\xs) := \sum_{j \in S_i(\xs)} (n-j), \quad \cB_i(\xs) = P_i(\xs) + R_i(\xs), \quad \cC_i(\xs) = R_i(\xs) - W_i(\xs).
\]

\begin{lem}\label{lem:W-S-continuation}
    With the assumptions of Theorem~\ref{thm:braid-invariance}, let $\xs \in \bT_{\as_0} \cap \bT_{\bs}$, $\ys \in \bT_{\as_1} \cap \bT_{\bs}$, and let $u \in \pi_2^{\sI}(\xs, \ys)$ be a continuation class satisfying $n_w(u) = 0$. Then
    \[
    \begin{split}
    n_{\qs}(u) &= R_1(\ys) - R_0(\xs); \\
    n_{q_{n+1}}(u) &= |S_1(\ys)| - |S_0(\xs)|; \\
    W_1(\ys) - W_0(\xs) &= 3n_{\qs}(u) + |S_1(\ys)| - |S_0(\xs)|. \\
    \end{split}
    \]
\end{lem}

\begin{proof}
    The proof is analogous to that of Proposition~\ref{prop:W-S} because the isotopy is fixed near $q_0,\dots,q_{n+1}$, the canonical intersection point $\xs^0$, and the paths from $w$ to the $q_i$. For $j \in \{1,\dots,n\}$, let
    \[
    \omega_j(u) := \one_{S_1(\ys)}(j) - \one_{S_0(\xs)}(j).
    \]
    Since the isotopy is fixed along the paths from $w$ to the $q_i$, the same argument computing the change of multiplicity as in Proposition~\ref{prop:W-S} gives
    \[
    n_{q_i}(u) - n_w(u) = \sum_{j=1}^{i-1} \omega_j(u)
    \]
    for $i \in \{1,\dots,n+1\}$. As $n_w(u) = 0$, the stated identities follow.
\end{proof}

\begin{prop}
    The $\bF_2[U,V]$-equivariant chain homotopy type of
    \[
    (\qCFB^+(\cH^q_\beta), \partial^q_+)
    \]
    is independent of the representative of the braid action and of the choice of $J \in \cJ(\cH^q_\beta)$.
\end{prop}

\begin{proof}
    By the proof of Theorem~\ref{thm:braid-invariance}, the continuation maps, their reverses, and the corresponding chain homotopies can be chosen to be filtered with respect to the bifiltration $(\cB, \cC)$. Indeed, if a continuation strip $u$ connects $\xs$ to $\ys$, and if we write $\fb(u) = n_\Delta(u) - 2n_{\hs}(u)$ and $\fc(u) = W_1(\ys) - W_0(\xs) - n_{\qs}(u)$, then
    \[
    \begin{split}
    \cB_0(\xs) - \cB_1(\ys) &= \fb(u); \\
    \cC_0(\xs) - \cC_1(\ys) &= \fc(u) = 2n_{\qs}(u) + |S_1(\ys)| - |S_0(\xs)|
    \end{split}
    \]
    by Lemma~\ref{lem:W-S-continuation}. Lemma~\ref{lem:continuation-surgery} and the moving boundary neck-stretching argument give $\fb(u) \ge 0$, while positivity at $V_{q_i}$ for $i \in \{0,\dots,n\}$ gives $n_{\qs}(u) \ge 0$, and positivity at $V_{q_{n+1}}$ gives $|S_1(\ys)| - |S_0(\xs)| \ge 0$. Thus, $\fc(u) \ge 0$. The same calculation applies to the two-parameter continuation homotopies. The diffeomorphism maps tautologically preserve $\cB$ and $\cC$.
    
    Given a $(\cB, \cC)$-filtered map $\Phi$, we define its homogenisation via
    \[
    \Phi^+(\xs) := \sum_{\ys \in \cI_1} \phi_{\xs,\ys}\, U^{\cB_0(\xs)-\cB_1(\ys)} V^{\cC_0(\xs)-\cC_1(\ys)} \ys,
    \]
    where $\xs \in \cI_0$ and $\Phi(\xs) = \sum_{\ys \in \cI_1} \phi_{\xs,\ys}\, \ys$. Homogenisation respects compositions, differentials, and filtered chain homotopies. Hence, the continuation equivalences and their filtered homotopies induce mutually inverse $\bF_2[U,V]$-equivariant chain homotopy equivalences of the $\qCFB^+$ complexes.
\end{proof}

\begin{prop}\label{prop:theta}
    Let $\beta \in B_{n+1}$ be a braid whose closure is a knot. Then $\partial^q \xs^0 = 0$. Its homology class
    \[
    \Theta(\beta) := [\xs^0] \in \qHFB_{\frac{w(\beta)+n}{2}, w(\beta)+n, m}(\beta)
    \]
    is non-zero, where $m = M(\Theta(\beta))$.
\end{prop}

\begin{proof}
    For $\xs \in \cI$, we have $S(\xs) = \{1,\dots,n\}$ if and only if $\xs = \xs^0$. By Proposition~\ref{prop:quantum-differential}, if $\xs$, $\ys \in \cI$, and $\ys$ appears in $\partial^q\xs$ with a non-zero coefficient, then $S(\xs) = S(\ys)$. Hence, if $\ys$ appears in $\partial^q \xs^0$, then $\ys = \xs^0$, but that is impossible as $M(\ys) = M(\xs^0)-1$. So, $\partial^q \xs^0 = 0$. Similarly, if $\xs^0$ appears in $\partial^q \xs$, then $\xs = \xs^0$, which leads to a contradiction. According to Corollary~\ref{coro:canonical-HF-grading} and Definition~\ref{def:canonical-quantum-grading}, we have 
    \[
    A^{\HF}(\xs^0) = \frac{w(\beta) + n}{2}, \quad A^{\qHF}(\xs^0) = w(\beta)+n.
    \]
    
    The class $[\xs^0]$ is preserved by the continuation map $\Gamma_q^{\sI, \cK}$ constructed in Theorem~\ref{thm:braid-invariance}. Indeed, $\Gamma_q^{\sI, \cK}$ preserves $W$, and $\xs^0$ is the unique generator of maximal $W$. Hence, there is $c \in \bF_2$ such that
    \[
    \Gamma_q^{\sI, \cK}(\xs^0) = c\,\xs^0.
    \]
    Since $[\xs^0] \neq 0$ and $\Gamma_q^{\sI, \cK}$ induces an isomorphism on homology, $c = 1$. Hence,
    \[
    (\Gamma_q^{\sI, \cK})_*[\xs^0] = [\xs^0],
    \] 
    and $[\xs^0]$ is independent of all the choices appearing in Theorem~\ref{thm:braid-invariance}.
\end{proof}

\section{Examples}\label{sec:examples}

In this section, we compute the quantum intersection polynomial $\Omega^q_{\sigma^3}$ associated with the braid $\sigma^3 \in B_2$. We will compute $\Omega^q_{\sigma^3}(x,1)$ and verify that it agrees with the Alexander polynomial of the right-handed trefoil $T$, which is the braid closure of $\sigma^3$. Evaluating $\Omega^q$ at $d = -x^{-1}$, we will obtain the Jones polynomial of the left-handed trefoil. We then compute $\qHFB(\sigma^3)$ and the spectral sequence to $\HFKh(T)$. Finally, we give an example that shows that $\Omega^q$ is not invariant under Markov conjugation or stabilisation.

\subsection{The right-handed trefoil}
\begin{figure}
\centering
\def\svgwidth{1.3\textwidth}
\begingroup%
  \makeatletter%
  \providecommand\color[2][]{%
    \errmessage{(Inkscape) Color is used for the text in Inkscape, but the package 'color.sty' is not loaded}%
    \renewcommand\color[2][]{}%
  }%
  \providecommand\transparent[1]{%
    \errmessage{(Inkscape) Transparency is used (non-zero) for the text in Inkscape, but the package 'transparent.sty' is not loaded}%
    \renewcommand\transparent[1]{}%
  }%
  \providecommand\rotatebox[2]{#2}%
  \newcommand*\fsize{\dimexpr\f@size pt\relax}%
  \newcommand*\lineheight[1]{\fontsize{\fsize}{#1\fsize}\selectfont}%
  \ifx\svgwidth\undefined%
    \setlength{\unitlength}{768bp}%
    \ifx\svgscale\undefined%
      \relax%
    \else%
      \setlength{\unitlength}{\unitlength * \real{\svgscale}}%
    \fi%
  \else%
    \setlength{\unitlength}{\svgwidth}%
  \fi%
  \global\let\svgwidth\undefined%
  \global\let\svgscale\undefined%
  \makeatother%
  \begin{picture}(1,0.2265625)%
    \lineheight{1}%
    \setlength\tabcolsep{0pt}%
    \put(0,0){\includegraphics[width=\unitlength,page=1]{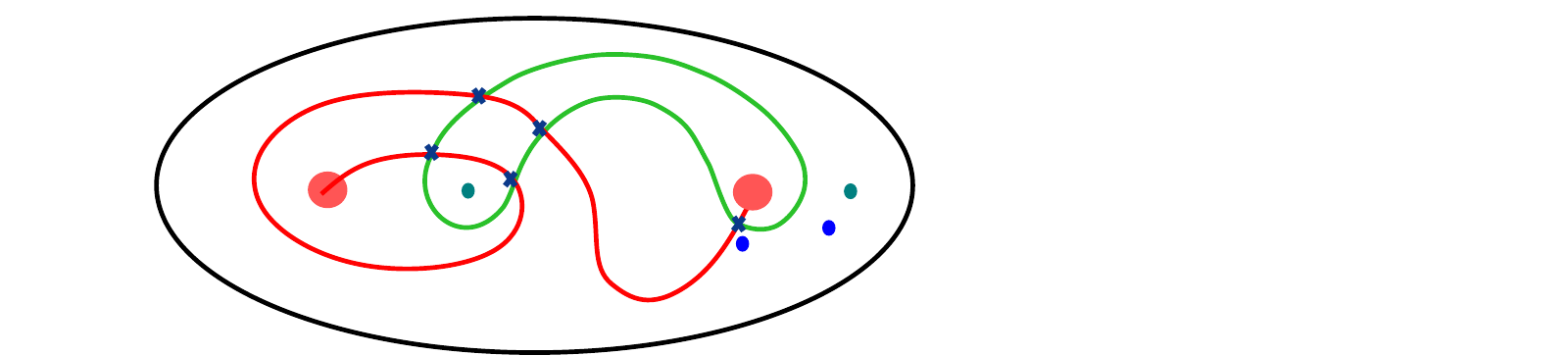}}%
    \put(0.22910945,0.1081232){\color[rgb]{0,0,0}\makebox(0,0)[lt]{\lineheight{1.25}\smash{\begin{tabular}[t]{l}\tiny $x^{-1}d^3$\end{tabular}}}}%
    \put(0.3346068,0.10926004){\color[rgb]{0,0,0}\makebox(0,0)[lt]{\lineheight{1.25}\smash{\begin{tabular}[t]{l}\tiny $-d^3$\end{tabular}}}}%
    \put(0.28744029,0.17856839){\color[rgb]{0,0,0}\makebox(0,0)[lt]{\lineheight{1.25}\smash{\begin{tabular}[t]{l}\tiny $xd^3$\end{tabular}}}}%
    \put(0.35495643,0.13870205){\color[rgb]{0,0,0}\makebox(0,0)[lt]{\lineheight{1.25}\smash{\begin{tabular}[t]{l}\tiny $-x^2d^3$\end{tabular}}}}%
    \put(0.43060006,0.0804815){\color[rgb]{0,0,0}\makebox(0,0)[lt]{\lineheight{1.25}\smash{\begin{tabular}[t]{l}\tiny $x^2d^4$\end{tabular}}}}%
  \end{picture}%
\endgroup%

\caption{The quantum Heegaard diagram associated with the braid $\sigma^3$ and the terms of the quantum intersection polynomial $\Omega^q_{\sigma^3}$.}\label{fig:trefoil}
\end{figure}

The disc model $\bD_{\sigma^3}^q$ consists of $D^2$ with punctures $z_0$, $z_1$, $w_0$, $w_1$, $q_0$, and $q_1$. There is one arc $a_1$ and one curve $b_1$. After applying the action of $\sigma^3 \cup \bI_4$ to $a_1$ and closing it up along the tube added at its endpoints, we obtain the curve $\alpha_1$ shown in Figure~\ref{fig:trefoil}.

Note that $n = 1$ and $w(\sigma^3)=3$; hence, $(-1)^\frac{n(n+1)}{2} = -1$. The intersection $\alpha_1 \cap \beta_1$ consists of five points $\xs$, labelled in Figure~\ref{fig:trefoil} by $-\varepsilon_{\xs} \, x^{A^{\HF}(\xs)} d^{A^{\qHF}(\xs)}$. This gives us the quantum intersection polynomial
\begin{equation}
\begin{split}
\Omega^q_{\sigma^{3}}(x,d) &= x^{-1}d^3 - d^3 + xd^3 - x^2d^3 + x^2d^4 \\
&= d^3 \left(x^{-1} - 1 + x -x^2 + x^2d\right).
\end{split}
\end{equation}

\begin{figure}
\centering
\def\svgwidth{\textwidth}
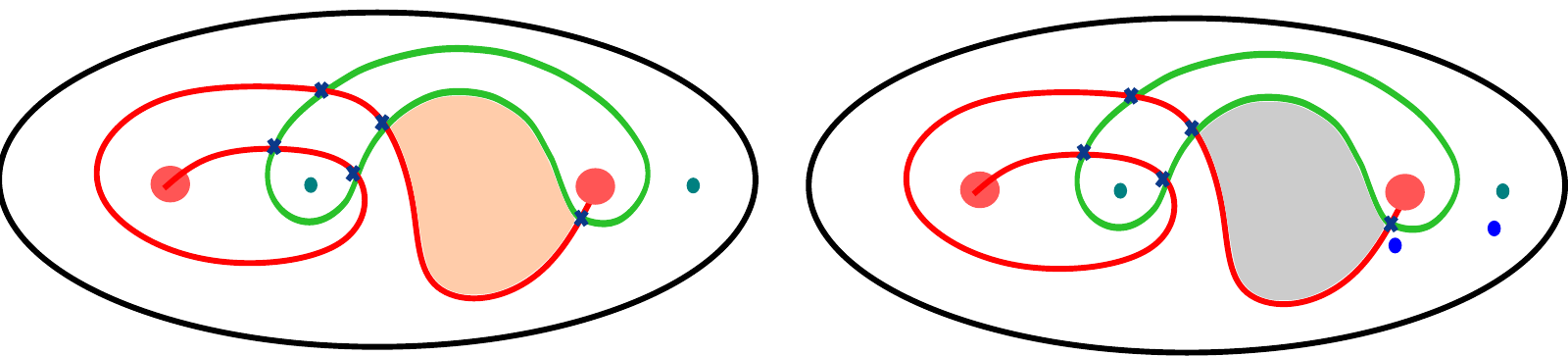
\caption{Left: The terms of the Alexander polynomial obtained after the substitution $d=1$. The only disc avoiding the base points $w$ and $z$ and hence contributing to the knot Floer boundary map is shaded. Right: The terms of the Jones polynomial obtained after the substitution $d=-x^{-1}$.}
\label{fig:trefoil-Alexander-Jones}
\end{figure}

If we substitute $d=1$, we get
\begin{equation}
\Omega^q_{\sigma^3}(x,1) = x^{-1}-1+x,
\end{equation}
which is indeed the Alexander polynomial $\Delta_T(x)$ of the right-handed trefoil; see the left of Figure~\ref{fig:trefoil-Alexander-Jones}. If we substitute $d=-x^{-1}$, we obtain
\begin{equation}
\Omega^q_{\sigma^3}(x,-x^{-1})= -x^{-4}+x^{-3}+x^{-1},
\end{equation}
which is indeed the Jones polynomial $V_{\overline{T}}(x)$ of the left-handed trefoil knot $\overline{T}$ when we use the conventions of Remark~\ref{rmk:Jones-convention}; see the right-hand side of Figure~\ref{fig:trefoil-Alexander-Jones}.

We now consider the knot Floer homology that arises from the Heegaard surface on the left of Figure~\ref{fig:trefoil-Alexander-Jones}. We have five intersection points. The shaded disc $u_1$ connects the intersection points labelled $-x^2$ to $x^2$, has $n_z(u_1) = n_w(u_1) = 0$, and a single rigid pseudo-holomorphic representative. There are no other discs that contribute to the boundary map: There is a disc $u_2$ from $-1$ to $x^{-1}$ and $u_3$ from $-x^2$ to $x$, but $n_z(u_2) = n_z(u_3) = 1$. So $x^2$ and $-x^2$ cancel in the homology of the Floer complex, and the other three generators survive. We obtain
\begin{equation}\label{eqn:HFK}
  \HFKh_{a, m}(T)=
  \begin{cases}
        \bF_2 \text{ for } (a, m) \in \{(-1,-2),(0,-1),(1,0)\};\\
        0 \ \text{otherwise}.
    \end{cases}
\end{equation}
This is the knot Floer homology of the right-handed trefoil.

Finally, we compute $\qHFB(\sigma^3)$. Consider the left-hand side of Figure~\ref{fig:trefoil-Alexander-Jones}. As $W(-x^2) = 0$ and $W(x^2) = 1$, the only rigid pseudo-holomorphic disc $u_1$ in the knot Floer complex does not contribute to $\partial^q$. Hence, $\partial^q = 0$ and 
\[
\dim\left(\qHFB(\sigma^3)\right) = 5.
\]
We have already computed the gradings $A^{\HF}$ and $A^{\qHF}$ of all five intersection points; see Figure~\ref{fig:trefoil}. For the Maslov gradings of the intersection points labelled $-1$, $x^{-1}$, and $x$ on the left-hand side of Figure~\ref{fig:trefoil-Alexander-Jones}, see equation~\eqref{eqn:HFK}. Since $\mu(u_3) = 1$ and $M(x) = 0$, we have $M(-x^2) = 1$. Using that $\mu(u_1) = 1$, we then obtain that $M(x^2) = 0$. In summary,
\[
  \qHFB_{a,q,m}(\sigma^3)=
  \begin{cases}
        \bF_2 \text{ for } (a, q, m) \in \{(-1,3,-2), (0,3,-1), (1,3,0), (2,3,1),  (2,4,0)\};\\
        0 \ \text{otherwise}.
    \end{cases}
\]
The graded Euler characteristic is indeed $\Omega^q_{\sigma^3}$. Now, consider the spectral sequence in Theorem~\ref{thm:spectral-sequence}. This is induced by the filtration $G = A^{\qHF} - 2W$. As 
\[
A^{\qHF}(-x^2) - 2W(-x^2) = 3 - 2 \cdot 0 = 3, \quad A^{\qHF}(x^2) - 2W(x^2) = 4 - 2 \cdot 1 = 2, 
\]
the disc $u_1$ decreases the filtration $G$ by one. Hence, the $E^1$ page does not count $u_1$ and is $\qHFB(\sigma^3)$, while the $E^2$ page counts $u_1$ and is $\HFKh(T)$. The element $\Theta(\beta) = [\xs^0]$ is represented by the point labelled $x^2$ of grading $(2,4,0)$. It is non-zero in $\qHFB(\sigma^3)$ but vanishes in $\HFKh(T)$. 

\subsection{Non-invariance under Markov moves}
We illustrate that $\Omega^q$, and hence $\qHFB$, are not invariant under Markov stabilisation or conjugation, based on the work of the first author~\cite{Anghel2024AIF}. By Theorem~\ref{qINT}, we have
\[
\Omega^q_\beta(x,d) = \Omega'_\beta(x,-d),
\]
so it suffices to show this for $\Omega'$. By a simple computation,
\[
\Omega'_{\sigma_1}(x,d) = xd^2 + xd - d.
\]
The positive Markov stabilisation $\sigma_1\sigma_2$ of $\sigma_1$ satisfies
\[
\Omega'_{\sigma_1\sigma_2}(x,d) = x^2d^4 + x^2d^3 - xd^3 + xd - d.
\]
Furthermore, the conjugate $\sigma_1^2\sigma_2\sigma_1^{-1}$ has
\[
\Omega'_{\sigma_1^2\sigma_2\sigma_1^{-1}}(x,d) = x^2d^2(d+1-x^{-1})^2.
\]
By \cite[Theorem~1.6]{Anghel2024AIF}, we obtain a knot invariant when we pass to the quotient 
\[
\bZ[x^{\pm 1}, d^{\pm 1}]/((d+1)(xd-1)). 
\]
Indeed, the differences
\[
\begin{split}
\Omega'_{\sigma_1\sigma_2}(x,d) - \Omega'_{\sigma_1}(x,d) &= xd^2(d+1)(xd-1); \\
\Omega'_{\sigma_1^2\sigma_2\sigma_1^{-1}}(x,d) - \Omega'_{\sigma_1\sigma_2}(x,d) &= d(x-1)(d+1)(xd-1)
\end{split}
\]
vanish when we set $(d+1)(xd-1) = 0$.

\bibliographystyle{plain}
\bibliography{references}

\end{document}